\documentclass[pmlr,11pt]{jmlr}
\usepackage{booktabs}
\usepackage{multirow}
\usepackage{paralist,amsmath, amssymb, bm}
\usepackage{algorithm,algorithmic}
\usepackage{graphicx}
\usepackage{subfigure}
\usepackage{times}
\usepackage{afterpage}

\usepackage{caption}
\makeatletter
\newcounter{ALC@tempcntr}
\makeatother

\makeatletter
\AtBeginDocument{\Hy@breaklinkstrue}
\makeatother

\makeatletter 
  \newcommand\figcaption{\def\@captype{figure}\caption} 
  \newcommand\tabcaption{\def\@captype{table}\caption} 
\makeatother

\newtheorem{thm}{Theorem}

\newtheorem{lem}{Lemma}
\newtheorem{cor}[thm]{Corollary}
\newtheorem{ass}{Assumption}

\newtheorem{remark}{Remark}

\def \y {\mathbf{y}}
\def \E {\mathrm{E}}
\def \x {\mathbf{x}}
\def \bv {\mathbf{v}}
\def \g {\mathbf{g}}

\def \O {\widetilde{O}}

\def \bo {\mathbf{o}}

\def \z {\mathbf{z}}
\def \Z {\mathcal{Z}}

\def \u {\mathbf{u}}

\def \w {\mathbf{w}}
\def \bv {\mathbf{v}}
\def \R {\mathbb{R}}
\def \zn {\mathbb{Z}_+}

\def \P {\mathcal{P}}

\def \W {\mathcal{W}}

\def \q {\mathbf{q}}
\def \bv {\mathbf{v}}

\def \q {\mathbf{q}}

\def \SS {\mathcal{S}}

\def \zh {\hat{\z}}

\def \wb {\bar{\w}}
\def \qb {\bar{\q}}

\def \wh {\widehat{\w}}
\def \qh {\widehat{\q}}

\def \bv {\mathbf{v}}

\def \Lt  {\widetilde{L}}

\DeclareMathOperator*{\argmin}{argmin}
\DeclareMathOperator*{\argmax}{argmax}
\DeclareMathOperator*{\sgn}{sign}

\begin{document}

\title[Minimax Excess Risk Optimization]{Efficient Stochastic Approximation of \\ Minimax Excess Risk Optimization}

\coltauthor{%
 \Name{Lijun Zhang} \Email{zhanglj@lamda.nju.edu.cn}\\
 \Name{Haomin Bai} \Email{baihm@lamda.nju.edu.cn}\\
 \addr National Key Laboratory for Novel Software Technology, Nanjing University, Nanjing, China
 \AND
  \Name{Wei-Wei Tu} \Email{tuwwcn@gmail.com}\\
 \addr Artificial Productivity Inc., Beijing, China
   \AND
  \Name{Ping Yang} \Email{jiadi@xiaohongshu.com}\\
    \Name{Yao Hu} \Email{yaoohu@gmail.com}\\
 \addr Xiaohongshu Inc., Beijing, China
}

\maketitle

\begin{abstract}
While traditional distributionally robust optimization (DRO) aims to minimize the maximal risk over a set of distributions,  \citet{Regret:RML:DS} recently proposed a variant that replaces risk with \emph{excess risk}. Compared to DRO, the new formulation---minimax excess risk optimization (MERO) has the advantage of suppressing the effect of heterogeneous noise in different distributions. However, the choice of excess risk leads to a very challenging minimax optimization problem, and currently there exists only an inefficient algorithm for empirical MERO. In this paper, we develop efficient stochastic approximation approaches which directly target MERO. Specifically, we leverage techniques from stochastic convex optimization to estimate the minimal risk of every distribution, and solve MERO as a stochastic convex-concave optimization (SCCO) problem with biased gradients. The presence of bias makes existing theoretical guarantees of SCCO inapplicable, and fortunately, we demonstrate that the bias, caused by the estimation error of the minimal risk, is under-control. Thus, MERO can still be optimized with a nearly optimal convergence rate. Moreover, we investigate a practical scenario where the quantity of samples drawn from each distribution may differ, and propose a stochastic approach that delivers \emph{distribution-dependent} convergence rates.
\end{abstract}

\section{Introduction}
With the widespread application of machine learning, it is common to encounter situations where the test distribution differs from the training distribution \citep{JMLR:v8:sugiyama07a,pmlr-v28-zhang13d,10.1145/2523813}.  When faced with distribution shifts, learning models trained by conventional techniques  (e.g., empirical risk minimization) tend to suffer significant performance degradation \citep{pmlr-v139-koh21a}. Distributionally robust optimization (DRO) seeks a model with minimax risk over a set of potential distributions, offering a principled approach for generalization across distributions \citep{10.1214/20-AOS2004}. In general, DRO can be formulated as the following minimax optimization problem
\begin{equation} \label{eqn:dro}
\min_{\w \in \W}     \sup_{\P \in \SS} \ \left\{\E_{\z \sim \P}\big[\ell(\w;\z)\big]\right\}
\end{equation}
where $\SS$ is a set of distributions,  $\z \in \Z$ denotes a random sample drawn from $\P$, $\W$ is a hypothesis class, and $\ell(
\cdot;\cdot)$ denotes a loss function that measures the performance. When $\SS$ contains a finite number of distributions, (\ref{eqn:dro}) is referred to as Group DRO  (GDRO) \citep{Gouop_DRO}. One motivating example involves deploying a common classifier across multiple hospitals, each serving populations with varied demographics, to predict disease occurrences \citep{NIPS2017_186a157b}.

Although the minimax formulation enhances the model's robustness to distribution shifts,  it carries the potential of being sensitive to heterogeneous noise, especially in GDRO where the candidate distributions might exhibit substantial differences. The rationale is evident: a single distribution with high levels of noise could easily dominate the maximum in (\ref{eqn:dro}), causing other distributions to be essentially overlooked. To overcome this limitation, \citet{Regret:RML:DS} recently proposed a variant of DRO, namely minimax regret optimization (MRO), that replaces the raw risk $\E_{\z \sim \P}[\ell(\w;\z)]$ in (\ref{eqn:dro}) with \emph{excess risk} (or regret):
\[
\E_{\z \sim \P}\big[\ell(\w;\z)\big] - \min_{\w \in \W}  \E_{\z \sim \P}\big[\ell(\w;\z)\big] .
\]
 In particular, they consider the setting of GDRO, and seek to minimize the worst-case excess risk over $m$ distributions $\P_1,\ldots,\P_m$:
\begin{equation} \label{eqn:mro}
\min_{\w \in \W}   \max_{i\in[m]}  \  \Big \{\smash[b]{\underbrace{\E_{\z \sim \P_i}\big[\ell(\w;\z)\big]}_{:=R_i(\w)}} - \smash[b]{\underbrace{\min_{\w \in \W}  \E_{\z \sim \P_i}\big[\ell(\w;\z)\big]}_{:=R_i^*}} \vphantom{\underbrace{\min_{\w \in \W}  \E_{\z \sim \P_i}\big[\ell(\w;\z)\big]}_{:=R_i^*}}\Big\}.
\end{equation}
To avoid confusion with the term ``regret'' commonly used in online learning \citep{bianchi-2006-prediction}, in this paper, we refer to MRO as minimax excess risk optimization (MERO).

MERO can be understood  as subtracting the intrinsic difficulty of each distribution (i.e., the minimal risk $R_i^*$) from the risk, making the remaining quantities more comparable.  However, since the value of $R_i^*$ is typically unknown, the optimization problem in (\ref{eqn:mro}) poses a significant challenge. In the previous work, \citet{Regret:RML:DS}  choose to optimize the empirical counterpart of MERO, but the resulting algorithm requires solving an empirical risk  minimization problem in each iteration, which is computationally prohibitive when the number of samples is large. Along this line of research, this paper aims to develop efficient stochastic approximation approaches that optimize (\ref{eqn:mro}) \emph{directly}. It is worth mentioning that by doing so, we also eliminate the need to analyze the discrepancy between the empirical MERO and the population MERO. To exploit techniques from stochastic approximation \citep{nemirovski-2008-robust}, we cast (\ref{eqn:mro}) as a stochastic convex-concave saddle-point problem:
\begin{equation} \label{eqn:convex:concave}
\min_{\w \in \W} \max_{\q \in \Delta_m}  \  \left\{\phi(\w,\q)= \sum_{i=1}^m q_i \big [R_i(\w)- R_i^*\big] \right\}
\end{equation}
where $\Delta_m=\{\q \in \R^m: \q \geq 0, \sum_{i=1}^m q_i=1\}$ is the ($m{-}1$)-dimensional simplex.

If the value of $R_i^*$ is known, we can readily employ stochastic mirror descent (SMD) to optimize (\ref{eqn:convex:concave}), mirroring the application of SMD to GDRO  \citep[Algorithm 1]{SA:GDRO}. An intuitive approach thus emerges, which first estimates the value of $R_i^*$, denoted by $\widehat{R}_i$, replaces the function $\phi(\w,\q)$ in (\ref{eqn:convex:concave}) with $\sum_{i=1}^m q_i  [R_i(\w)- \widehat{R}_i]$, and  subsequently applies SMD to the revised problem.  As elaborated later, this multi-stage approach is indeed feasible, but it needs to fix the total number of iterations (or the total number of random samples), and cannot return a solution upon demand. To possess the \emph{anytime} ability, i.e., the capability to deliver a solution at any  round, we propose a stochastic approximation approach that alternates between estimating $R_i^*$ and optimizing (\ref{eqn:convex:concave}) in an iterative manner. Specifically, we execute $m$ instances of SMD, each designed to minimize the risk of an individual distribution.  In this way, at each round $t$, we have a model $\wb_t^{i}$ that attains a low risk for each distribution $\P_i$. Leveraging these $\wb_t^{i}$s, we are able to construct biased stochastic gradients for (\ref{eqn:convex:concave}), and proceed to optimize  the problem by SMD. Our theoretical analysis demonstrates that the bias is under-control, and the optimization error reduces at an $\O(\sqrt{(\log m)/t})$ rate,\footnote{We use the $\O$ notation to hide constant factors as well as polylogarithmic factors in $t$.} which is optimal up to a polylogarithmic factor. 

Furthermore, we investigate MERO under the imbalanced setting where the number of samples drawn from each distribution is different. To this end, we develop a stochastic approximation approach for a weighted formulation of MERO, which yields a convergence rate for each distribution depending on the budget of that particular distribution. Under appropriate conditions, we demonstrate that for smooth risk functions, the excess risk of the $i$-th distribution decreases at an $O((\log m)/\sqrt{n_i})$ rate, where $n_i$ is the sample budget. Finally, we conduct experiments to validate the efficiency and effectiveness of our algorithms.

\section{Related Work}
There is an extensive body of literature on DRO, with various problem formulations and research focuses \citep{OJMO:DRO:Review}. Below, we provide a brief review of related work, with a specific emphasis on studies examining GDRO and MERO.

For the DRO problem in (\ref{eqn:dro}), the uncertain set $\SS$ is typically chosen to be the neighborhood surrounding a target distribution $\P_0$, and constructed by using certain distance functions between distributions, such as the $f$-divergence \citep{DRO:MS:2013}, the Wasserstein distance \citep{DRO:MP:2018,doi:10.1287/educ.2019.0198}, and the maximum mean discrepancy \citep{NEURIPS2019_1770ae9e}. Other ways for defining $\SS$ include moment constraints \citep{10.2307/40792682,opre.2014.1314} and hypothesis testing of goodness-of-fit  \citep{RSAA:MP:2018}.  The nature of the loss function can exhibit variability: it could assume a convex form \citep{doi:10.1287/opre.2015.1374,DRSP:SIAM:2017}, a non-convex structure \citep{NEURIPS2021_164bf317,NEURIPS2021_533fa796}, or potentially incorporate a regularizer \citep{Sinha:ICLR:2018}. Research efforts may be directed towards diverse objectives, including the development of the optimization algorithms \citep{NIPS2016_4588e674,NEURIPS2020_64986d86,Tianbao:22,DRO:Composite,NEURIPS2022_8a54a80f}, the exploration of finite sample and asymptotic properties of the empirical solution \citep{NIPS2017_5a142a55,10.1214/20-AOS2004}, the determination of confidence intervals for the risk \citep{Statistics:RO}, or the approximation of nonparametric likelihood \citep{NEURIPS2019_b49b0133}.

When the number of distributions in $\SS$ is finite,  (\ref{eqn:dro}) becomes GDRO  \citep{Gouop_DRO}, which can be formulated as the following stochastic convex-concave problem:
\begin{equation} \label{eqn:GDRO:convex:concave}
\min_{\w \in \W} \max_{\q \in \Delta_m}  \   \left\{ \sum_{i=1}^m q_i R_i(\w) \right\}.
\end{equation}
By drawing $1$ random sample in each round, \citet{Gouop_DRO} have applied SMD to (\ref{eqn:GDRO:convex:concave}), and established an $O(m \sqrt{(\log m)/T})$ convergence rate. As a result, their algorithm achieves an $O(m^2 (\log m)/\epsilon^2)$ sample complexity for finding an $\epsilon$-optimal solution, which is suboptimal according to the $\Omega(m/\epsilon^2)$ lower bound \citep[Theorem 5]{DRO:Online:Game}.  In the literature, there exist $3$ different ways to obtain a (nearly) optimal $O(m (\log m)/\epsilon^2)$ sample complexity.
\begin{compactenum}
\item In each round, we draw $m$ random samples, one from each distribution, and then update $\w$ and $\q$ in (\ref{eqn:GDRO:convex:concave}) by SMD. In this way, the optimization error reduces at an $O(\sqrt{(\log m)/T})$ rate \citep[\S 3.2]{nemirovski-2008-robust}, implying an $O(m (\log m)/\epsilon^2)$ sample complexity.
\item In each round, we only draw $1$ sample from one specific distribution, update $\w$ by SMD, and update $\q$ by an online algorithm for non-oblivious multi-armed bandits \citep{2020:bandit-book}. It can be shown that the optimization error reduces at an $O(\sqrt{m (\log m)/T})$ rate \citep{DRO:Online:Game,SA:GDRO}, keeping the same sample complexity.
\item The third approach is similar to the second one, but it updates $\q$ by integrating  SMD with gradient clipping \citep[\S A.1]{NEURIPS2022_e90b00ad}.
\end{compactenum}

Recently, \citet{SA:GDRO} have considered a practical scenario where the number of samples that can be drawn from each distribution is different. To this end, they introduce a weighted formulation of GDRO:
\begin{equation} \label{eqn:convex:concave:weight:GDRO}
\min_{\w \in \W} \max_{\q \in \Delta_m}  \   \left\{\sum_{i=1}^m q_i p_i R_i(\w) \right\}
\end{equation}
where $p_i>0$ is a weight assigned to distribution $\P_i$. Then, they develop two stochastic algorithms to solve (\ref{eqn:convex:concave:weight:GDRO}), and establish distribution-wise convergence rates. In this way, the  rate is not dominated by the distribution with the smallest budget. Inspired by \citet{SA:GDRO}, we will investigate MERO with varying sample sizes across distributions in Section~\ref{sec:W:MERO}. The technique of introducing scale factors has been previously presented in the study of MERO under heterogeneous distributions \citep[\S 5]{Regret:RML:DS}, where all distributions have the same number of samples, but with different complexities.

To optimize MERO,  \citet{Regret:RML:DS} consider the empirical version: 
\begin{equation} \label{eqn:mro:empirical}
\min_{\w \in \W}   \max_{i\in[m]}  \  \Bigg \{\frac{1}{n} \sum_{j=1}^n \ell(\w;\z^{(i,j)}) - \min_{\w \in \W}  \frac{1}{n} \sum_{j=1}^n \ell(\w;\z^{(i,j)}) \Bigg\}
\end{equation}
where $\{\z^{(i,j)}: j=1,\ldots,n\}$ are random samples drawn from distribution $\P_i$, and analyze the generalization performance using classical tools from learning theory. They have developed an iterative method for solving (\ref{eqn:mro:empirical}), which needs to address an empirical risk minimization problem in each iteration, rendering the process inefficient. It's worth noting that the principle of minimizing the worst-case excess risk has surfaced in various other fields \citep{Minimax:Parameter,JMLR:v8:alaiz-rodriguez07a,Minimax:Commitment,Minimax:Attacks}. 

\section{Stochastic Approximation of MERO}
We first describe preliminaries of stochastic approximation, including the setup and assumptions, then develop a multi-stage approach, and finally propose an anytime approach. 
\subsection{Preliminaries}\label{sec:pre}
We first present the standard setup of mirror descent \citep{nemirovski-2008-robust}.  We endow the domain $\W$ with a distance-generating function $\nu_w(\cdot)$, which is $1$-strongly convex w.r.t.~a specific norm $\|\cdot\|_w$.  We define the Bregman distance corresponding to  $\nu_w(\cdot)$ as
\[
B_w(\u,\bv)= \nu_w(\u) - \big[\nu_w(\bv) + \langle \nabla \nu_w(\bv) , \u-\bv\rangle\big].
\]
For the simplex $\Delta_m$, we select the entropy function $\nu_q(\q) = \sum_{i=1}^m q_i \ln q_i$, which demonstrates $1$-strong convexity w.r.t.~the vector $\ell_1$-norm $\|\cdot\|_1$, as the distance-generating function. In a similar manner, $B_q(\cdot,\cdot)$ represents the Bregman distance associated with $\nu_q(\cdot)$, which is the Kullback--Leibler divergence between probability distributions.

Next, we introduce standard assumptions used in our analysis.
\begin{ass}\label{ass:1}
All the risk functions $R_1(\cdot),\cdots,R_m(\cdot)$ and the domain $\W$ are convex.
\end{ass}
\begin{ass}\label{ass:2}
The domain $\W$ is bounded in the sense that
\begin{equation} \label{eqn:domain:W}
 \max_{\w \in \W}  B_w(\w,\bo_w) \leq D^2.
\end{equation}
where $\bo_w =\argmin_{\w \in \W} \nu_w(\w)$.
\end{ass}
For the simplex $\Delta_m$, we have $\max_{\q \in \Delta_m}  B_q(\q,\bo_q) \leq \ln m$ where $\bo_q= \frac{1}{m} \mathbf{1}_m \in \R^m$ and $\mathbf{1}_m$ is the $m$-dimensional vector of all ones \citep[Proposition 5.1]{Beck_Marc}.

We assume that the gradient is bounded, and the loss  belongs to $[0,1]$.
\begin{ass}\label{ass:3}
For all $i \in [m]$, we have
\begin{equation} \label{eqn:gradient}
\|\nabla \ell(\w;\z) \|_{w,*}\leq G, \ \forall \w \in \W,  \ \z \sim \P_i
\end{equation}
where $\|\cdot\|_{w,*}$ is the dual norm of $\|\cdot\|_w$.
\end{ass}
\begin{ass}\label{ass:4}
For all $i \in [m]$, we have
\begin{equation} \label{eqn:value}
0 \leq \ell(\w;\z) \leq 1, \ \forall \w \in \W,  \ \z \sim \P_i.
\end{equation}
\end{ass}

Given an solution $(\wb, \qb)$ to (\ref{eqn:convex:concave}), the optimization error is defined as
\begin{equation} \label{eqn:per:measure}
\epsilon_{\phi}(\wb, \qb) = \max_{\q\in \Delta_m}  \phi(\wb,\q)- \min_{\w\in \W}  \phi(\w,\qb),
\end{equation}
which provides an upper bound for the excess risk on any distribution, since
\begin{equation} \label{eqn:relation:error}
\begin{split}
&\max_{i\in[m]}  \big[R_i(\wb)- R_i^*\big] -  \min_{\w \in \W}   \max_{i\in[m]}   \big[R_i(\w)- R_i^*\big]\\
 \leq & \max_{\q \in \Delta_m} \sum_{i=1}^m q_i \big[R_i(\wb) - R_i^*\big]- \min_{\w\in \W} \sum_{i=1}^m \bar{q}_i \big[R_i(\w)- R_i^*\big] =  \epsilon_{\phi}(\wb, \qb).
\end{split}
\end{equation}
\subsection{A Multi-Stage Stochastic Approximation Approach for MERO} \label{sec:multi:stage}
As mentioned in the introduction, we can design a multi-stage stochastic approach for MERO. Notably, analogous methodologies have found their application in the empirical research conducted on language modeling \citep{oren-etal-2019-distributionally,DRO:DoReMi}.

\paragraph{Stage 1: Minimizing the risk} For each distribution $\P_i$, we run an instance of SMD to minimize the risk $R_i(\cdot)$, and obtain an approximate solution $\wb^{(i)}$. We execute each SMD  for $T$ iterations, and thus consume $m T$ samples. From the theoretical guarantee of SMD \citep{nemirovski-2008-robust}, with probability at least $1-\delta$, we have $R_i(\wb^{(i)}) - R_i^* = O( \sqrt{\log(1/\delta)/T} )$, for each $i \in [m]$. By the union bound, with high probability, we have $\max_{i \in [m]} [R_i(\wb^{(i)}) - R_i^*] = O( \sqrt{(\log m)/T} )$.

\paragraph{Stage 2: Estimating the minimal risk} To estimate the value of $R_i(\wb^{(i)})$, we draw  $T$ samples $\z_1^{(i)}, \ldots, \z_T^{(i)}$ from each distribution $\P_i$, and calculate $\widehat{R}_i(\wb^{(i)}) = \frac{1}{T}\sum_{j=1}^{T} \ell(\wb^{(i)};\z_j^{(i)})$. From standard concentration inequalities \citep{Concentration_inequalities} and the union bound, with high probability, we have $ \max_{i \in [m]}|\widehat{R}_i(\wb^{(i)})- R_i(\wb^{(i)})| = O( \sqrt{(\log m)/T} )$. In this step, we also use $m T$ samples.

\paragraph{Stage 3: Applying SMD to an approximate problem} From the above two steps, with high probability, we have
\begin{equation} \label{eqn:appror:risk}
 \max_{i \in [m]} |\widehat{R}_i(\wb^{(i)})- R_i^*| = O( \sqrt{(\log m)/T} ).
\end{equation}
Then, we formulate  the following problem
\begin{equation} \label{eqn:convex:concave:appror}
  \min_{\w \in \W} \max_{\q \in \Delta_m}  \  \left\{\hat{\phi}(\w,\q)= \sum_{i=1}^m q_i \big [R_i(\w)- \widehat{R}_i(\wb^{(i)})\big] \right\}
\end{equation}
which serves as an approximation to (\ref{eqn:convex:concave}). Since $\widehat{R}_i(\wb^{(i)})$ is a constant, we can directly apply SMD to (\ref{eqn:convex:concave:appror}). After $T$ iterations, we obtain solutions $\wb$ and $\qb$ such that, with high probability $\max_{\q\in \Delta_m}  \hat{\phi}(\wb,\q)- \min_{\w\in \W}  \hat{\phi}(\w,\qb)= O( \sqrt{(\log m)/T})$ \citep{SA:GDRO}. From (\ref{eqn:appror:risk}), it is easy to prove that $\max_{\q\in \Delta_m}  \phi(\wb,\q)- \min_{\w\in \W}  \phi(\w,\qb) = O( \sqrt{(\log m)/T} )$. This step also requires $m T$  samples.

In summary, the above approach reduces the optimization error to $O( \sqrt{(\log m)/T} )$, at the cost of $3mT$ samples. In other words, it attains an $O(m (\log m)/\epsilon^2)$ sample complexity, which is nearly optimal according to the  lower bound of GDRO \citep{DRO:Online:Game}. We would like to highlight that the 2nd stage is not essential, and is included to facilitate understanding. In fact, we can omit the 2nd stage, and define $\hat{\phi}(\w,\q)= \sum_{i=1}^m q_i  [R_i(\w)- R_i(\wb^{(i)})]$ in (\ref{eqn:convex:concave:appror}). The resulting optimization problem can still be solved by SMD, and the sample complexity remains in the same order. An illustrative example of this two-stage approach can be found in Section~\ref{sec:W:MERO}.

\begin{remark}
\textnormal{While the multi-stage approach achieves a nearly optimal sample complexity, it suffers two \emph{limitations}: (i) the total number of iterations must be predetermined; (ii) a solution is available only when the algorithm enters the final stage. To circumvent these drawbacks, we put forth a stochastic approximation approach that interleaves the aforementioned three stages together, being able to return a solution at any time.}
\end{remark}

\subsection{An Anytime Stochastic Approximation Approach for MERO}
We maintain $m$ instances of SMD to minimize all the risk functions $R_1(\cdot),\cdots,R_m(\cdot)$, and meanwhile utilize their solutions to optimize (\ref{eqn:convex:concave}) according to SMD.

For the purpose of minimizing $R_i(\cdot)$, we denote by  $\w_t^{(i)}$ the solution in the $t$-th iteration. We first draw $1$ sample $\z_t^{(i)}$ from each distribution $\P_i$, and calculate the stochastic gradient $\nabla \ell(\w_t^{(i)};\z_t^{(i)})$ which is an unbiased estimator of $\nabla R_i(\w_t^{(i)})$. According to SMD \citep{nemirovski-2008-robust}, we update $\w_t^{(i)}$ by
\begin{equation} \label{eqn:wti:smd}
\w_{t+1}^{(i)}= \argmin_{\w \in \W} \left\{ \eta_t^{(i)} \langle \nabla \ell(\w_t^{(i)};\z_t^{(i)}) , \w -\w_t^{(i)}\rangle + B_w(\w,\w_t^{(i)}) \right\} , \ \forall i \in [m]
\end{equation}
where $ \eta_t^{(i)} >0$ is the step size. Due to technical reasons, we will use the weighted average of iterates
\begin{equation} \label{eqn:wti:average}
\wb_t^{(i)}= \sum_{j=1}^t \frac{\eta_j^{(i)} \w_j^{(i)} }{\sum_{k=1}^t \eta_k^{(i)}} = \frac{(\sum_{j=1}^{t-1} \eta_j^{(i)})\wb_{t-1}^{(i)}  + \eta_t^{(i)} \w_t^{(i)} }{\sum_{k=1}^t \eta_k^{(i)}}
\end{equation}
as an approximate solution to $\min_{\w \in \W} R_i(\w)$. While selecting the last iterate $\w_t^{(i)}$ is also an option  \citep{ICML2013Shamir:ICML,pmlr-v99-harvey19a,pmlr-v99-jain19a}, this choice leads to a more complex analysis. Therefore, we prefer to employ $\wb_t^{(i)}$.

We proceed to minimize (\ref{eqn:convex:concave}) by SMD. Let $\w_t$ and $\q_t$ be the solutions in the $t$-th round. Based on the random samples $\z_t^{(1)},\ldots,\z_t^{(m)}$, we define the stochastic gradient of $\phi(\cdot,\cdot)$ at $(\w_t,\q_t)$ w.r.t.~$\w$ as
\begin{equation}\label{eqn:stoch:grad:1}
\g_{w}(\w_t,\q_t) = \sum_{i=1}^m q_{t,i} \nabla \ell(\w_{t};\z_t^{(i)}),
\end{equation}
which is an unbiased estimator of the true gradient $\nabla_\w \phi(\w_t,\q_t) = \sum_{i=1}^m q_{t,i} \nabla R_i(\w_t)$. The challenge lies in the construction of the stochastic gradient w.r.t.~$\q$. To this end, we define
\begin{equation}\label{eqn:stoch:grad:2}
\g_{q}(\w_t,\q_t)= \left[\ell(\w_{t};\z_t^{(1)}) - \ell(\wb_t^{(1)};\z_t^{(1)}), \ldots, \ell(\w_{t};\z_t^{(m)})-\ell(\wb_t^{(m)};\z_t^{(m)}) \right]^\top
\end{equation}
which is a \emph{biased} estimator of the true gradient $\nabla_\q \phi(\w_t,\q_t) = [R_1(\w_t)-R_1^*, \ldots, R_m(\w_t)-R_m^*]^\top$, since
\begin{equation}\label{eqn:exp:grad:2}
\E_{t-1}\left[\g_{q}(\w_t,\q_t) \right] =\left [R_1(\w_t)-R_1(\wb_t^{(1)}), \ldots, R_m(\w_t)-R_m(\wb_t^{(m)})\right]^\top \neq \nabla_\q \phi(\w_t,\q_t)
\end{equation}
where $\E_{t-1} [\cdot]$ represents the expectation conditioned on the randomness until round $t-1$. Thanks to the SMD update in (\ref{eqn:wti:smd}), we know that $R_i(\wb_t^{(i)})$ is close to $R_i^*$, for all $i\in[m]$. As a result, the bias in $\g_{q}(\w_t,\q_t)$, determined by $R_i(\wb_t^{(i)})-R_i^*$, is effectively managed, making it possible to maintain a (nearly) optimal convergence rate.

Equipped with the stochastic gradients in (\ref{eqn:stoch:grad:1}) and (\ref{eqn:stoch:grad:2}), we update $\w_t$ and $\q_t$ by SMD:
\begin{align}
\w_{t+1}= &\argmin_{\w \in \W} \big\{ \eta_t^w \langle \g_{w}(\w_t,\q_t) , \w -\w_t\rangle + B_w(\w,\w_t) \big\} , \label{eqn:upate:wt}  \\
\q_{t+1}=& \argmin_{\q \in \Delta_m} \big\{ \eta_t^q \langle -\g_{q}(\w_t,\q_t) , \q -\q_t\rangle + B_q(\q,\q_t) \big\} \label{eqn:upate:qt}
\end{align}
where $\eta_t^w>0$ and $\eta_t^q>0$ are step sizes. We will maintain the weighted averages of iterates:
\begin{equation} \label{eqn:two:average}
\wb_t= \sum_{j=1}^t \frac{\eta_j^w \w_j }{\sum_{k=1}^t \eta_k^w}, \textrm{ and } \qb_t= \sum_{j=1}^t \frac{\eta_j^q \q_j }{\sum_{k=1}^t \eta_k^q}
\end{equation}
which can be returned as  solutions if necessary. The completed procedure is given in Algorithm~\ref{alg:1}.

\begin{algorithm}[t]
\caption{An Anytime Stochastic Approximation Approach for MERO}
\begin{algorithmic}[1]
\STATE Initialize $\w_1=\w_1^{(1)}=\cdots=\w_1^{(m)}=\argmin_{\w \in \W}\nu_w(\w)$, and $\q_1=\frac{1}{m} \mathbf{1}_m \in \R^m$
\FOR{$t=1$ to $T$}
\STATE For each $i\in[m]$, draw a sample $\z_{t}^{(i)}$ from distribution $\P_i$
\STATE For each $i\in[m]$, calculate $\nabla \ell(\w_t^{(i)};\z_t^{(i)})$ and update $\w_t^{(i)}$ according to (\ref{eqn:wti:smd})
\STATE For each $i\in[m]$, calculate the  weighted average $\wb_t^{(i)}$ in (\ref{eqn:wti:average})
\STATE Construct the stochastic gradients in  (\ref{eqn:stoch:grad:1}) and (\ref{eqn:stoch:grad:2})
\STATE Update $\w_t$ and $\q_t$ according to (\ref{eqn:upate:wt}) and (\ref{eqn:upate:qt}), respectively
\STATE Calculate the the weighted averages $\wb_t$ and $\qb_t$ in (\ref{eqn:two:average})
\ENDFOR
\end{algorithmic}\label{alg:1}
\end{algorithm}

Next, we discuss the theoretical guarantee of Algorithm~\ref{alg:1}. To this end, we first present the optimization error of $\wb_t^{(i)}$ in (\ref{eqn:wti:average}) for each risk function \citep[\S 2.3]{nemirovski-2008-robust}.
\begin{thm} \label{thm:1} Under Assumptions~\ref{ass:1}, \ref{ass:2}, and \ref{ass:3}, by setting $\eta_t^{(i)} = \frac{D}{G\sqrt{t}}$ in Algorithm~\ref{alg:1}, we have
\begin{equation} \label{eqn:exp:risk}
\E \left[R_i(\wb_t^{(i)}) - R_i^*  \right] \leq  \frac{DG(3+\ln t)}{4(\sqrt{t+1} - 1)}, \quad  \forall  i \in [m], \ t \in \zn.
\end{equation}
Furthermore, with probability at least $1-\delta$,
\begin{equation} \label{eqn:exp:high}
R_i(\wb_t^{(i)}) - R_i^* \leq \frac{DG\big[3+\ln t + 16\sqrt{(1+ \ln t) \ln (2 m t^2/\delta) } \big]}{4(\sqrt{t+1} - 1)}, \quad \forall  i \in [m], \ t \in \zn.
\end{equation}
\end{thm}
\begin{remark}
\textnormal{With a varying step size, the expected excess risk in (\ref{eqn:exp:risk}) is  of the order $O((\log t)/\sqrt{t})$, which is higher than the $O(1/\sqrt{t})$ bound achieved with a fixed step size \citep[(2.46)]{nemirovski-2008-robust}. We remark that it is possible to remove the additional $\ln t$ factor in (\ref{eqn:exp:risk}), by using a stronger condition that $\max_{\w , \w'\in \W}  B_w(\w,\w')  \leq D^2$ to replace (\ref{eqn:domain:W}) and adapting the proof of online gradient descent with varying step sizes \citep[ Theorem 1]{zinkevich-2003-online}. However, the $\ln t$ factor appears three times in the high probability bound (\ref{eqn:exp:high}), and  the last occurrence is due to the application of the union bound, which seems unavoidable.}
\end{remark}

Finally, we examine the optimization error of $\wb_t$ and $\qb_t$ for (\ref{eqn:convex:concave}). Because of the biased stochastic gradient in (\ref{eqn:stoch:grad:2}), we cannot apply existing  guarantees of SMD for stochastic convex-concave optimization \citep[\S 3.1]{nemirovski-2008-robust}. Therefore, we provide a novel analysis of SMD with biased gradients, and utilize Theorem~\ref{thm:1} to demonstrate that the bias does not significantly affect the convergence behavior.

\begin{thm} \label{thm:2} Under Assumptions~\ref{ass:1}, \ref{ass:2}, \ref{ass:3}, and \ref{ass:4}, by setting
\begin{equation}\label{eqn:size:all}
\eta_t^{(i)} = \frac{D}{G\sqrt{t}}, \  \eta_t^w = \frac{2  D^2}{\sqrt{(2 D^2 G^2  +  2 \ln m)t} } , \textrm{ and } \eta_t^q = \frac{2  \ln m }{\sqrt{(2 D^2 G^2  +  2 \ln m)t}}
\end{equation}
 in Algorithm~\ref{alg:1}, we have
\begin{equation} \label{eqn:mero:risk}
\begin{split}
 & \E \left[\epsilon_{\phi}(\wb_t, \qb_t)=\max_{\q\in \Delta_m}  \phi(\wb_t,\q)- \min_{\w\in \W}  \phi(\w,\qb_t)  \right] \\
\leq &  \frac{(5+ 3\ln t)\sqrt{2 D^2 G^2  +  2 \ln m} + 2 DG\left[3+\ln t + 16 (1+\sqrt{\ln m}) \sqrt{2(1+\ln t)}\right] (1+ \ln t)}{2 (\sqrt{t+1} - 1)}\\
=& O\left( \frac{\log^2 t + \log^{1/2} m \log^{3/2} t}{\sqrt{t}} \right)
\end{split}
\end{equation}
for all  $t \in \zn$. Furthermore, with probability at least $1-2\delta$,
\begin{equation} \label{eqn:mero:high}
\begin{split}
& \epsilon_{\phi}(\wb_t, \qb_t) = \max_{\q\in \Delta_m}  \phi(\wb_t,\q)- \min_{\w\in \W}  \phi(\w,\qb_t)\\
 \leq&\frac{1}{2 (\sqrt{t+1} - 1)} \left[ \sqrt{2 D^2 G^2  +  2 \ln m}\left(5+3\ln t + 8 \sqrt{  (1+\ln t) \ln \frac{2 t^2}{\delta}} \right)\right.\\
& \left.+ 2DG\left( 3+\ln t + 16\sqrt{(1+ \ln t) \ln \frac{2 m t^2}{\delta} }\right )(1+\ln t) \right] = O\left( \frac{\log^2 t + \log^{1/2} m \log^{3/2} t}{\sqrt{t}} \right)
\end{split}
\end{equation}
for all  $t \in \zn$.
\end{thm}
\begin{remark}
\textnormal{The above theorem demonstrates that Algorithm~\ref{alg:1} converges at a (nearly) optimal rate of $O([\log^2 t + \log^{1/2} m \log^{3/2} t]/\sqrt{t})$, which holds both in expectation and with high probability. Given a fixed number of iterations $T$, this rate is slightly slower than the $O(\sqrt{(\log m)/T})$ rate achieved by the multi-stage approach in Section~\ref{sec:multi:stage}. Nonetheless, a considerable advantage of Algorithm~\ref{alg:1} is its anytime characteristic, indicating it is capable of returning a solution at any round.}\end{remark}

\begin{remark}
\textnormal{Algorithm~\ref{alg:1} also enjoys an interesting interpretation from the perspective of prediction with expert advice (PEA) \citep{bianchi-2006-prediction}. We can consider the $m$ instances of SMD as $m$ experts, each generating a model to minimize its own loss, i.e., the risk (Steps 3-5), and regard the other two instances of SMD as a meta-algorithm, which produces a model by aggregating predictions from all experts (Steps 6-8). It's important to note that existing PEA techniques cannot be directly applied here, as they typically employ the \emph{same} loss  for all experts. PEA becomes  more challenging when each expert selects a loss that may differ from those used by others, and previous studies usually impose stringent conditions on  these losses \citep{Expert_Evaluators}. In this context,  Algorithm~\ref{alg:1} essentially provides a novel way for learning from heterogeneous experts. }
\end{remark}

\begin{remark}
\textnormal{Compared to stochastic algorithms for GDRO \citep{nemirovski-2008-robust,SA:GDRO}, our method requires the maintenance and updating of $m$ additional models, resulting in a higher computational complexity.  When $m$ is very large and the computational cost becomes prohibitive,  we may opt for a simplified version of MERO \citep{DRO:DoReMi}:
\begin{equation}\label{eqn:simple:meor}
\min_{\w \in \W} \max_{\q \in \Delta_m}  \   \left\{\sum_{i=1}^m q_i \left[ R_i(\w) - R_i(\w_r)\right]\right\}
\end{equation}
where $\w_r$ is a reference model shared by all distributions. In (\ref{eqn:simple:meor}), we use the risk of $\w_r$ on each distribution to replace the minimal risk of that distribution, and thus only need to estimate one additional model, namely $\w_r$.  A possible way  for selecting $\w_r$ is  to choose the model that minimizes the average risk across all distributions. To efficiently solve (\ref{eqn:simple:meor}), we could develop a multi-stage stochastic method, akin to the approach in Section~\ref{sec:multi:stage}, or an anytime stochastic method, similar to Algorithm~\ref{alg:1}. }
\end{remark}

\section{Stochastic Approximation of Weighted MERO} \label{sec:W:MERO}
In practice, the costs associated with gathering samples can differ among distributions, and thus we may allocate diverse budgets for different distributions. Inspired by the study of GDRO with varying sample sizes across distributions \citep{SA:GDRO}, we also investigate MERO under the imbalanced setting where the number of samples drawn from each distribution is different. 

\subsection{Preliminaries} \label{sec:wmero:basline}
Let $n_i$ represent the sample budget of distribution $\P_i$. For the sake of simplicity, we assume an order of $n_1 \geq n_2 \geq \cdots \geq n_m$. To satisfy the budget, we can simply run the multi-stage approach in Section~\ref{sec:multi:stage} with $T=n_m/3$, or Algorithm~\ref{alg:1} for $n_m$ iterations. However, the optimization error decreases only at an  $O(\sqrt{(\log m)/n_m})$ or $\O(\sqrt{(\log m)/n_m})$ rate, which is determined by the smallest budget $n_m$. In other words, for a distribution $\P_i$ with $n_i > n_m$, the extra $n_i-n_m$ samples are wasted.

Analogous to the weighted GDRO in (\ref{eqn:convex:concave:weight:GDRO}), we also formulate a weighted version of MERO:
\begin{equation} \label{eqn:convex:concave:weight}
\min_{\w \in \W} \max_{\q \in \Delta_m}  \  \left\{\varphi(\w,\q)= \sum_{i=1}^m q_i p_i \big [R_i(\w)- R_i^*\big] \right\}
\end{equation}
where the value of the weight $p_i$ will be determined later. As demonstrated by \citet{Regret:RML:DS}, the heterogeneous scaling in (\ref{eqn:convex:concave:weight}) allows us to establish distribution-specific bounds for the excess risk. In our paper, the goal is to make the excess risk of distribution $\P_i$ reducing at an $O(\sqrt{(\log m)/n_i})$ rate. Again, the optimization problem in (\ref{eqn:convex:concave:weight}) is more challenging than the counterpart in (\ref{eqn:convex:concave:weight:GDRO}), due to the existence of $R_i^*$s. We notice that because the budgets are fixed and known, there is no need to pursue the anytime ability. In the following, we will develop a two-stage stochastic approach for weighted MERO.

\subsection{A Two-Stage Stochastic Approximation Approach for Weighted MERO} \label{sec:two:stage}
Our approach consists of two stages: minimizing each risk and minimizing an approximate problem. The complete procedure is summarized in Algorithm~\ref{alg:2}, where Steps 1-8 belong to the first stage, and the rest correspond to the second stage.

\paragraph{Stage 1: Minimizing the risk} Similar to the first stage in Section~\ref{sec:multi:stage}, we will deploy an instance of SMD to minimize each individual risk $R_i(\cdot)$. The difference is that the number of iterations is set to be $n_i/2$  for distribution $\P_i$. Consequently, a larger budget yields a smaller error. Because the total number of iterations is fixed, we will use a fixed step size for each SMD. Specifically, the update rule for the $i$-th distribution at the $t$-th round is given by
\begin{equation} \label{eqn:wti:smd:fixed}
\w_{t+1}^{(i)}= \argmin_{\w \in \W} \left\{\eta^{(i)} \langle \nabla \ell(\w_t^{(i)};\z_t^{(i)}) , \w -\w_t^{(i)}\rangle + B_w(\w,\w_t^{(i)}) \right\}
\end{equation}
where $\w_t^{(i)}$ is the current solution, $\z_t^{(i)}$ is a random sample drawn from $\P_i$, and $ \eta^{(i)} >0$ is the step size. After $n_i/2$ iterations, we will use the average of iterates $\wb^{(i)}= \frac{1}{n_i/2} \sum_{t=1}^{n_i/2} \w_{t}^{(i)}$ as an approximate solution to $\min_{\w \in \W} R_i(\w)$.

Similar to Theorem~\ref{thm:1}, we have the following guarantee for the excess risk of each $\wb^{(i)}$.
\begin{thm} \label{thm:3} Under Assumptions~\ref{ass:1}, \ref{ass:2}, and \ref{ass:3}, by setting $\eta^{(i)}= \frac{2D}{G\sqrt{n_i}}$ in Algorithm~\ref{alg:2}, we have
\begin{equation} \label{eqn:exp:risk:fixed}
\E \left[R_i(\wb^{(i)}) - R_i^*  \right] \leq  \frac{2DG}{\sqrt{n_i}}, \quad  \forall  i \in [m].
\end{equation}
Furthermore,  with probability at least $1-\delta$, we have
\begin{equation} \label{eqn:exp:high:fixed}
R_i(\wb^{(i)}) - R_i^* \leq \frac{2DG}{\sqrt{n_i}} \left(1+4 \sqrt{2\ln \frac{m}{\delta}}\right), \quad \forall  i \in [m].
\end{equation}
\end{thm}
Compared with Theorem~\ref{thm:1}, the upper bounds do not contain the additional $\ln t$ factor.

\paragraph{Stage 2: Minimizing an approximate problem} Based on the solutions from the first stage, we construct the  problem below
\begin{equation} \label{eqn:convex:concave:weight:appro}
\min_{\w \in \W} \max_{\q \in \Delta_m}  \  \left\{\widehat{\varphi}(\w,\q)= \sum_{i=1}^m q_i p_i \big [R_i(\w)- R_i(\wb^{(i)})\big] \right\}
\end{equation}
to approximate (\ref{eqn:convex:concave:weight}). The following lemma shows that, the optimization error of any $(\wb,\qb)$ for  (\ref{eqn:convex:concave:weight}) is close to that for (\ref{eqn:convex:concave:weight:appro}), provided $R_i(\wb^{(i)})-R_i^*$ is small, for all $i\in[m]$. As a result, we can focus on optimizing (\ref{eqn:convex:concave:weight:appro}), which is more manageable, as the stochastic gradient of $\widehat{\varphi}(\w,\q)$ can be easily constructed.
\begin{lem} \label{lem:diff:opt}
For any $(\wb,\qb)$, we have
\[
\epsilon_{\varphi}(\wb, \qb) \leq  \epsilon_{\widehat{\varphi}}(\wb, \qb)+ 2\max_{i \in [m]} \left\{p_i \big[ R_i(\wb^{(i)})- R_i^*\big]\right\}.
\]
\end{lem}
From Theorem~\ref{thm:3}, we know that a larger $n_i$ leads to a smaller excess risk $R_i(\wb^{(i)})- R_i^*$. Therefore, when the budget $n_i$ is  large, we can set a large $p_i$ without influencing $\max_{i \in [m]} \left\{p_i [ R_i(\wb^{(i)})- R_i^*]\right\}$, which is crucial for achieving faster rates in distributions with a larger number of samples.
\begin{algorithm}[t]
\caption{A Two-Stage Stochastic Approximation Approach for Weighted MERO}
{\bf Input}: Step sizes: $\eta^{(1)},\ldots,\eta^{(m)}$, $\eta_w$ and $\eta_q$
\begin{algorithmic}[1]
\STATE Initialize $\w_1^{(1)}=\cdots=\w_1^{(m)}=\argmin_{\w \in \W}\nu_w(\w)$
\FOR{$i=1$ to $m$}
\FOR{$t=1$ to $n_i/2$}
\STATE Draw a sample $\z_{t}^{(i)}$ from distribution $\P_i$
\STATE Calculate $\nabla \ell(\w_t^{(i)};\z_t^{(i)})$ and update $\w_t^{(i)}$ according to (\ref{eqn:wti:smd:fixed})
\ENDFOR
\STATE Calculate the average of iterates $\wb^{(i)}= \frac{1}{n_i/2} \sum_{t=1}^{n_i/2} \w_{t+1}^{(i)}$
\ENDFOR \vspace{.5ex}
\STATE Initialize $\w_1'=\argmin_{\w \in \W}\nu_w(\w)$, and $\q_1'=\frac{1}{m} \mathbf{1}_m \in \R^m$
\FOR{$t=1$ to $n_m/4$}
\STATE For each $i\in[m]$, draw $n_i/n_m$ samples $\{\z_t^{(i,j)}: j=1,\ldots,n_i/n_m\}$ from distribution $\P_i$
\STATE Construct the stochastic gradients defined in (\ref{eqn:stoch:grad:weight1})
\STATE Calculate $\w_{t+1}$ and $\q_{t+1}$ according to (\ref{eqn:upate:wt:mirror}) and (\ref{eqn:upate:qt:mirror}), respectively
\STATE For each $i\in[m]$, draw $n_i/n_m$ samples $\{\zh_t^{(i,j)}: j=1,\ldots,n_i/n_m\}$ from distribution $\P_i$
\STATE Construct the stochastic gradients defined in (\ref{eqn:stoch:grad:weight2})
\STATE Calculate $\w_{t+1}'$ and $\q_{t+1}'$ according to (\ref{eqn:upate:wt:mirror:2}) and (\ref{eqn:upate:qt:mirror:2}), respectively
\ENDFOR
\RETURN $\wb=\frac{4}{n_m} \sum_{t=2}^{1+n_m/4} \w_t$  and $\qb=\frac{4}{n_m} \sum_{t=2}^{1+n_m/4} \q_t$
\end{algorithmic}\label{alg:2}
\end{algorithm}

Inspired by  \citet[Algorithm~4]{SA:GDRO}, we solve (\ref{eqn:convex:concave:weight:appro}) by stochastic mirror-prox algorithm (SMPA) \citep{Nemirovski:SMP}. One notable merit of SMPA is that its optimization error depends on the variance of the gradient. For a distribution $\P_i$ with a larger $n_i$, we can leverage  mini-batches \citep{NIPS2007_Topmoumoute,ICML:13:Zhang:logT} to estimate $R_i(\w)- R_i(\wb^{(i)})$ in (\ref{eqn:convex:concave:weight:appro}) more accurately, (i.e., with a smaller variance), which again makes it possible to use a larger $p_i$.

Within the framework of SMPA, we keep two sets of solutions: $(\w_t, \q_t)$ and $(\w_t', \q_t')$. In the $t$-th iteration, we first draw $n_i/n_m$ samples from every distribution $\P_i$, denoted by $\z_t^{(i,1)}, \ldots, \z_t^{(i,n_i/n_m)}$. Then, we use them to construct stochastic gradients of $\widehat{\varphi}(\w,\q)$  at $(\w_t', \q_t')$:
\begin{equation}\label{eqn:stoch:grad:weight1}
\begin{split}
&\g_{w}(\w_t',\q_t') = \sum_{i=1}^m q_{t,i}' p_i \left(\frac{n_m}{n_i} \sum_{j=1}^{n_i/n_m} \nabla \ell(\w_{t}';\z_t^{(i,j)}) \right), \\
&\g_{q}(\w_t',\q_t')\\
= &\left[ p_1 \frac{ n_m}{n_1} \sum_{j=1}^{n_1/n_m} \big[ \ell(\w_{t}';\z_t^{(1,j)}) -  \ell(\wb^{(1)};\z_t^{(1,j)}) \big], \ldots, p_m \big[\ell(\w_{t}';\z_t^{(m)})-  \ell(\wb^{(m)};\z_t^{(m)}) \big] \right]^\top .
\end{split}
\end{equation}
It is easy to verify that $\E_{t-1}[\g_{w}(\w_t',\q_t')]=\nabla_\w \widehat{\varphi}(\w_t',\q_t')$ and $\E_{t-1}[\g_{q}(\w_t',\q_t')]=\nabla_\q \widehat{\varphi}(\w_t',\q_t')$. Based on (\ref{eqn:stoch:grad:weight1}), we use SMD to update $(\w_t', \q_t')$, and obtain $(\w_{t+1}, \q_{t+1})$:
\begin{align}
\w_{t+1}= &\argmin_{\w \in \W} \big\{ \eta_w \langle \g_{w}(\w_t',\q_t') , \w -\w_t'\rangle + B_w(\w,\w_t') \big\}, \label{eqn:upate:wt:mirror}  \\
\q_{t+1}=& \argmin_{\q \in \Delta_m} \big\{ \eta_q \langle -\g_{q}(\w_t',\q_t') , \q -\q_t'\rangle + B_q(\q,\q_t') \big\} \label{eqn:upate:qt:mirror}
\end{align}
where $\eta_w>0$ and $\eta_q>0$ are step sizes. Next, we draw another $n_i/n_m$ samples from every distribution $\P_i$  to construct stochastic gradients at $(\w_{t+1}, \q_{t+1})$:
\begin{equation}\label{eqn:stoch:grad:weight2}
\begin{split}
&\g_{w}(\w_{t+1},\q_{t+1}) =\sum_{i=1}^m q_{t+1,i} p_i \left(\frac{n_m}{n_i} \sum_{j=1}^{n_i/n_m} \nabla \ell(\w_{t+1};\zh_t^{(i,j)}) \right), \\
&\g_{q}(\w_{t+1},\q_{t+1})\\
= &\left[p_1 \frac{n_m}{n_1}  \sum_{j=1}^{n_1/n_m} \big[\ell(\w_{t+1};\zh_t^{(1,j)})-  \ell(\wb^{(1)};\zh_t^{(1,j)}) \big],  \ldots,  p_m \big[\ell(\w_{t+1};\zh_t^{(m)}) -  \ell(\wb^{(m)};\zh_t^{(m)}) \big]\right]^\top
\end{split}
\end{equation}
where $\zh_t^{(i,1)}, \ldots, \zh_t^{(i,n_i/n_m)}$ are random samples from distribution $\P_i$.  Then, we use them to update $(\w_t', \q_t')$, and obtain $(\w_{t+1}', \q_{t+1}')$:
\begin{align}
\w_{t+1}'= &\argmin_{\w \in \W} \big\{ \eta_w \langle \g_{w}(\w_{t+1},\q_{t+1}) , \w -\w_t'\rangle + B_w(\w,\w_t') \big\} ,\label{eqn:upate:wt:mirror:2}  \\
\q_{t+1}'=& \argmin_{\q \in \Delta_m} \big\{ \eta_q \langle -\g_{q}(\w_{t+1},\q_{t+1}) , \q -\q_t'\rangle + B_q(\q,\q_t') \big\}. \label{eqn:upate:qt:mirror:2}
\end{align}
Recall that after the first stage, we have $n_i/2$ samples left for each distribution $\P_i$. So, we repeat the above process for $n_m/4$ iterations to meet the budget constraints. Finally, we return $\wb=\frac{4}{n_m} \sum_{t=2}^{1+n_m/4} \w_t$  and $\qb=\frac{4}{n_m} \sum_{t=2}^{1+n_m/4} \q_t$ as solutions.

To analyze the performance of Algorithm~\ref{alg:2}, we  introduce two additional assumptions.
\begin{ass}\label{ass:5} All the risk functions are $L$-smooth, i.e.,
\begin{equation}\label{eqn:smooth:R}
\|\nabla R_i(\w) - \nabla R_i(\w')\|_{w,*} \leq  L \|\w-\w'\|_w, \ \forall  \w,  \w' \in \W, i \in [m].
\end{equation}
\end{ass}
The assumption of smoothness is essential for achieving a convergence rate that depends on the variance \citep{Lan:SCO}.
\begin{ass} \label{ass:6} The dual norm $\|\cdot\|_{w,*}$ is $\kappa$-regular for some small constant $\kappa \geq 1$.
\end{ass}
The condition of regularity plays a role when examining the impact of mini-batches on stochastic gradients. For a comprehensive definition, please consult the work of \citet{Vector:Martingales}.

Then, we have the following theorem regarding the excess risk of  $\wb$ on every distribution.
\begin{thm} \label{thm:6}Define
\begin{equation} \label{eqn:mirror:parameters}
\begin{split}
&p_{\max}=\max_{i\in[m]} p_i, \quad \omega_{\max}=  \max_{i \in [m]}  \frac{p_i^2 n_m}{n_i} , \quad  r_{\max}=\max_{i \in [m]} \frac{p_i }{\sqrt{n_i}}\\
&\Lt= 2\sqrt{2} p_{\max} (D^2 L +  D^2  G \sqrt{\ln m}), \textrm{ and } \sigma^2=2 c \omega_{\max}(\kappa D^2  G^2 +   \ln^2 m)
\end{split}
\end{equation}
where $c>0$ is an absolute constant. Under Assumptions~\ref{ass:1}, \ref{ass:2}, \ref{ass:3}, \ref{ass:4}, \ref{ass:5}, and \ref{ass:6}, and setting
\[
\eta^{(i)}= \frac{2D}{G\sqrt{n_i}}, \ \eta_w =2  D^2 \min \left( \frac{1}{\sqrt{3} \Lt}, 2\sqrt{\frac{2}{7 \sigma^2 n_m}} \right), \textrm{ and } \eta_q = 2  \min \left( \frac{1}{\sqrt{3} \Lt}, 2\sqrt{\frac{2}{7 \sigma^2 n_m}} \right) \ln m
\]
in Algorithm~\ref{alg:2}, with probability at least  $1-2\delta$, we have
\begin{equation} \label{eqn:final:risk}
\begin{split}
&R_i(\wb)- R_i^* \leq  \frac{1}{p_i} p_{\varphi}^* \\
& + \frac{1}{p_i} \left[  \frac{14 \Lt}{n_m} + \sqrt{\frac{\sigma^2}{n_m}}  \left(\frac{28}{\sqrt{3}} + 7\sqrt{6 \log \frac{2}{\delta}} + \frac{28\sqrt{2}}{n_m} \log \frac{2}{\delta} \right) + 4DG \left(1+4 \sqrt{2\ln \frac{m}{\delta}}\right) r_{\max} \right]
\end{split}
\end{equation}
where $p_{\varphi}^*$ is the optimal value of (\ref{eqn:convex:concave:weight}).
\end{thm}

To simplify the upper bound in (\ref{eqn:final:risk}), we set $p_i$ as
\begin{equation} \label{eqn:pi:mirror}
p_i = \frac{1/\sqrt{n_m} + 1}{1/\sqrt{n_m} + \sqrt{n_m/n_i}}
\end{equation}
which is proposed by \citet{SA:GDRO} for weighted GDRO, and obtain the following corollary.
\begin{cor} \label{cor:1} Under the condition of Theorem~\ref{thm:6} and (\ref{eqn:pi:mirror}), with high probability, we have
\[
R_i(\wb)- R_i^*  =  \frac{1}{p_i} p_{\varphi}^* +  O\left(\left( \frac{1}{n_m} + \frac{1}{\sqrt{n_i}} \right) \sqrt{\kappa+ \ln^2 m}\right).
\]
\end{cor}
\begin{remark}
\textnormal{We observe that for a distribution $\P_i$ with $n_i \leq n_m^2$, the excess risk diminishes at a rate of $O((\log m)/\sqrt{n_i})$, a significant improvement over the $\O(\sqrt{(\log m)/n_m})$ rate outlined in Section~\ref{sec:wmero:basline}, until it approaches $p_{\varphi}^*/p_i$.  For  any $\P_i$ with a very large budget, i.e., $n_i  > n_m^2$, it attains an $O((\log m) /n_m)$ rate, which almost matches the convergence rate of deterministic convex-concave saddle-point optimization \citep{nemirovski-2005-prox}.  Our convergence rate is of the same order as that of  \citet[Theorem 6]{SA:GDRO}; however, the implications differ due to variations in the objective functions, specifically between weighted GDRO and weighted MERO.}
\end{remark}

Although the exact value of $p_{\varphi}^*$ is generally unknown, we can expect it to be relatively small when there exists a single model that performs well on all distributions. In particular, if all distributions are aligned, we can prove that $p_{\varphi}^*=0$ \citep[Corollary 9]{Regret:RML:DS}, leading to the following corollary.
\begin{cor} \label{cor:2} Suppose there exists a model $\w_* \in \W$ such that $R_i(\w_*)=R_i^*$ for all $i\in[m]$. Under the condition of Theorem~\ref{thm:6} and (\ref{eqn:pi:mirror}), with high probability, we have
\[
R_i(\wb)- R_i^*  =   O\left(\left( \frac{1}{n_m} + \frac{1}{\sqrt{n_i}} \right) \sqrt{\kappa+ \ln^2 m}\right).
\]
\end{cor}
\begin{remark}
\textnormal{When all distributions are aligned, the above corollary offers upper bounds for the \emph{standard} excess risk, making it more interpretable than the theoretical guarantee provided by \citet{SA:GDRO}.}
\end{remark}

\section{Analysis}
In this section, we present proofs of main theorems.
\subsection{Proof of Theorem~\ref{thm:1}}
The analysis closely adheres to the content in Section~2.3 of \citet{nemirovski-2008-robust}, and for the sake of completeness, we present the proof here.

Let $\w_*^{(i)} \in \argmin_{\w \in \W} R_i(\w)$ be the optimal solution that minimizes $R_i(\cdot)$. From the property of mirror descent, e.g., Lemma 2.1 of  \citet{nemirovski-2008-robust}, we have
\begin{equation} \label{eqn:mr:risk:1}
\begin{split}
&\eta_j^{(i)} \langle \nabla \ell(\w_j^{(i)};\z_j^{(i)}) , \w_j^{(i)} - \w_*^{(i)} \rangle \\
\leq & B_w(\w_*^{(i)}, \w_j^{(i)})-B_w(\w_*^{(i)}, \w_{j+1}^{(i)})   + \frac{(\eta_j^{(i)})^2 }{2} \| \nabla \ell(\w_j^{(i)};\z_j^{(i)})  \|_{w,*}^2 \\
\overset{\text{(\ref{eqn:gradient})}}{\leq} & B_w(\w_*^{(i)}, \w_j^{(i)}) - B_w(\w_*^{(i)}, \w_{j+1}^{(i)}) + \frac{(\eta_j^{(i)})^2 G^2 }{2}.
\end{split}
\end{equation}
Thus, we have
\[
\begin{split}
&\eta_j^{(i)} \langle \nabla R_i(\w_j^{(i)}) , \w_j^{(i)} - \w_*^{(i)} \rangle \\
=& \eta_j^{(i)} \langle \nabla \ell(\w_j^{(i)};\z_j^{(i)}) , \w_j^{(i)} - \w_*^{(i)} \rangle + \eta_j^{(i)} \langle \nabla R_i(\w_j^{(i)}) -\nabla \ell(\w_j^{(i)};\z_j^{(i)}) , \w_j^{(i)} - \w_*^{(i)} \rangle\\
\overset{\text{(\ref{eqn:mr:risk:1})}}{\leq} & B_w(\w_*^{(i)}, \w_j^{(i)}) - B_w(\w_*^{(i)}, \w_{j+1}^{(i)}) + \frac{(\eta_j^{(i)})^2 G^2 }{2} + \eta_j^{(i)} \langle \nabla R_i(\w_j^{(i)}) -\nabla \ell(\w_j^{(i)};\z_j^{(i)}) , \w_j^{(i)} - \w_*^{(i)} \rangle .
\end{split}
\]
Summing the above inequality over $j=1,\ldots,t$, we have
\begin{equation} \label{eqn:mr:risk:2}
\begin{split}
&\sum_{j=1}^t \eta_j^{(i)} \langle \nabla R_i(\w_j^{(i)}) , \w_j^{(i)} - \w_*^{(i)} \rangle \\
\leq & B_w(\w_*^{(i)}, \w_1^{(i)}) + \frac{ G^2 }{2} \sum_{j=1}^t  (\eta_j^{(i)})^2+ \sum_{j=1}^t\eta_j^{(i)} \langle \nabla R_i(\w_j^{(i)}) -\nabla \ell(\w_j^{(i)};\z_j^{(i)}) , \w_j^{(i)} - \w_*^{(i)} \rangle \\
\overset{\text{(\ref{eqn:domain:W})}}{\leq}& D^2 + \frac{ G^2 }{2} \sum_{j=1}^t  (\eta_j^{(i)})^2+ \sum_{j=1}^t\eta_j^{(i)} \langle \nabla R_i(\w_j^{(i)}) -\nabla \ell(\w_j^{(i)};\z_j^{(i)}) , \w_j^{(i)} - \w_*^{(i)} \rangle.
\end{split}
\end{equation}

From the convexity of the risk function, we have
\begin{equation} \label{eqn:mr:risk:3}
\begin{split}
& R_i(\wb_t^{(i)}) - R_i(\w_*^{(i)})  = R_i\left(\sum_{j=1}^t \frac{\eta_j^{(i)} \w_j^{(i)} }{\sum_{k=1}^t \eta_k^{(i)}}\right) - R_i(\w_*^{(i)}) \\
\leq &  \left(\sum_{j=1}^t \frac{\eta_j^{(i)}  }{\sum_{k=1}^t \eta_k^{(i)}} R_i(\w_j^{(i)})\right) - R_i(\w_*^{(i)}) =  \sum_{j=1}^t \frac{1}{\sum_{k=1}^t \eta_k^{(i)}}  \eta_j^{(i)} \left( R_i(\w_j^{(i)}) - R_i(\w_*^{(i)}) \right) \\
\leq & \sum_{j=1}^t \frac{1}{\sum_{k=1}^t \eta_k^{(i)}}  \eta_j^{(i)} \langle \nabla R_i(\w_j^{(i)}) , \w_j^{(i)} - \w_*^{(i)} \rangle \\
\overset{\text{(\ref{eqn:mr:risk:2})}}{\leq} & \frac{D^2 + \frac{ G^2 }{2} \sum_{j=1}^t  (\eta_j^{(i)})^2 + \sum_{j=1}^t\eta_j^{(i)} \langle \nabla R_i(\w_j^{(i)}) -\nabla \ell(\w_j^{(i)};\z_j^{(i)}) , \w_j^{(i)} - \w_*^{(i)} \rangle}{\sum_{j=1}^t \eta_j^{(i)}}.
\end{split}
\end{equation}

Define
\begin{equation} \label{eqn:delta:mart}
\delta_j^{(i)}= \eta_j^{(i)} \langle \nabla R_i(\w_j^{(i)}) -\nabla \ell(\w_j^{(i)};\z_j^{(i)}) , \w_j^{(i)} - \w_*^{(i)} \rangle.
\end{equation}
Recall that $\w_j^{(i)}$ and $\w_*^{(i)}$ do not depend on $\z_j^{(i)}$, and thus
\[
\E_{j-1} [\delta_j^{(i)}] = \eta_j^{(i)} \left\langle \E_{j-1} \left[\nabla R_i(\w_j^{(i)}) -\nabla \ell(\w_j^{(i)};\z_j^{(i)}) \right] , \w_j^{(i)} - \w_*^{(i)} \right \rangle=0.
\]
So, $\delta_1^{(i)},\ldots,\delta_t^{(i)}$ is a martingale difference sequence.
\subsubsection{The Expectation Bound}
Taking expectation over (\ref{eqn:mr:risk:3}), we have
\begin{equation} \label{eqn:mr:risk:4}
\E \left[R_i(\wb_t^{(i)}) - R_i(\w_*^{(i)})  \right] \leq \frac{2D^2 + G^2  \sum_{j=1}^t  (\eta_j^{(i)})^2 }{2\sum_{j=1}^t \eta_j^{(i)}}.
\end{equation}
By setting $\eta_j^{(i)} = \frac{D}{G\sqrt{j}}$,  we have
\begin{equation} \label{eqn:mr:risk:5}
\E \left[R_i(\wb_t^{(i)}) - R_i(\w_*^{(i)})  \right] \leq \frac{2DG + DG  \sum_{j=1}^t \frac{1}{j}}{2 \sum_{j=1}^t \frac{1}{\sqrt{j}}}.
\end{equation}
Notice that
\begin{equation} \label{eqn:mr:risk:6}
\begin{split}
 \sum_{j=1}^t \frac{1}{j} &\leq 1 + \int_1^t \frac{1}{x} dx = 1 + \ln x|_1^t = 1+ \ln t, \\
\sum_{j=1}^t \frac{1}{\sqrt{j}} & \geq  \int_1^{t+1} \frac{1}{\sqrt{x}} dx = 2\sqrt{x} |_{1}^{t+1} = 2 (\sqrt{t+1} - 1).
\end{split}
\end{equation}
From (\ref{eqn:mr:risk:5}) and (\ref{eqn:mr:risk:6}), we have
\[
\E \left[R_i(\wb_t^{(i)}) - R_i(\w_*^{(i)})  \right] \leq  \frac{DG(3+\ln t)}{4(\sqrt{t+1} - 1)}
\]
which proves (\ref{eqn:exp:risk}).

\subsubsection{The High Probability Bound} \label{sec:high:pro}
To establish the high probability bound, we make use of the Hoeffding-Azuma inequality for martingales \citep{bianchi-2006-prediction}.
\begin{lem} \label{lem:azuma}
Let $V_1, V_2,  \ldots$  be a martingale difference sequence with respect to some sequence $X_1, X_2, \ldots$ such that $V_i \in [A_i , A_i + c_i ]$ for some random variable $A_i$, measurable with respect to $X_1, \ldots , X_{i-1}$ and a positive constant $c_i$. If $S_n = \sum_{i=1}^n V_i$, then for any
$t > 0$,
\[
\Pr[ S_n > t] \leq \exp \left( -\frac{2t^2}{\sum_{i=1}^n c_i^2} \right).
\]
\end{lem}
To apply the above lemma, we need to show that $|\delta_j^{(i)}|$ is bounded. We have
\begin{equation} \label{eqn:mr:risk:7}
\begin{split}
& \left\| \nabla R_i(\w_j^{(i)}) -\nabla \ell(\w_j^{(i)};\z_j^{(i)}) \right \|_{w,*} \leq  \left\| \nabla R_i(\w_j^{(i)})\right \|_{w,*} + \left\|\nabla \ell(\w_j^{(i)};\z_j^{(i)}) \right \|_{w,*} \\
\leq &  \E_{t-1} \left[\left\| \nabla \ell(\w_j^{(i)};\z_j^{(i)}) \right \|_{w,*} \right] + \left\|\nabla \ell(\w_j^{(i)};\z_j^{(i)}) \right \|_{w,*} \overset{\text{(\ref{eqn:gradient})}}{\leq} 2G,
\end{split}
\end{equation}
\begin{equation} \label{eqn:mr:risk:8}
\begin{split}
& \left\| \w_j^{(i)} - \w_*^{(i)} \right\|_{w} \leq \left\| \w_j^{(i)} - \w_1^{(i)}   \right\|_{w} + \left\| \w_1^{(i)} - \w_*^{(i)} \right\|_{w} \\
\leq & \left( \sqrt{2 B_w(\w_j^{(i)}, \w_1^{(i)})} + \sqrt{2 B_w(\w_*^{(i)} , \w_1^{(i)})} \right) \overset{\text{(\ref{eqn:domain:W})}}{\leq} 2\sqrt{2} D.
\end{split}
\end{equation}
As a result,
\begin{equation} \label{eqn:bound:delta}
\begin{split}
|\delta_j^{(i)}|= & \eta_j^{(i)} \left| \langle \nabla R_i(\w_j^{(i)}) -\nabla \ell(\w_j^{(i)};\z_j^{(i)}) , \w_j^{(i)} - \w_*^{(i)} \rangle \right| \\
\leq &  \eta_j^{(i)} \left\| \nabla R_i(\w_j^{(i)}) -\nabla \ell(\w_j^{(i)};\z_j^{(i)}) \right \|_{w,*}  \left\| \w_j^{(i)} - \w_*^{(i)} \right\|_{w} \overset{\text{(\ref{eqn:mr:risk:7}),(\ref{eqn:mr:risk:8})}}{\leq} 4 \sqrt{2}\eta_j^{(i)}  DG .
\end{split}
\end{equation}
From Lemma~\ref{lem:azuma}, with probability at least $1-\delta/[2m t^2]$, we have
\begin{equation} \label{eqn:mr:risk:9}
\sum_{j=1}^t \delta_j^{(i)} \leq 8 DG \sqrt{ \sum_{j=1}^t (\eta_j^{(i)})^2 \left(\ln \frac{2 mt^2}{\delta}\right)} .
\end{equation}
Substituting (\ref{eqn:mr:risk:9}) into (\ref{eqn:mr:risk:3}), with probability at least $1-\delta/[2 m t^2]$, we have
\[
\begin{split}
 R_i(\wb_t^{(i)}) - R_i(\w_*^{(i)})  \leq & \frac{D^2 + \frac{ G^2 }{2} \sum_{j=1}^t  (\eta_j^{(i)})^2 + 8 DG \sqrt{ \sum_{j=1}^t (\eta_j^{(i)})^2 \left(\ln \frac{2 m t^2}{\delta}\right)} }{\sum_{j=1}^t \eta_j^{(i)}} \\
=& \frac{2DG + DG  \sum_{j=1}^t \frac{1}{j} +16 DG \sqrt{ \sum_{j=1}^t \frac{1}{j} \left( \ln \frac{2 m t^2}{\delta} \right)}}{2 \sum_{j=1}^t \frac{1}{\sqrt{j}}}\\
\overset{\text{(\ref{eqn:mr:risk:6})}}{\leq} &  \frac{DG \big[3+\ln t + 16\sqrt{ (1+ \ln t) \ln (2 m t^2/\delta) } \big]}{4(\sqrt{t+1} - 1)}.
\end{split}
\]
We complete the proof by taking the union bound over all $i\in[m]$ and $t \in \zn$,  and using the well-known fact
\[
\sum_{t=1}^\infty \frac{1}{t^2} = \frac{\pi^2}{6} \leq 2.
\]

\subsection{Proof of Theorem~\ref{thm:2}}
Following the analysis of \citet[\S 3.1]{nemirovski-2008-robust}, we will first combine the two update rules in (\ref{eqn:upate:wt}) and (\ref{eqn:upate:qt}) into a single one.

\subsubsection{Merging the Two Update Rules in (\ref{eqn:upate:wt}) and (\ref{eqn:upate:qt})}
Let $\mathcal{E}$ be the space in which $\W$ resides. We equip the Cartesian product $\mathcal{E} \times \R^m$ with the following norm and dual norm:
\begin{equation} \label{eqn:norm:new}
\big\|(\w,\q) \big\|= \sqrt{ \frac{1}{2 D^2} \|\w\|_{w}^2 + \frac{1}{2 \ln m} \|\q\|_1^2 } , \textrm{ and } \big\|(\u,\bv)\big\|_*= \sqrt{2 D^2 \|\u\|_{w,*}^2 +  2 \|\bv\|_\infty^2   \ln m}.
\end{equation}
We use the notation $\x=(\w,\q)$, and equip the set $\W \times \Delta_m$  with  the distance-generating function
\begin{equation} \label{eqn:dis:fun}
\nu(\x) =\nu(\w,\q)  = \frac{1}{2 D^2} \nu_w(\w) +\frac{1}{2 \ln m} \nu_q(\q) .
\end{equation}
It is easy to verify that $\nu(\x)$ is $1$-strongly convex w.r.t.~the norm $\|\cdot\|$ in (\ref{eqn:norm:new}). Let $B(\cdot,\cdot)$ be the Bregman distance associated with $\nu(\cdot)$:
\begin{equation} \label{eqn:Bregman:merge}
\begin{split}
B(\x,\x')= &\nu(\x) - \big[\nu(\x') + \langle \nabla \nu(\x') , \x-\x'\rangle\big] \\
=& \frac{1}{2 D^2} \left(  \nu_w(\w) - \big[\nu_w(\w') + \langle \nabla \nu_w(\w') , \w-\w'\rangle\big]
 \right)\\
&+\frac{1}{2 \ln m} \left( \nu_q(\q) - \big[\nu_q(\q') + \langle \nabla \nu_q(\q') , \q-\q'\rangle\big]
\right)\\
=& \frac{1}{2 D^2} B_w(\w,\w') + \frac{1}{2 \ln m} B_q(\q,\q')
\end{split}
\end{equation}
where $\x'=(\w',\q')$. Recall the definitions of $\bo_w$ and $\bo_q$ in Section~\ref{sec:pre}, and we have
\[
(\bo_w, \bo_q) =\argmin_{(\w,\q) \in \W \times \Delta_m} \nu(\w,\q).
\]
Then, we can show that the domain $\W \times \Delta_m$ is bounded since
\begin{equation} \label{eqn:domain:merge}
\max_{(\w,\q) \in \W \times \Delta_m} B([\w,\q],[\bo_w, \bo_q]) = \frac{1}{2 D^2}  \max_{\w \in \W} B_w(\w,\bo_w) + \frac{1}{2 \ln m} \max_{\q \in \Delta_m} B_q(\q,\bo_q) \overset{\text{(\ref{eqn:domain:W})}}{\leq} 1.
\end{equation}

With the above configurations, (\ref{eqn:upate:wt}) and (\ref{eqn:upate:qt}) are equivalent to
\begin{equation} \label{eqn:update:xt}
\x_{t+1}= \argmin_{\x \in \W \times \Delta_m} \Big\{ \eta_t \big \langle [ \g_{w}(\w_t,\q_t), -\g_{q}(\w_t,\q_t)] , \x -\x_t \big\rangle + B(\x,\x_t) \Big\}
\end{equation}
where $\eta_t>0$ is the step size that satisfies
\begin{equation}\label{eqn:relation:step}
\eta_t^w =2 \eta_t D^2, \textrm{ and } \eta_t^q = 2 \eta_t \ln m .
\end{equation}
And in the beginning, we set $\x_1=\argmin_{\x \in \W \times \Delta_m}\nu(\x) = (\w_1,\q_1)=(\bo_w, \bo_q)$.

\subsubsection{Analysis of SMD with Biased Stochastic Gradients}
To simplify the notation, we define
\begin{equation}\label{eqn:true:gradient}
\begin{split}
&F(\w_t,\q_t)=[\nabla_\w \phi(\w_t,\q_t), -\nabla_\q \phi(\w_t,\q_t) ]\\
=&\left[ \sum_{i=1}^m q_{t,i} \nabla R_i(\w_t) , -\big[R_1(\w_t)-R_1^*, \ldots,  R_m(\w_t)-R_m^* \big]^\top \right]
\end{split}
\end{equation}
which contains the true gradient of $\phi(\cdot,\cdot)$ at $(\w_t,\q_t)$, and
\begin{equation}\label{eqn:stochastic:gradient}
\begin{split}
&\g(\w_t,\q_t)=[ \g_{w}(\w_t,\q_t); -\g_{q}(\w_t,\q_t)] \\
\overset{\text{(\ref{eqn:stoch:grad:1}),(\ref{eqn:stoch:grad:2})}}{=} &\left[ \sum_{i=1}^m q_{t,i} \nabla \ell(\w_{t};\z_t^{(i)}) ,  -\big[\ell(\w_{t};\z_t^{(1)}) - \ell(\wb_t^{(1)};\z_t^{(1)}), \ldots, \ell(\w_{t};\z_t^{(m)})-\ell(\wb_t^{(m)};\z_t^{(m)})\big]^\top \right]
\end{split}
\end{equation}
which contains the stochastic gradient used in (\ref{eqn:update:xt}). The norm of the stochastic gradient is well-bounded:
\[
\begin{split}
\|\g_{w}(\w_t,\q_t)\|_{w,*} = &\left\|\sum_{i=1}^m q_{t,i} \nabla \ell(\w_{t};\z_t^{(i)})\right\|_{w,*} \leq \sum_{i=1}^m  q_{t,i} \left\|\nabla \ell(\w_{t};\z_t^{(i)})\right\|_{w,*}\overset{\text{(\ref{eqn:gradient})}}{\leq} \sum_{i=1}^m  q_{t,i} G =G,\\
\|\g_{q}(\w_t,\q_t)\|_\infty = &\left \|\big[\ell(\w_{t};\z_t^{(1)}) - \ell(\wb_t^{(1)};\z_t^{(1)}), \ldots, \ell(\w_{t};\z_t^{(m)})-\ell(\wb_t^{(m)};\z_t^{(m)})\big]^\top \right\|_\infty \overset{\text{(\ref{eqn:value})}}{\leq} 1
\end{split}
\]
and thus
\begin{equation} \label{eqn:gradient:bound}
\big\|\g(\w_t,\q_t) \big\|_*= \sqrt{2 D^2 \|\g_{w}(\w_t,\q_t)\|_{w,*}^2 +  2 \|\g_{q}(\w_t,\q_t)\|_\infty^2   \ln m}  \leq  \underbrace{\sqrt{2 D^2 G^2  +  2 \ln m}}_{:=M} .
\end{equation}
The bias of  $\g(\w_t,\q_t)$  is characterized by
\begin{equation} \label{eqn:gradient:bias}
F(\w_t,\q_t) - \E_{t-1} \left[ \g(\w_t,\q_t) \right]= \left[0 , -\big[R_1(\wb_t^{(1)})-R_1^*, \ldots, R_m(\wb_t^{(m)})-R_m^*\big]^\top \right].
\end{equation}

From the convexity-concavity of $\phi(\cdot,\cdot)$, we have \citep[(3.9)]{nemirovski-2008-robust}
\begin{equation} \label{eqn:cvx:cav}
\begin{split}
& \max_{\q\in \Delta_m}  \phi(\wb_t,\q)- \min_{\w\in \W}  \phi(\w,\qb_t) \\
\overset{\text{(\ref{eqn:two:average})}}{=} & \max_{\q\in \Delta_m}  \phi\left(\sum_{j=1}^t \frac{\eta_j^w \w_j }{\sum_{k=1}^t \eta_k^w},\q \right)- \min_{\w\in \W}  \phi\left(\w,\sum_{j=1}^t \frac{\eta_j^q \q_j }{\sum_{k=1}^t \eta_k^q} \right) \\
\overset{\text{(\ref{eqn:relation:step})}}{=} & \max_{\q\in \Delta_m}  \phi\left(\sum_{j=1}^t \frac{\eta_j \w_j }{\sum_{k=1}^t \eta_k},\q \right)- \min_{\w\in \W}  \phi\left(\w,\sum_{j=1}^t \frac{\eta_j \q_j }{\sum_{k=1}^t \eta_k} \right) \\
\leq & \left( \sum_{j=1}^t \eta_j\right)^{-1} \left[ \max_{\q\in \Delta_m}  \sum_{j=1}^t \eta_j \phi(\w_j,\q) -  \min_{\w\in \W} \sum_{j=1}^t \eta_j \phi(\w,\q_j) \right] \\
\leq & \left( \sum_{j=1}^t \eta_j\right)^{-1}  \max_{(\w,\q) \in \W \times \Delta_m} \sum_{j=1}^t \eta_j \Big[ \big\langle  \nabla_\w \phi(\w_j,\q_j), \w_j - \w\rangle - \langle \nabla_\q \phi(\w_j,\q_j) , \q_j -\q \big \rangle \Big]\\
\overset{\text{(\ref{eqn:true:gradient})}}{=} & \left( \sum_{j=1}^t \eta_j\right)^{-1} \max_{\x \in \W \times \Delta_m} \sum_{j=1}^t \eta_j \big\langle  F(\w_j,\q_j), \x_j - \x \big\rangle .
\end{split}
\end{equation}
As a result, we can decompose the optimization error as follows:
\begin{equation} \label{eqn:error:decompo}
\begin{split}
& \max_{\q\in \Delta_m}  \phi(\wb_t,\q)- \min_{\w\in \W}  \phi(\w,\qb_t) \\
\leq & \left( \sum_{j=1}^t \eta_j\right)^{-1} \underbrace{\max_{\x \in \W \times \Delta_m} \sum_{j=1}^t \eta_j \big\langle  \g(\w_j,\q_j), \x_j - \x \big\rangle}_{:=E_1} \\
&+\left( \sum_{j=1}^t \eta_j\right)^{-1} \underbrace{\max_{\x \in \W \times \Delta_m}  \sum_{j=1}^t \eta_j \big \langle  \E_{j-1} \left[ \g(\w_j,\q_j) \right]-\g(\w_j,\q_j) , \x_j - \x \big\rangle}_{:=E_2}\\
&+\left( \sum_{j=1}^t \eta_j\right)^{-1} \underbrace{\max_{\x \in \W \times \Delta_m}  \sum_{j=1}^t \eta_j \big \langle  F(\w_j,\q_j) -\E_{j-1} \left[ \g(\w_j,\q_j) \right], \x_j - \x \big \rangle}_{:=E_3} .
\end{split}
\end{equation}

We proceed to bound the three terms in (\ref{eqn:error:decompo}). To bound the first term $E_1$, we follow the derivation of (\ref{eqn:mr:risk:1}), and have
\begin{equation} \label{eqn:thm2:1}
\begin{split}
\eta_j \big\langle  \g(\w_j,\q_j), \x_j - \x \big\rangle \leq & B(\x, \x_j)-B(\x, \x_{j+1})   + \frac{\eta_j^2 }{2} \|\g(\w_j,\q_j) \|_*^2 \\
\overset{\text{(\ref{eqn:gradient:bound})}}{\leq} &  B(\x, \x_j)-B(\x, \x_{j+1})   +   \frac{M^2}{2}  \eta_j^2 .
\end{split}
\end{equation}
Summing the above inequality over $j=1,\ldots,t$, we have
\begin{equation} \label{eqn:thm2:2}
\sum_{j=1}^t \eta_j \big\langle  \g(\w_j,\q_j), \x_j - \x \big\rangle  \leq B(\x, \x_1)   + \frac{M^2}{2}  \sum_{j=1}^t\eta_j^2 \overset{\text{(\ref{eqn:domain:merge})}}{\leq}  1+ \frac{M^2}{2}  \sum_{j=1}^t\eta_j^2 .
\end{equation}

Next, we consider the second term $E_2$. Because of the maximization operation, the variable $\x$ depends on the randomness of the algorithm, and thus we cannot treat $\eta_j \langle  \E_{j-1} [ \g(\w_j,\q_j) ]-\g(\w_j,\q_j) , \x_j - \x \rangle$,  $j=1,\ldots,t$ as a  martingale difference sequence. To address this challenge, we  make use of the  ``ghost iterate'' technique of \citet{nemirovski-2008-robust}, and develop the following lemma.
\begin{lem} \label{lem:1} Under the condition of Theorem~\ref{thm:2}, we have
\begin{equation} \label{lem:1:exp}
\begin{split}
 \E \left[ \max_{\x \in \W \times \Delta_m}  \sum_{j=1}^t \eta_j \big \langle  \E_{j-1} \left[ \g(\w_j,\q_j) \right]-\g(\w_j,\q_j) , \x_j - \x \big\rangle \right]   \leq  1+ 2M^2  \sum_{j=1}^t\eta_j^2 , \quad \forall  t \in \zn.
\end{split}
\end{equation}
Furthermore, with probability at least $1-\delta$,
\begin{equation} \label{lem:1:high}
\begin{split}
& \max_{\x \in \W \times \Delta_m}  \sum_{j=1}^t \eta_j \big \langle  \E_{j-1} \left[ \g(\w_j,\q_j) \right]-\g(\w_j,\q_j) , \x_j - \x \big\rangle  \\
 \leq &1+ 2M^2  \sum_{j=1}^t\eta_j^2 + 8 M \sqrt{ \sum_{j=1}^t \eta_j^2 \left(\ln \frac{2 t^2}{\delta}\right)}, \quad \forall  t \in \zn.
\end{split}
\end{equation}
\end{lem}

To bound the last term $E_3$, we make use of Theorem~\ref{thm:1}, and prove the following lemma.
\begin{lem} \label{lem:2} Under the condition of Theorem~\ref{thm:2}, we have
\begin{equation} \label{lem:2:exp}
\begin{split}
 &\E \left[\max_{\x \in \W \times \Delta_m}  \sum_{j=1}^t \eta_j \big \langle  F(\w_j,\q_j) -\E_{j-1} \left[ \g(\w_j,\q_j) \right], \x_j - \x \big \rangle \right] \\
\leq &\sum_{j=1}^t  \eta_j  \frac{DG\big[3+\ln j + 16 (1+\sqrt{\ln m}) \sqrt{2(1+\ln j)} \big]}{2(\sqrt{j+1} - 1)}, \quad \forall  t \in \zn.
\end{split}
\end{equation}
Furthermore, with probability at least $1-\delta$,
\begin{equation} \label{lem:2:high}
\begin{split}
&\max_{\x \in \W \times \Delta_m} \sum_{j=1}^t \eta_j \big \langle  F(\w_j,\q_j) -\E_{j-1} \left[ \g(\w_j,\q_j) \right], \x_j - \x \big \rangle \\
\leq& \sum_{j=1}^t \eta_j \frac{DG\big[3+\ln j + 16\sqrt{ (1+ \ln j) \ln (2 m j^2/\delta) } \big]}{2(\sqrt{j+1} - 1)} , \quad \forall  t \in \zn.
\end{split}
\end{equation}
\end{lem}

Combining (\ref{eqn:error:decompo}), (\ref{eqn:thm2:2}),  (\ref{lem:1:exp}) and  (\ref{lem:2:exp}), we have
\[
\begin{split}
 & \E \left[\max_{\q\in \Delta_m}  \phi(\wb_t,\q)- \min_{\w\in \W}  \phi(\w,\qb_t)  \right] \\
\leq & \left( \sum_{j=1}^t \eta_j\right)^{-1} \left( 2+ \frac{5M^2}{2}  \sum_{j=1}^t\eta_j^2 + \sum_{j=1}^t  \eta_j  \frac{DG\big[3+\ln j + 16 (1+\sqrt{\ln m}) \sqrt{2(1+\ln j)}\big ]}{2(\sqrt{j+1} - 1)} \right).
\end{split}
\]
By setting
\begin{equation} \label{eqn:size:main}
\eta_j = \frac{1}{M \sqrt{j}},
\end{equation}
we have
\[
\begin{split}
 & \E \left[\max_{\q\in \Delta_m}  \phi(\wb_t,\q)- \min_{\w\in \W}  \phi(\w,\qb_t)  \right] \\
\leq & \left( \sum_{j=1}^t \frac{1}{\sqrt{j}} \right)^{-1}  \left( 2M+ \frac{5 M}{2} \sum_{j=1}^t \frac{1}{j} + \sum_{j=1}^t  \frac{1}{ \sqrt{j}}  \frac{DG\big[3+\ln j + 16 (1+\sqrt{\ln m}) \sqrt{2(1+\ln j)} \big]}{2(\sqrt{j+1} - 1)} \right).
\end{split}
\]
It is easy to verify that
\begin{equation} \label{eqn:simple:ienqulity}
2(\sqrt{j+1} - 1) \geq \frac{\sqrt{j}}{2}, \ \forall j \in \zn.
\end{equation}
So, we have
\[
\begin{split}
 & \E \left[\max_{\q\in \Delta_m}  \phi(\wb_t,\q)- \min_{\w\in \W}  \phi(\w,\qb_t)  \right] \\
\leq & \left( \sum_{j=1}^t \frac{1}{\sqrt{j}} \right)^{-1}  \left( 2M+ \frac{5 M}{2} \sum_{j=1}^t \frac{1}{j} + 2 DG\left[3+\ln t + 16 (1+\sqrt{\ln m}) \sqrt{2(1+\ln t)}\right]\sum_{j=1}^t  \frac{1}{ j}  \right) \\
\overset{\text{(\ref{eqn:mr:risk:6})}}{\leq} &  \frac{1}{2 (\sqrt{t+1} - 1)}  \left( M  (5+ 3\ln t) + 2 DG\left[3+\ln t + 16 (1+\sqrt{\ln m}) \sqrt{2(1+\ln t)}\right] (1+ \ln t)  \right)
\end{split}
\]
which proves (\ref{eqn:mero:risk}). The setting of  step sizes  $\eta_t^w$ and $\eta_t^q$ in (\ref{eqn:size:all}) is derived by combining (\ref{eqn:relation:step}) and (\ref{eqn:size:main}).

From (\ref{eqn:error:decompo}) , (\ref{eqn:thm2:2}),  (\ref{lem:1:high}) and  (\ref{lem:2:high}), with probability at least $1-2\delta$, we have
\[
\begin{split}
 & \max_{\q\in \Delta_m}  \phi(\wb_t,\q)- \min_{\w\in \W}  \phi(\w,\qb_t)   \\
\leq & \left( \sum_{j=1}^t \eta_j\right)^{-1} \left( 2+ \frac{5M^2}{2}  \sum_{j=1}^t\eta_j^2 +8 M \sqrt{\sum_{j=1}^t \eta_j^2 \left(\ln \frac{2 t^2}{\delta}\right)}\right.\\
& \left.+\sum_{j=1}^t \eta_j \frac{DG\big[3+\ln j + 16\sqrt{ (1+ \ln j) \ln (2 m j^2/\delta) } \big]}{2(\sqrt{j+1} - 1)} \right) \\
\overset{\text{(\ref{eqn:size:main}),(\ref{eqn:simple:ienqulity})}}{\leq} & \left( \sum_{j=1}^t \frac{1}{\sqrt{j}} \right)^{-1} \left( 2M+ \frac{5 M}{2} \sum_{j=1}^t \frac{1}{j} +8 M \sqrt{ \left(\ln \frac{2 t^2}{\delta}\right) \sum_{j=1}^t \frac{1}{j}}\right.\\
& \left.+ 2DG\left( 3+\ln t + 16\sqrt{(1+ \ln t) \ln \frac{2 m t^2}{\delta} }\right )\sum_{j=1}^t \frac{1}{j} \right)\\
\overset{\text{(\ref{eqn:mr:risk:6})}}{\leq} &  \frac{1}{2 (\sqrt{t+1} - 1)} \left[ M\left(5 +3 \ln t+8 \sqrt{ (1+\ln t) \ln \frac{2 t^2}{\delta}} \right)\right.\\
& \left.+ 2DG\left(3+\ln t + 16\sqrt{(1+ \ln t) \ln \frac{2 m t^2}{\delta} }\right )(1+\ln t) \right]
\end{split}
\]
which proves (\ref{eqn:mero:high}).
\subsection{Proof of Theorem~\ref{thm:3}}
We will follow the analysis of Theorem~\ref{thm:1}, and use a fixed step size to simplify the results.

Let $t=\frac{n_i}{2}$. From (\ref{eqn:mr:risk:4}), we have
\[
\begin{split}
\E \left[R_i(\wb^{(i)}) - R_i(\w_*^{(i)})  \right] \leq & \frac{2D^2 + G^2  \sum_{j=1}^t  (\eta^{(i)})^2 }{2\sum_{j=1}^t \eta^{(i)}} =\frac{D^2}{t \eta^{(i)}} + \frac{G^2 \eta^{(i)}}{2} = DG \sqrt{\frac{2}{t}} = \frac{2DG}{\sqrt{n_i}}
\end{split}
\]
where we set $\eta^{(i)}= \frac{D}{G}\sqrt{\frac{2}{t}} = \frac{2D}{G\sqrt{n_i}}$, which proves  (\ref{eqn:exp:risk:fixed}).

Repeating the proof in Section~\ref{sec:high:pro}, with probability at least $1-\delta/m$, we have
\[
\begin{split}
R_i(\wb^{(i)}) - R_i(\w_*^{(i)})  \leq  & \frac{D^2 + \frac{ G^2 }{2} \sum_{j=1}^t  (\eta^{(i)})^2 + 8 DG \sqrt{ \sum_{j=1}^t (\eta^{(i)})^2 \ln \frac{m}{\delta}} }{\sum_{j=1}^t \eta^{(i)}} \\
=&\frac{D^2}{t \eta^{(i)}} + \frac{G^2 \eta^{(i)}}{2} + \frac{8DG \sqrt{\ln \frac{m}{\delta}}}{\sqrt{t}}\\
= & DG \sqrt{\frac{2}{t}}\left(1+4 \sqrt{2\ln \frac{m}{\delta}}\right)= \frac{2DG}{\sqrt{n_i}} \left(1+4 \sqrt{2\ln \frac{m}{\delta}}\right)
\end{split}
\]
for any  $i \in [m]$. We obtain (\ref{eqn:exp:high:fixed}) by taking the union bound over all $i\in[m]$.

\subsection{Proof of Theorem~\ref{thm:6}}
Following the analysis of \citet[Theorem 6]{SA:GDRO}, we bound the optimization error of $(\wb,\qb)$ for (\ref{eqn:convex:concave:weight:appro}) below.
\begin{thm} \label{thm:5} Under the condition of Theorem~\ref{thm:6},  with probability at least  $1-\delta$, we have
\[
\epsilon_{\widehat{\varphi}}(\wb, \qb) \leq  \frac{14 \Lt}{n_m} + \sqrt{\frac{\sigma^2}{n_m}}  \left(\frac{28}{\sqrt{3}} + 7\sqrt{6 \log \frac{2}{\delta}} + \frac{28\sqrt{2}}{n_m} \log \frac{2}{\delta} \right) .
\]
\end{thm}

Then, we make use of Lemma~\ref{lem:diff:opt} to bound the optimization error of  $(\wb,\qb)$ for problem (\ref{eqn:convex:concave:weight}).  From Theorem~\ref{thm:3} and Lemma~\ref{lem:diff:opt}, with probability at least  $1-\delta$, we have
\[
\begin{split}
\epsilon_{\varphi}(\wb, \qb)\leq & \epsilon_{\widehat{\varphi}}(\wb, \qb)+ 2\max_{i \in [m]} \left\{p_i \big[ R_i(\wb^{(i)})- R_i^*\big]\right\} \\
\leq & \epsilon_{\widehat{\varphi}}(\wb, \qb)+ 2\max_{i \in [m]} \left\{p_i \ \frac{2DG}{\sqrt{n_i}} \left(1+4 \sqrt{2\ln \frac{m}{\delta}}\right) \right\}\\
=& \epsilon_{\widehat{\varphi}}(\wb, \qb)+ 4DG \left(1+4 \sqrt{2\ln \frac{m}{\delta}}\right)\max_{i \in [m]} \frac{p_i }{\sqrt{n_i}} .
\end{split}
\]
Combining with Theorem~\ref{thm:5}, with probability at least $1-2\delta$, we have
\begin{equation} \label{eqn:final:error}
\begin{split}
&\epsilon_{\varphi}(\wb, \qb) \\
 \leq & \frac{14 \Lt}{n_m} + \sqrt{\frac{\sigma^2}{n_m}}  \left(\frac{28}{\sqrt{3}} + 7\sqrt{6 \log \frac{2}{\delta}} + \frac{28\sqrt{2}}{n_m} \log \frac{2}{\delta} \right) + 4DG \left(1+4 \sqrt{2\ln \frac{m}{\delta}}\right) \max_{i \in [m]} \frac{p_i }{\sqrt{n_i}}.
\end{split}
\end{equation}

Next, we bound the excess risk of $\wb$ on every distribution. To this end, we have
\[
\begin{split}
& \max_{i\in [m]} \left\{ p_i \big [R_i(\wb)- R_i^*\big] \right\} - \min_{\w \in \W} \max_{\q \in \Delta_m} \varphi(\w,\q) \\
=& \max_{\q\in \Delta_m}  \left\{\sum_{i=1}^m q_i p_i \big [R_i(\wb)- R_i^*\big] \right\} - \min_{\w \in \W} \max_{\q \in \Delta_m} \varphi(\w,\q) \\
\leq & \max_{\q\in \Delta_m}  \left\{\sum_{i=1}^m q_i p_i \big [R_i(\wb)- R_i^*\big] \right\} - \min_{\w\in \W}  \left\{ \sum_{i=1}^m \bar{q}_i p_i \big [R_i(\w)- R_i^*\big] \right\}\\
= &  \max_{\q\in \Delta_m}  \varphi(\wb,\q)- \min_{\w\in \W}  \varphi(\w,\qb) = \epsilon_{\varphi}(\wb, \qb).
\end{split}
\]
Thus, for every distribution $\P_i$, the excess risk can be bounded in the following way:
\[
\begin{split}
&R_i(\wb)- R_i^* \\
\leq & \frac{1}{p_i} \min_{\w \in \W} \max_{\q \in \Delta_m} \varphi(\w,\q) + \frac{1}{p_i} \epsilon_{\varphi}(\wb, \qb)  \\
\overset{\text{(\ref{eqn:final:error})}}{\leq} & \frac{1}{p_i} \min_{\w \in \W} \max_{\q \in \Delta_m} \varphi(\w,\q) \\
& + \frac{1}{p_i} \left[  \frac{14 \Lt}{n_m} + \sqrt{\frac{\sigma^2}{n_m}}  \left(\frac{28}{\sqrt{3}} + 7\sqrt{6 \log \frac{2}{\delta}} + \frac{28\sqrt{2}}{n_m} \log \frac{2}{\delta} \right) + 4DG \left(1+4 \sqrt{2\ln \frac{m}{\delta}}\right) \max_{i \in [m]} \frac{p_i }{\sqrt{n_i}}\right].
\end{split}
\]
which proves (\ref{eqn:final:risk}).

\subsection{Proof of Corollary~\ref{cor:1}}
From (72) and (73) of \citet{SA:GDRO}, we have
\begin{equation} \label{eqn:pi:ine:1}
\frac{1}{p_i} \frac{\Lt}{n_m}  = O\left( \left(\frac{1}{n_m} + \frac{1}{\sqrt{n_i}}\right)\sqrt{\ln m}\right) \textrm{ and } \frac{1}{p_i} \sqrt{\frac{\sigma^2}{n_m}} =  O\left(\left( \frac{1}{n_m} + \frac{1}{\sqrt{n_i}} \right) \sqrt{\kappa+\ln^2 m}\right).
\end{equation}
Furthermore,
\begin{equation} \label{eqn:pi:ine:2}
\begin{split}
& \frac{r_{\max}}{p_i}=\frac{1}{p_i} \max_{i \in [m]} \frac{p_i }{\sqrt{n_i}} =\frac{1}{p_i} \max_{i \in [m]} \left(\frac{1/\sqrt{n_m} + 1}{1/\sqrt{n_m} + \sqrt{n_m/n_i}} \frac{1 }{\sqrt{n_i}}  \right) \\
\leq & \frac{1}{p_i} \max_{i \in [m]} \left( \left( \frac{1}{\sqrt{n_m}} + 1\right) \sqrt{\frac{n_i}{n_m}} \frac{1 }{\sqrt{n_i}}  \right)  =\frac{1}{p_i} \left( \frac{1}{\sqrt{n_m}} + 1\right)  \frac{1}{\sqrt{n_m}} \\
= &  \frac{1/\sqrt{n_m} + \sqrt{n_m/n_i}}{1/\sqrt{n_m} + 1} \left( \frac{1}{\sqrt{n_m}} + 1\right)  \frac{1}{\sqrt{n_m}} = \frac{1}{n_m} + \frac{1}{\sqrt{n_i}}.
\end{split}
\end{equation}
Combining (\ref{eqn:final:risk}), (\ref{eqn:pi:ine:1}), and (\ref{eqn:pi:ine:2}), we have
\[
R_i(\wb)- R_i^*  = \frac{1}{p_i} \min_{\w \in \W} \max_{\q \in \Delta_m} \varphi(\w,\q) +  O\left(\left( \frac{1}{n_m} + \frac{1}{\sqrt{n_i}} \right) \sqrt{\kappa+\ln^2 m}\right).
\]

\subsection{Proof of Corollary~\ref{cor:2} }
Under the assumption of this corollary, we have
\[
\begin{split}
0 \leq & p_{\varphi}^* = \min_{\w \in \W} \max_{\q \in \Delta_m} \varphi(\w,\q)  \\
\leq & \max_{\q \in \Delta_m} \varphi(\w_*,\q) =\max_{\q \in \Delta_m} \left\{ \sum_{i=1}^m q_i p_i \big [R_i(\w_*)- R_i^*\big] \right\} = 0.
\end{split}
\]

\subsection{Proof of Theorem~\ref{thm:5}}
Our analysis is similar to that of Theorem 6 of \citet{SA:GDRO}. For brevity, we will only highlight the differences.

We first introduce the monotone operator $F(\w,\q)$ associated with (\ref{eqn:convex:concave:weight:appro}):
\[
\begin{split}
F(\w,\q)=&[\nabla_\w \widehat{\varphi}(\w,\q); -\nabla_\q \widehat{\varphi}(\w,\q) ]\\
=&\left[ \sum_{i=1}^m q_i p_i \nabla R_i(\w)  ; -\Big[p_1 [R_1(\w)-R_1(\wb^{(1)})], \ldots, p_m \big[R_m(\w)-R_m(\wb^{(m)})\big] \Big]^\top \right] .
\end{split}
\]
It is slightly different from the monotone operator defined by \citet[(68)]{SA:GDRO}, attributed to the inclusion of $R_1(\wb^{(1)}),\ldots,R_m(\wb^{(m)})$. However, these additional terms do not alter the continuity of $F(\w,\q)$. In particular, Lemma 3 of \citet{SA:GDRO} remains applicable, leading to the following lemma.
\begin{lem} \label{lem:monotone} For the monotone operator $F(\w,\q)$, we have
\[
\| F(\w,\q) -F(\w',\q')\|_* \leq \Lt \big\|(\w-\w',\q-\q') \big\|
\]
where $\Lt$ is defined in (\ref{eqn:mirror:parameters}).
\end{lem}

Then, we investigate the stochastic oracle in Algorithm~\ref{alg:2}:
\[
\g(\w,\q) = \left[\g_w(\w,\q); -\g_q(\w,\q)\right]
\]
where
\[
\begin{split}
\g_w(\w,\q) &=\sum_{i=1}^m q_{i} p_i \left(\frac{n_m}{n_i} \sum_{j=1}^{n_i/n_m} \nabla \ell(\w;\z^{(i,j)}) \right),\\
\g_q(\w,\q) & =\left[p_1 \frac{n_m}{n_1}  \sum_{j=1}^{n_1/n_m} \big[\ell(\w;\z^{(1,j)}) - \ell(\wb^{(1)};\z^{(1,j)})\big], \ldots, p_m \big[\ell(\w;\z^{(m)})  - \ell(\wb^{(m)};\z^{(m)}) \big]\right]^\top
\end{split}
\]
and $\z^{(i,j)}$ is the $j$-th sample drawn from distribution $\P_i$. Again, $\g(\w,\q)$ is  different from that of \citet[(70)]{SA:GDRO}, because of the additional terms. However, it is easy to verify that the variance only changes by a constant factor, and Lemma 4 of \citet{SA:GDRO} still holds with a different constant.
\begin{lem} \label{lem:variance} For the stochastic oracle $\g(\w,\q)$, we have
\[
 \E \left[ \exp \left( \frac{\| F(\w,\q) -\g(\w,\q)\|_*^2}{\sigma^2 } \right ) \right] \leq  2
\]
where $\sigma^2$ is defined in (\ref{eqn:mirror:parameters}).
\end{lem}
The final difference lies in the number of iterations for SMPA, which is $n_m/4$ in our Algorithm~\ref{alg:2}, and $n_m/2$ in the work of \citet{SA:GDRO}.

From Corollary 1 of \citet{Nemirovski:SMP}, by setting
\[
\eta=\min \left( \frac{1}{\sqrt{3} \Lt}, 2\sqrt{\frac{2}{7 \sigma^2 n_m}} \right)
\]
we have
\[
\Pr\left[\epsilon_{\varphi}(\wb, \qb) \geq \frac{14 \Lt}{n_m} + 28\sqrt{\frac{\sigma^2}{ 3n_m}} + 7\Lambda \sqrt{\frac{2\sigma^2}{n_m}} \right] \leq \exp\left(-\frac{\Lambda^2}{3} \right) + \exp\left(-\frac{\Lambda n_m}{4}\right)
\]
for all $\Lambda >0$. Choosing $\Lambda$ such that $\exp(-\Lambda^2/3) \leq \delta/2$ and $\exp(-\Lambda n_m/4) \leq \delta/2$, we have with probability at least $1-\delta$
\[
\epsilon_{\varphi}(\wb, \qb) \leq  \frac{14 \Lt}{n_m} + 28\sqrt{\frac{\sigma^2}{3n_m}} +  7\left(\sqrt{3 \log \frac{2}{\delta}} + \frac{4}{n_m} \log \frac{2}{\delta}\right) \sqrt{\frac{2\sigma^2}{n_m}}.
\]
\section{Experiments}
In this section, we conduct empirical studies to evaluate our proposed algorithms.

\subsection{Datasets and Experimental Settings} \label{sec:setting}
Following previous work \citep{NIPS2016_4588e674,DRO:Online:Game}, we employ both synthetic and real-world datasets.

For the synthetic dataset, we construct $m=6$ distributions, each of which is associated with a  true classifier $\w_i^*\in\mathbb{R}^{1000}$. The selection process is as follows: we initially choose an arbitrary $\w_0^*$ on the unit sphere. Subsequently, we randomly pick $m$ points on a sphere with radius $d$, centered at $\w_0^*$. These points are then projected onto the unit sphere to form the set $\{\w_i^*\}_{i\in[m]}$. We set $d=0.2$ to keep the classifiers $\{\w_i^*\}_{i\in[m]}$ close, thereby emphasizing the optimization challenges due to the varying noise across the distributions.  For each distribution $i\in[m]$, a sample $(\x,y)$ is generated where $\x$ follows a standard normal distribution $\mathcal{N}(0,I)$, and $y=\sgn(\x^{\top}\w_i^*)$  or its inverse, each with respective probabilities $p_i=1-0.05 \times i$ and $1-p_i$. 

We additionally utilize the Adult dataset \citep{misc_adult_2}, which encompasses a variety of attributes, including age, gender, race, and educational background, for a total of $48842$ individuals. The samples are classified into $6$ distinct groups, based on a combination of race \{black, white, others\} and gender \{female, male\}. By employing one-hot encoding for $12$ attributes, $103$-dimensional feature vectors are generated to classify whether the income surpasses $\$5000$.

We designate the logistic loss as our loss function, and utilize various methods to train a linear model. In assessing their performance, it's essential to determine the model's risk for each distribution. To approximate the risk value, we will draw a specific number of samples and  use the average risk calculated from these samples as an estimation. The minimal risk $R_i^*$ for each distribution is estimated in a similar way: initially, a model is trained to minimize the empirical risk using a large set of samples; subsequently, the risk is calculated based on a freshly drawn sample set.

\subsection{Experiments on Balanced Data} 
In our experiments with the synthetic dataset, samples will be dynamically generated in real time, adhering to the generation protocol in Section~\ref{sec:setting}. Regarding the Adult dataset, we designate the distribution  $\P_i$  to represent the uniform distribution over the data within the $i$-th group, and thus the sample generalization process simplifies to randomly selecting samples from each group with replacement.

\begin{figure}[t]
\begin{center}
\subfigure[The synthetic dataset]{
    \label{fig:1:a} 
    \includegraphics[width=0.45\textwidth]{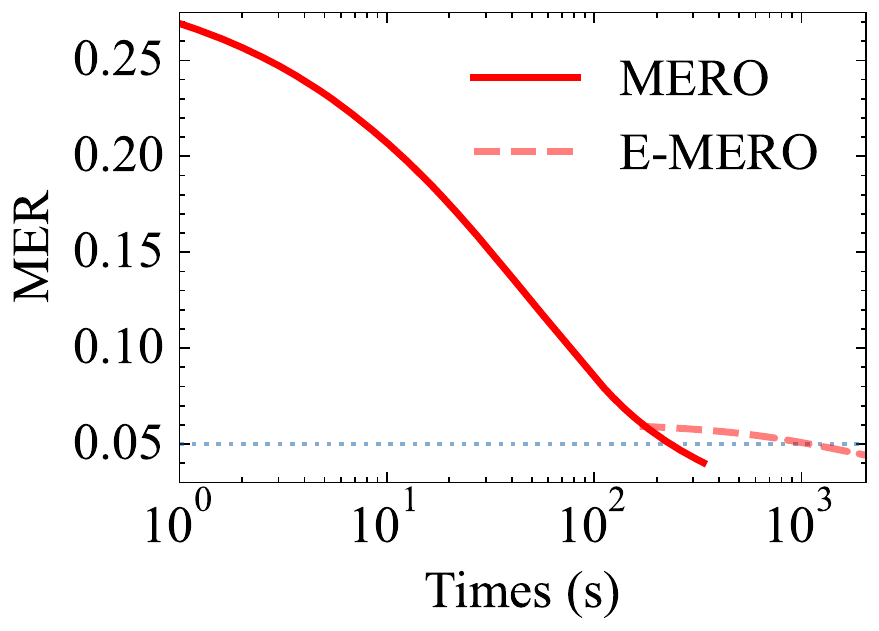}}%
\subfigure[The Adult dataset]{
    \label{fig:1:b} 
    \includegraphics[width=0.45\textwidth]{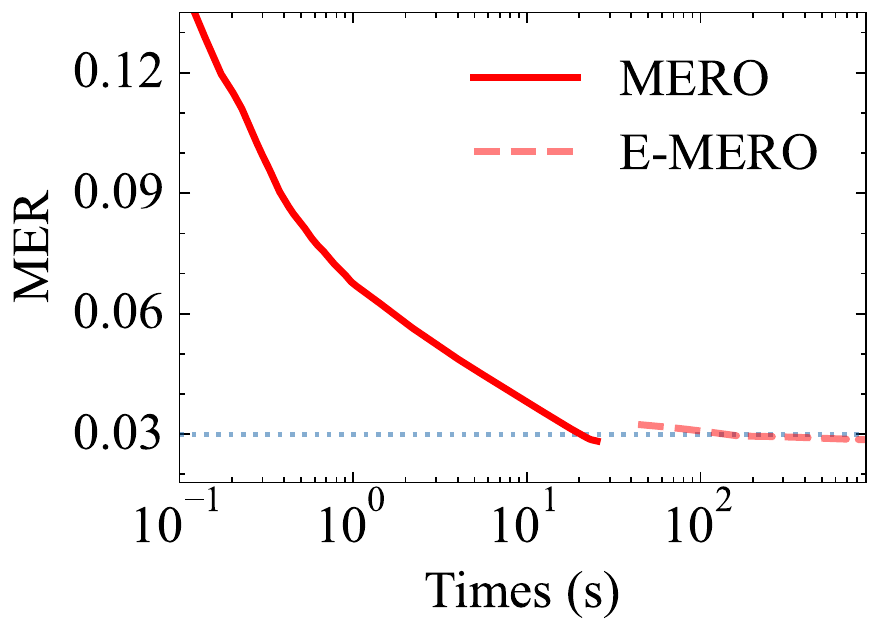}}%
	\caption{The maximal excess risk (MER) versus the running time.}
  \label{fig:1}
\end{center}
\begin{center}
\tabcaption{Running times of MERO and E-MERO.}
		\label{tab:1}
		\begin{tabular}{c|c|c|c}
			\hline
			Dataset & Algorithm &  MER Value & Times (s) \\
			\hline
			\multirow{2}{*}{Synthetic} & MERO & 0.05 & 255.0 \\
			\cline{2-4}
			\multirow{2}{*}{} & E-MERO & 0.05 & 1470.5 \\
			\hline
			\multirow{2}{*}{Adult} & MERO & 0.03 & 20.4 \\
			\cline{2-4}
			\multirow{2}{*}{}& E-MERO & 0.03 &  141.6 \\
			\hline
		\end{tabular}
\end{center}
\begin{center}
\subfigure[The synthetic dataset]{
    \label{fig:2:a} 
    \includegraphics[width=0.45\textwidth]{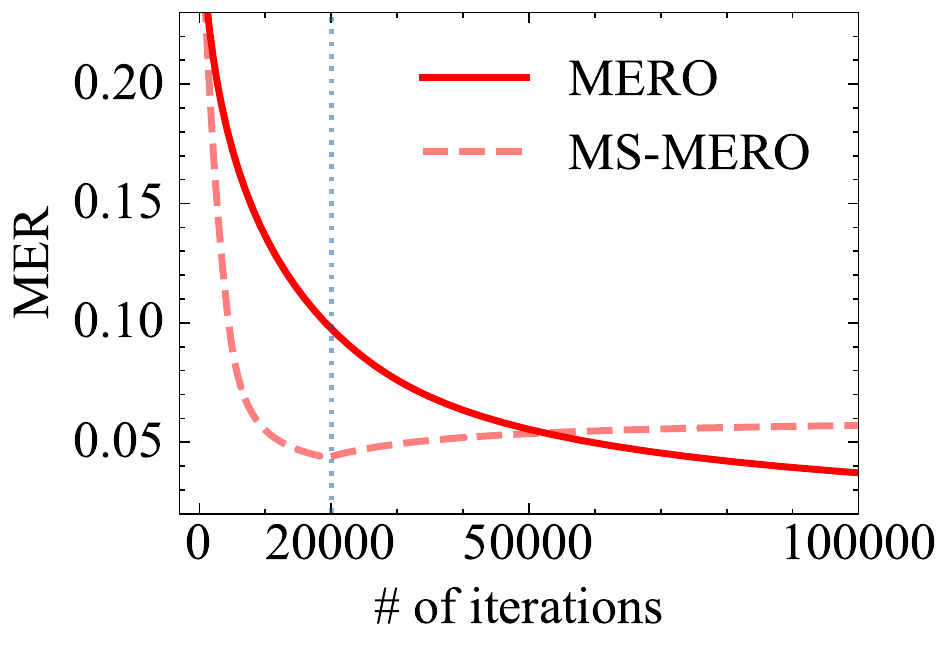}}%
\subfigure[The Adult dataset]{
    \label{fig:2:b} 
    \includegraphics[width=0.45\textwidth]{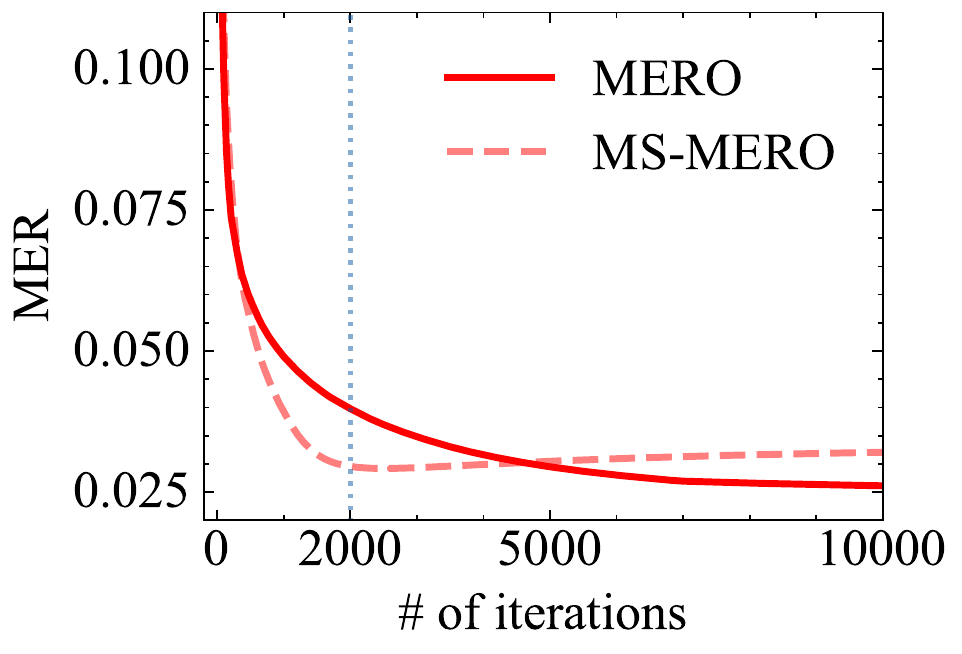}}%
	\caption{The maximal excess risk (MER) versus the number of iterations.}
  \label{fig:2}
\end{center}
\end{figure}

First, we demonstrate the efficiency of our anytime stochastic approximation approach in Algorithm~\ref{alg:1}, referred to as MERO, through a comparative analysis with the optimization procedure of \citet{Regret:RML:DS} for empirical MERO, denominated as E-MERO. Recall that Algorithm~\ref{alg:1} runs for $T$ rounds, making use of a cumulative total of $mT$ samples. For a fair comparison, we set the number of samples from each distribution in E-MERO to be $T$. We assign the value of $T$ as $10^5$ for the synthetic dataset and $10^4$ for the Adult dataset. In Fig.~\ref{fig:1}, we depict the relationship between the maximal excess risk (MER), i.e.,  $\max_{i\in[m]}[R_i(\w)-R_i^*]$, and the running time.  The lack of the E-MERO curve at the beginning is attributed to its initialization phase, a period focused on minimizing the empirical risk for each distribution. It's evident that MERO achieves convergence significantly quicker than E-MERO for both datasets, highlighting the computational efficiency of stochastic approximation. To be more clear, we present in Table~\ref{tab:1} the time required for MERO and E-MERO to attain a specified MER target---$0.05$ for the synthetic dataset and $0.03$ for the Adult dataset. The data reveals that MERO outpaces E-MERO, being $5.7$ times faster and $6.9$ times faster on the synthetic and Adult datasets, respectively. Additionally, it's important to note that MERO is more memory-efficient, as it eliminates the need to store training data.

Second, we illustrate the benefit of the anytime capability of Algorithm~\ref{alg:1} by conducting a comparison with the multi-stage approach detailed in Section~\ref{sec:multi:stage}, denoted by MS-MERO. In the case of MS-MERO, we assign a preset value of $T$ as $2\times10^4$ for the synthetic dataset and $2\times10^3$ for the Adult dataset. However, we continue the execution of the 3rd stage even when the count of iterations goes beyond the set value of $T$. In Fig.~\ref{fig:2}, the graph displays  the progression of MER in relation to the number of iterations. We exclude the initial two stages of MS-MERO from the illustration, because it only produces a model during the 3rd stage. In the beginning, MS-MERO rapidly reduces the MER, but as the actual number of executed rounds exceeds the predetermined limit, its parameter settings are not optimal. Consequently, we observe a stagnation or even a slight increase in the MER. In contrast,  MERO  demonstrates a consistent decrease in MER and ultimately outperforms MS-MERO. Such results  underscore the importance of anytime algorithms, especially when the total number of iterations is unknown.

\begin{figure}[t]
\begin{center}
\subfigure[Risk on $\P_1$]{
    \label{fig:3:a} 
    \includegraphics[width=0.33\textwidth]{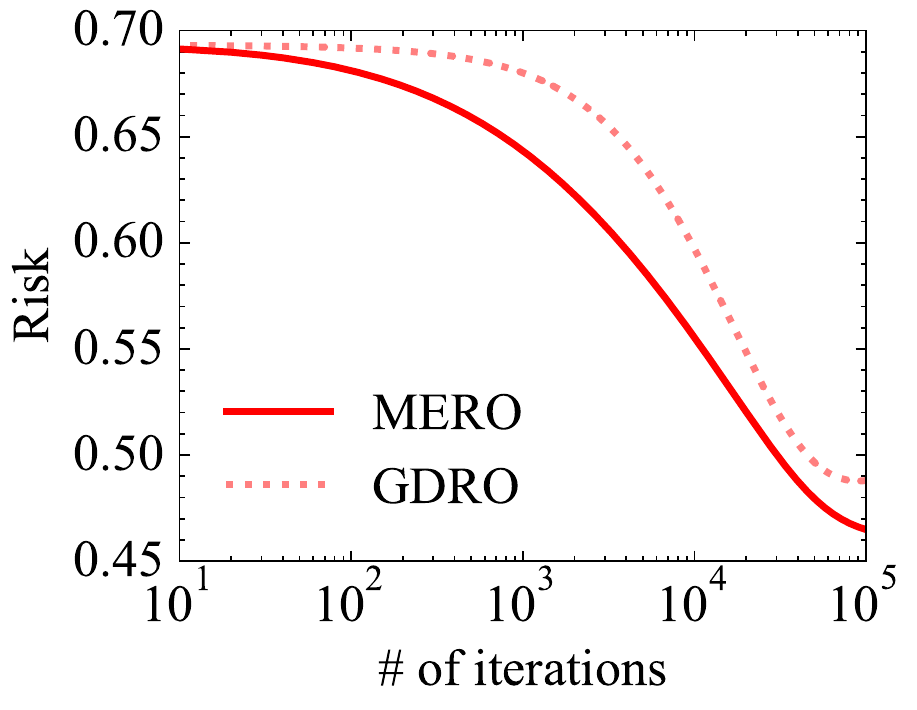}}%
\subfigure[Risk on $\P_{2}$]{
    \label{fig:3:b} 
    \includegraphics[width=0.33\textwidth]{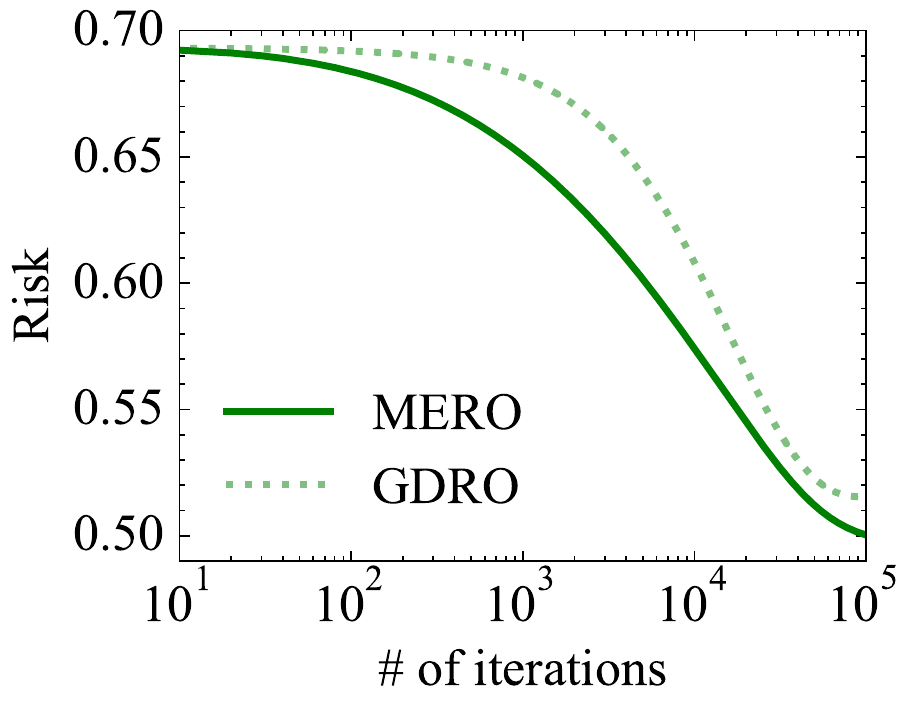}}%
    \subfigure[Risk on $\P_{3}$]{
    \label{fig:3:c} 
    \includegraphics[width=0.33\textwidth]{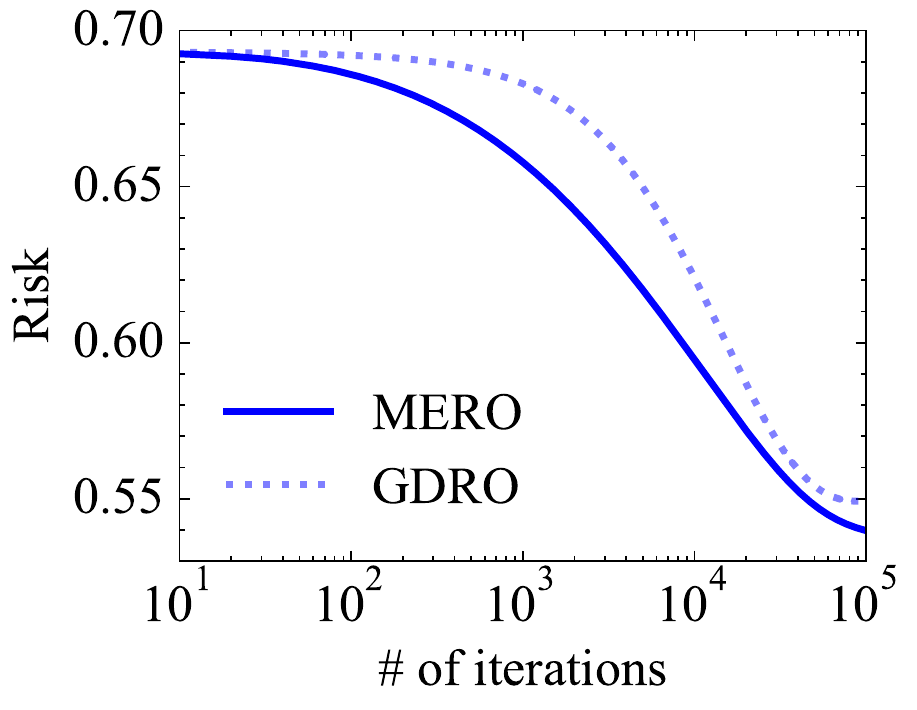}}\\%
        \subfigure[Risk on $\P_{4}$]{
    \label{fig:3:d} 
    \includegraphics[width=0.33\textwidth]{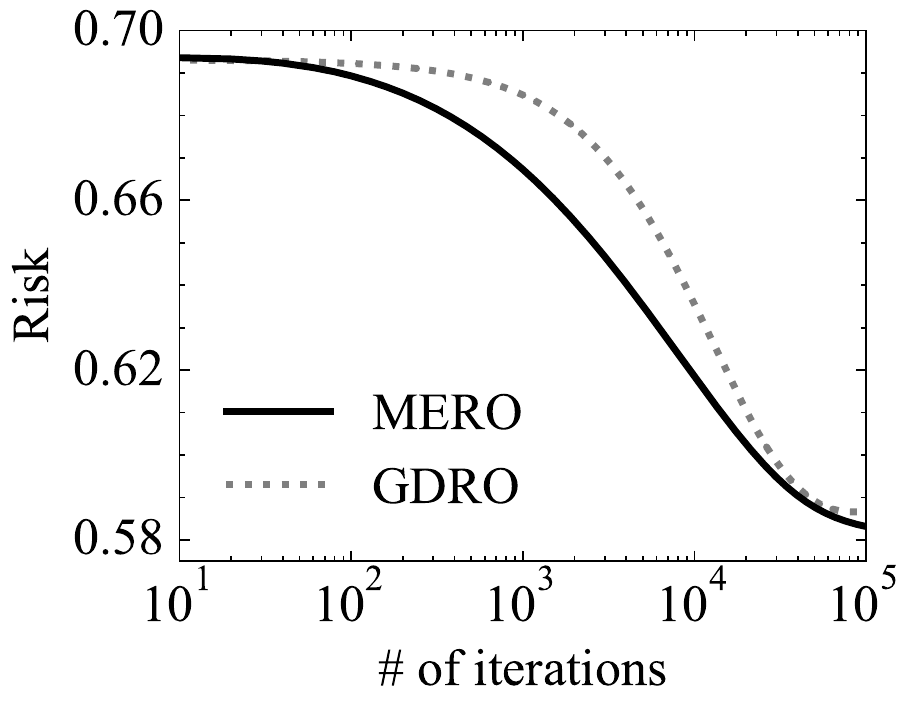}}%
    \subfigure[Risk on $\P_{5}$]{
    \label{fig:3:e} 
    \includegraphics[width=0.33\textwidth]{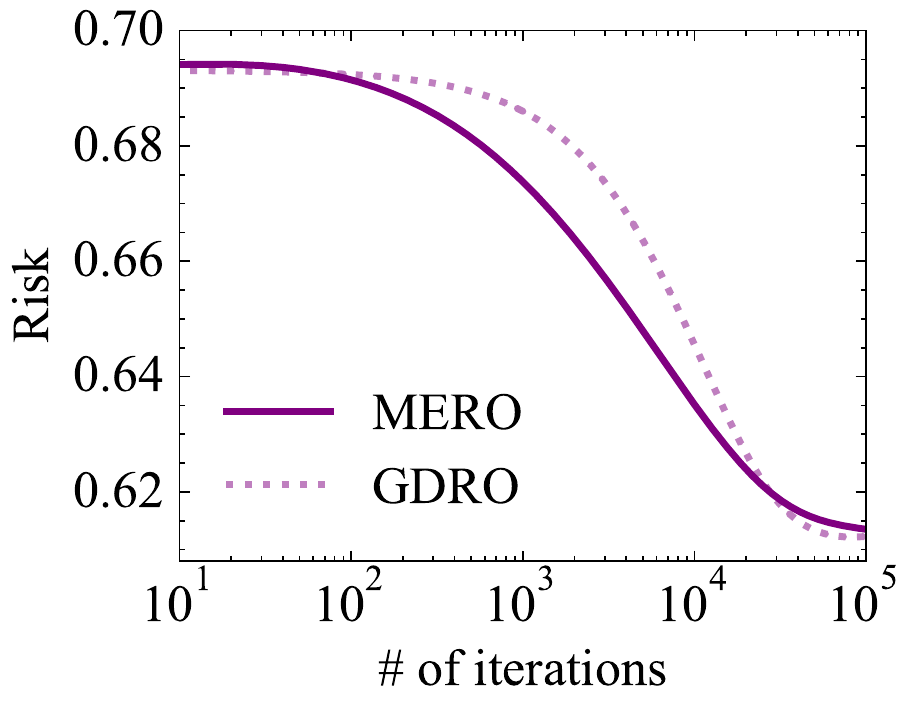}}%
    \subfigure[Risk on $\P_{6}$]{
    \label{fig:3:f} 
    \includegraphics[width=0.33\textwidth]{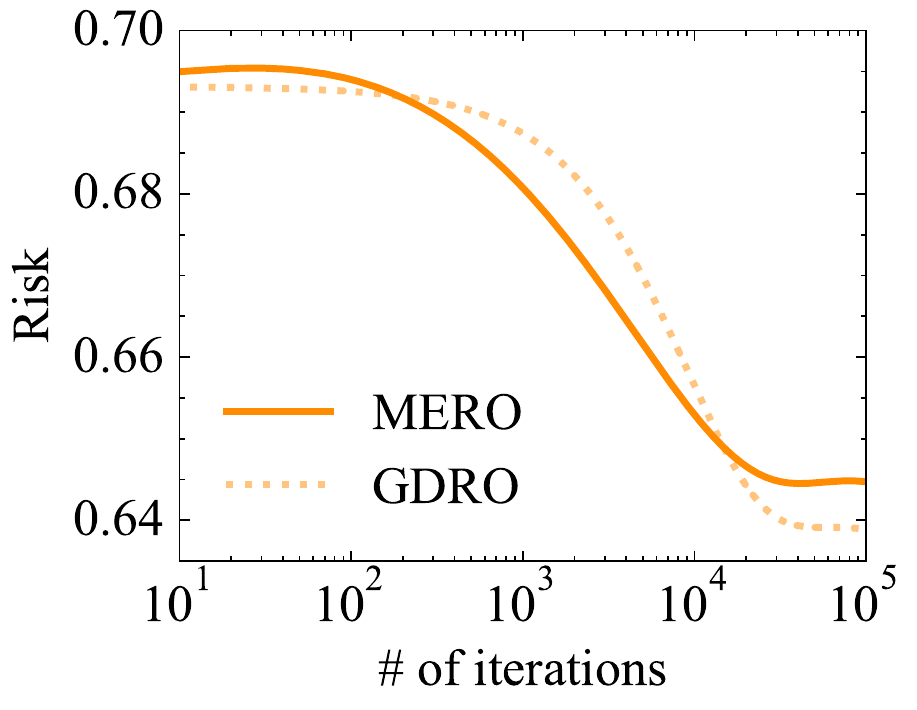}}\\%
  \caption{Individual risk versus the number of iterations on the synthetic dataset.}
  \label{fig:3} \vspace{1ex}%
\subfigure[Risk on $\P_1$]{
    \label{fig:4:a} 
    \includegraphics[width=0.33\textwidth]{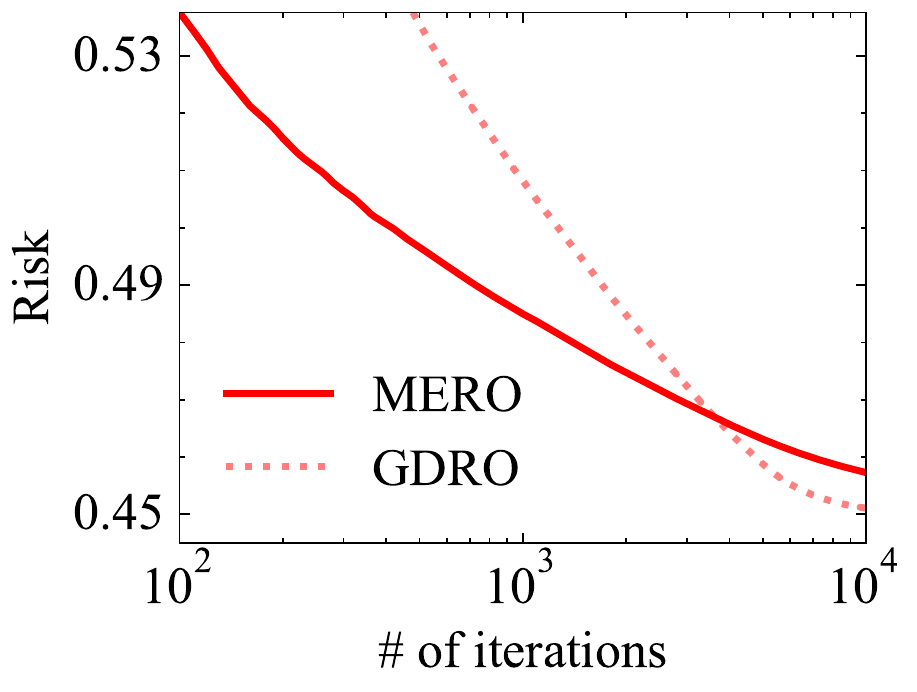}}%
\subfigure[Risk on $\P_{2}$]{
    \label{fig:4:b} 
    \includegraphics[width=0.33\textwidth]{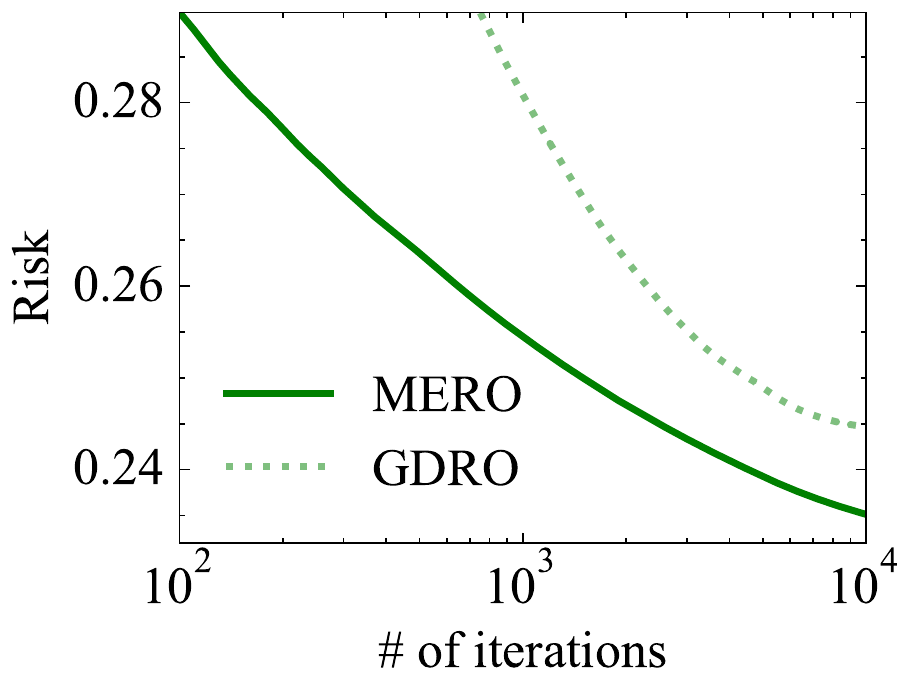}}%
    \subfigure[Risk on $\P_{3}$]{
    \label{fig:4:c} 
    \includegraphics[width=0.33\textwidth]{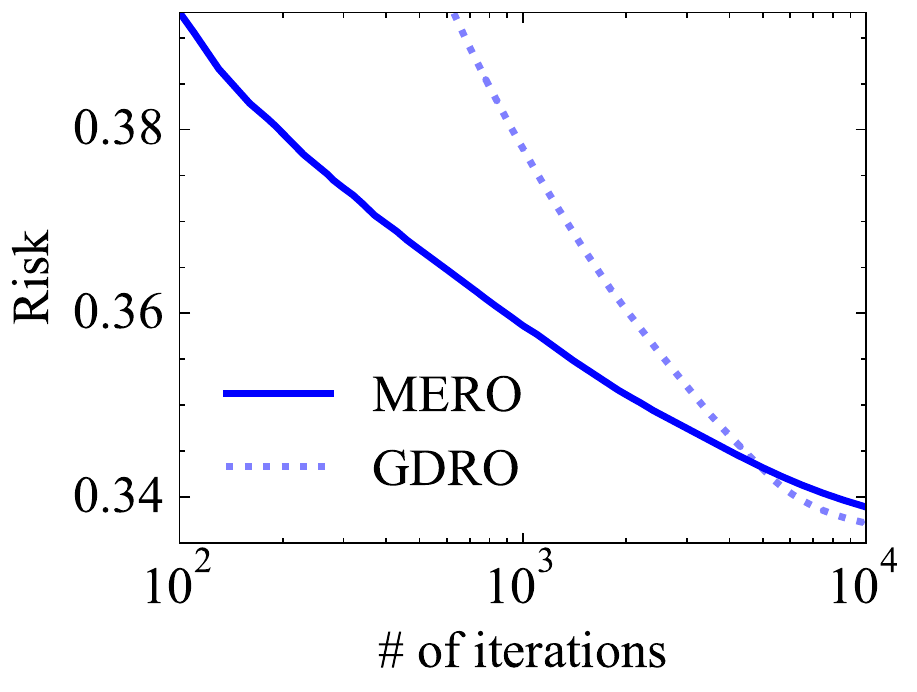}}\\%
        \subfigure[Risk on $\P_{4}$]{
    \label{fig:4:d} 
    \includegraphics[width=0.33\textwidth]{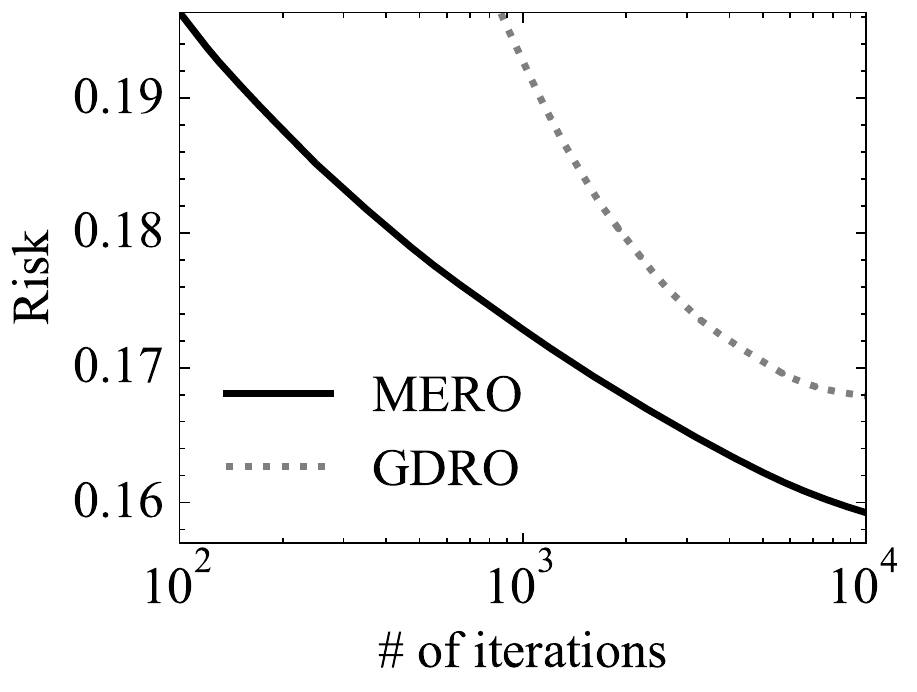}}%
    \subfigure[Risk on $\P_{5}$]{
    \label{fig:4:e} 
    \includegraphics[width=0.33\textwidth]{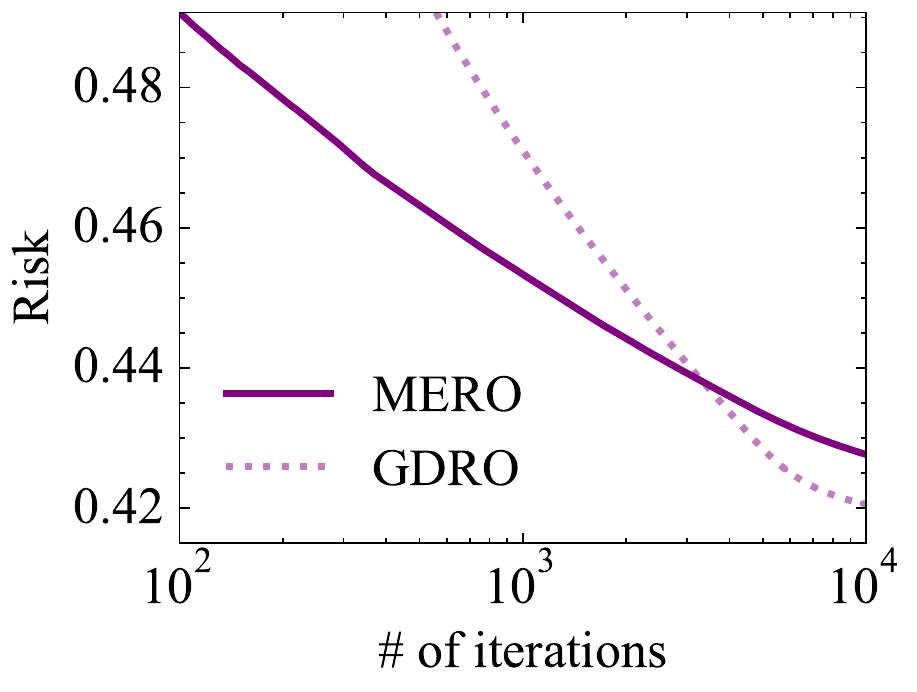}}%
    \subfigure[Risk on $\P_{6}$]{
    \label{fig:4:f} 
    \includegraphics[width=0.33\textwidth]{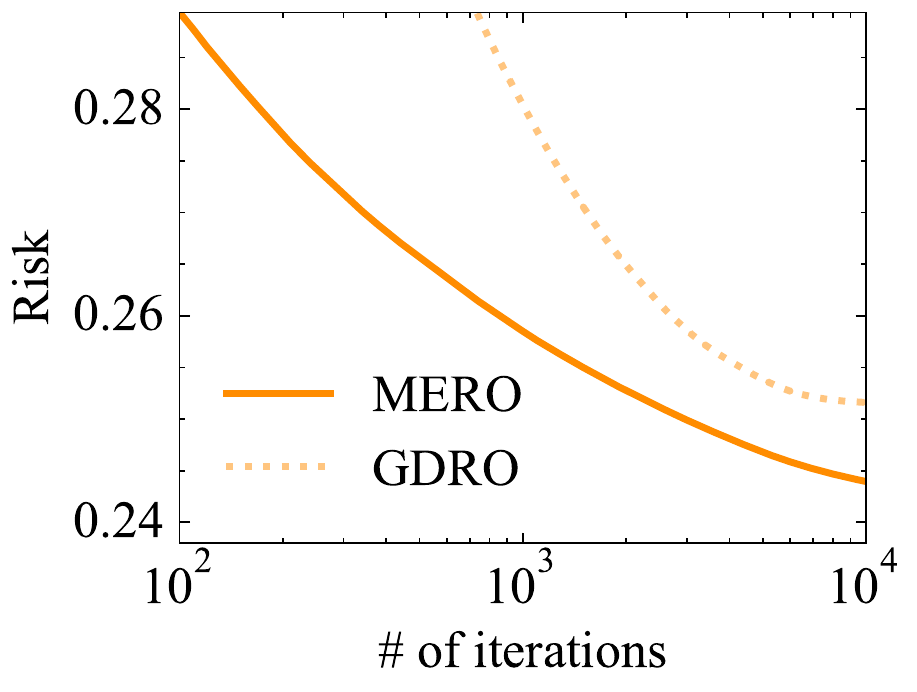}}\\%
  \caption{Individual risk versus the number of iterations on the Adult dataset.}
  \label{fig:4}
\end{center}
\end{figure}

Third, we contrast MERO with GDRO to analyze their unique characteristics.  Notice that MERO and GDRO adopt different objectives, and therefore neither holds absolute superiority  \citep[Proposition 1 and Example 1]{Regret:RML:DS}.  While the MER of MERO is always better than that of GDRO, it does not necessarily mean that  the model yielded by MERO is universally preferable across all distributions.  To illustrate this point, we compare our Algorithm~\ref{alg:1} with the stochastic approximation algorithm of \citet[Algorithm 1]{SA:GDRO} designed for GDRO, which is also referred to as GDRO for convenience. We assess the risks associated with each distribution for both MERO and GDRO, and in the initial comparison, we set the $x$-axis as the number of iterations to ensure that both algorithms consume the same number of samples. Experimental results on the synthetic dataset are presented in Fig.~\ref{fig:3}. As can be seen, GDRO exhibits strong performance with distributions $\P_5$ and $\P_6$, while MERO demonstrates superior results with the remaining $4$ distributions. This pattern is as expected, since GDRO targets the raw risk, and the last two distributions are characterized by the high level of noise, hence the large risk. Consequently, GDRO tends to concentrate its efforts on these distributions, achieving lower risks for them. By contrast, MERO effectively mitigates the impact of noise, achieving a more balanced performance across various distributions.  Experimental results on the Adult dataset, showcased in Fig.~\ref{fig:4}, lead us to similar conclusions: GDRO exhibits slightly better performance on distributions ($\P_1$, $\P_3$ and $\P_5$) with large risk.

\begin{figure}[t]
\begin{center}
\subfigure[Risk on $\P_1$]{
    \label{fig:5:a} 
    \includegraphics[width=0.33\textwidth]{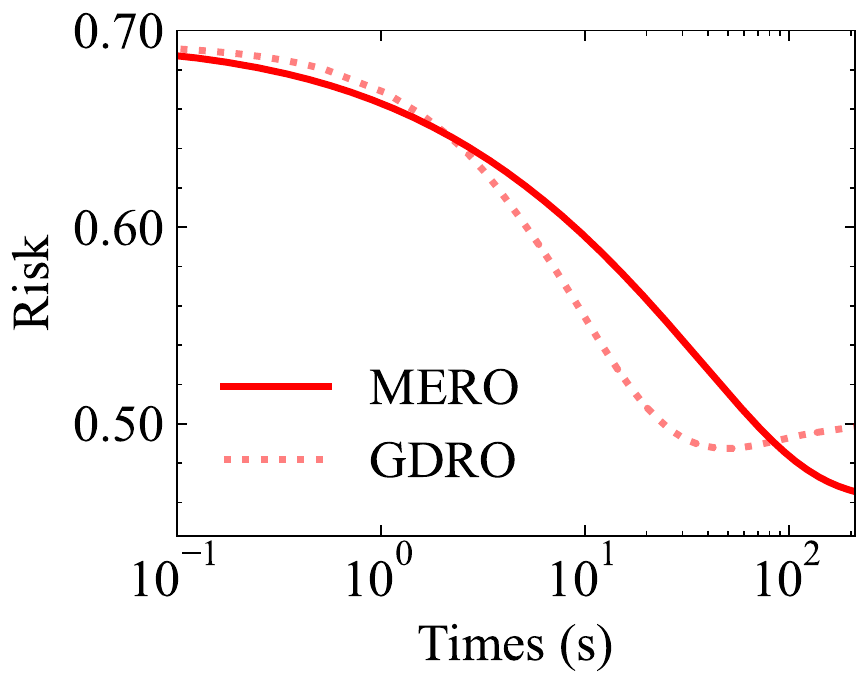}}%
\subfigure[Risk on $\P_{2}$]{
    \label{fig:5:b} 
    \includegraphics[width=0.33\textwidth]{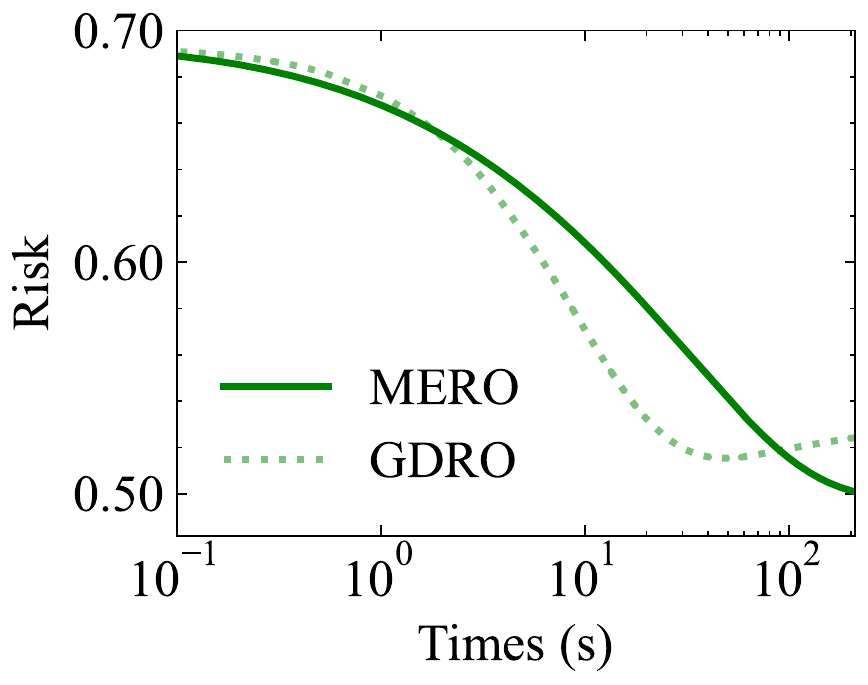}}%
    \subfigure[Risk on $\P_{3}$]{
    \label{fig:5:c} 
    \includegraphics[width=0.33\textwidth]{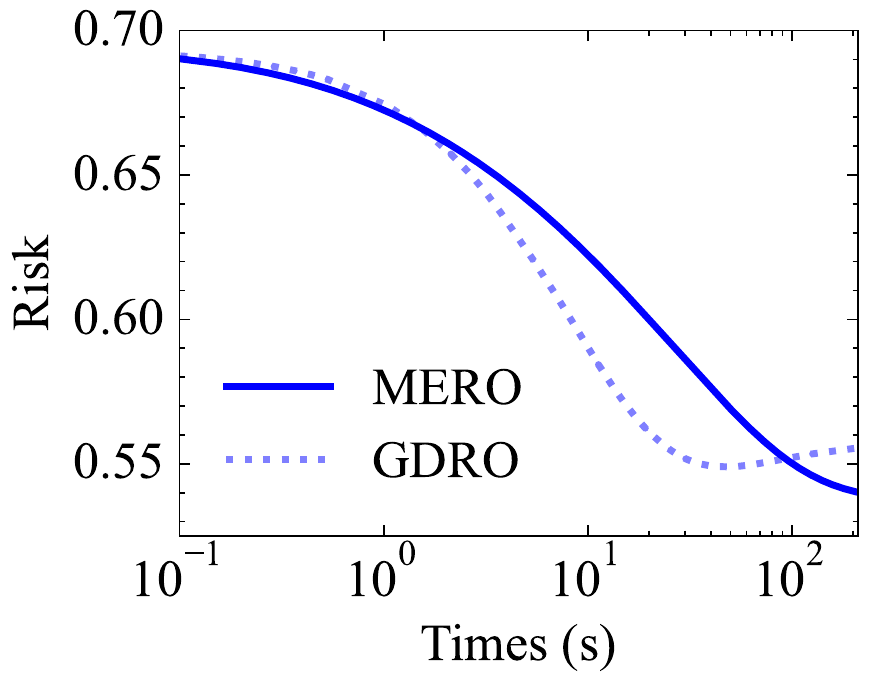}}\\%
        \subfigure[Risk on $\P_{4}$]{
    \label{fig:5:d} 
    \includegraphics[width=0.33\textwidth]{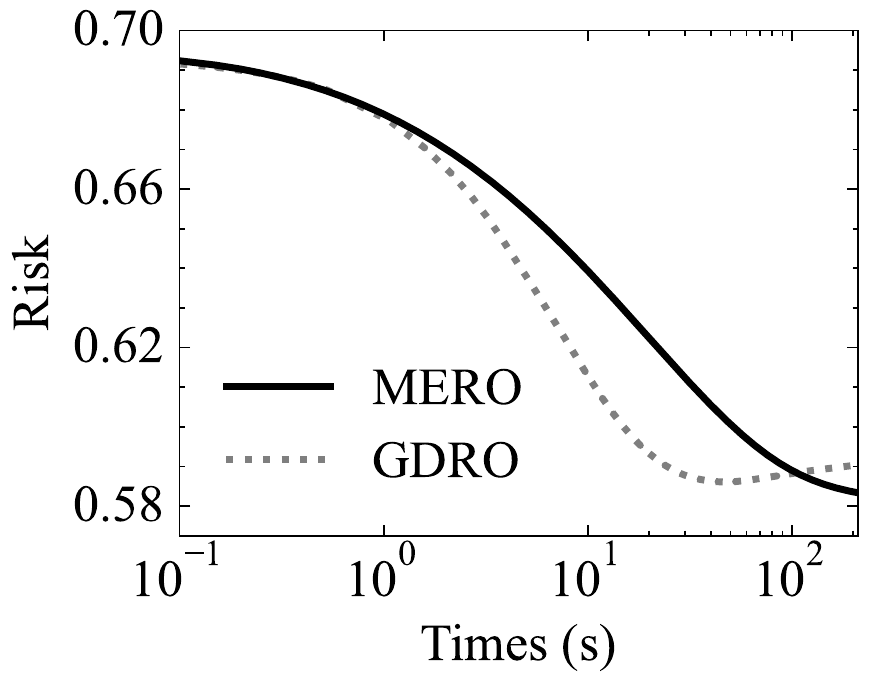}}%
    \subfigure[Risk on $\P_{5}$]{
    \label{fig:5:e} 
    \includegraphics[width=0.33\textwidth]{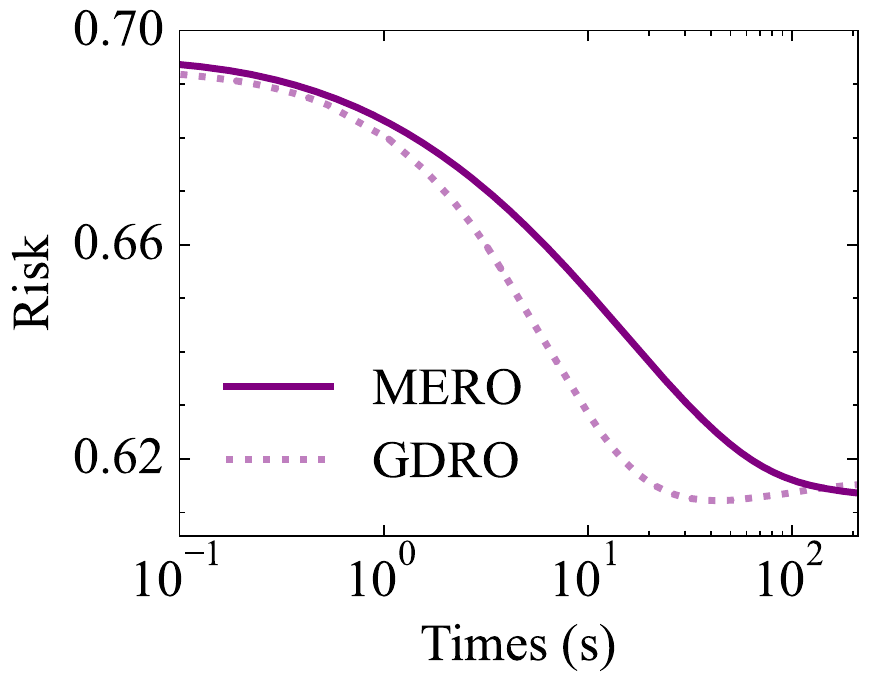}}%
    \subfigure[Risk on $\P_{6}$]{
    \label{fig:5:f} 
    \includegraphics[width=0.33\textwidth]{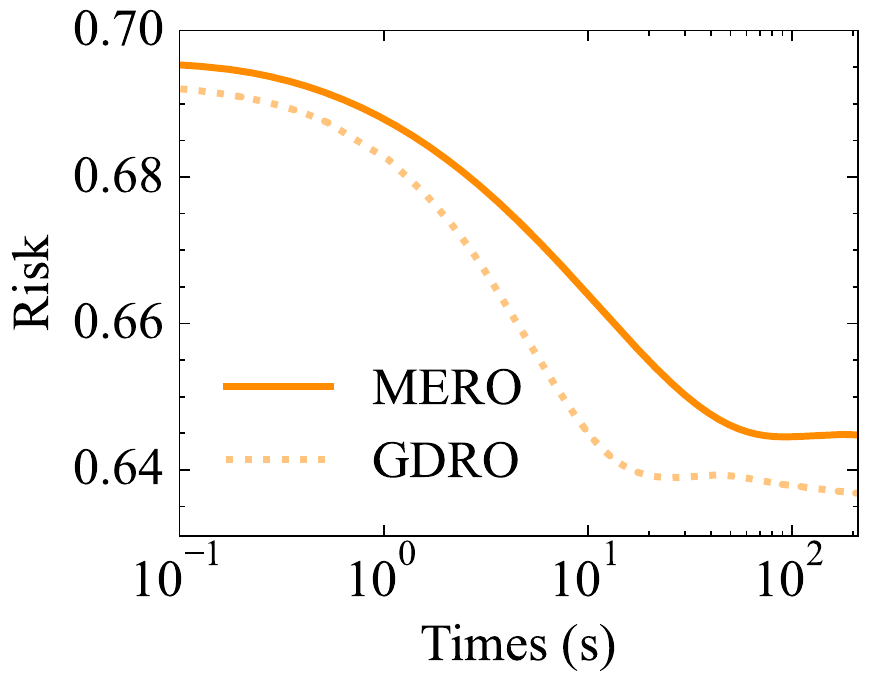}}\\%
  \caption{Individual risk versus the running time on the synthetic dataset.}
  \label{fig:5} \vspace{1ex}%
\subfigure[Risk on $\P_1$]{
    \label{fig:6:a} 
    \includegraphics[width=0.33\textwidth]{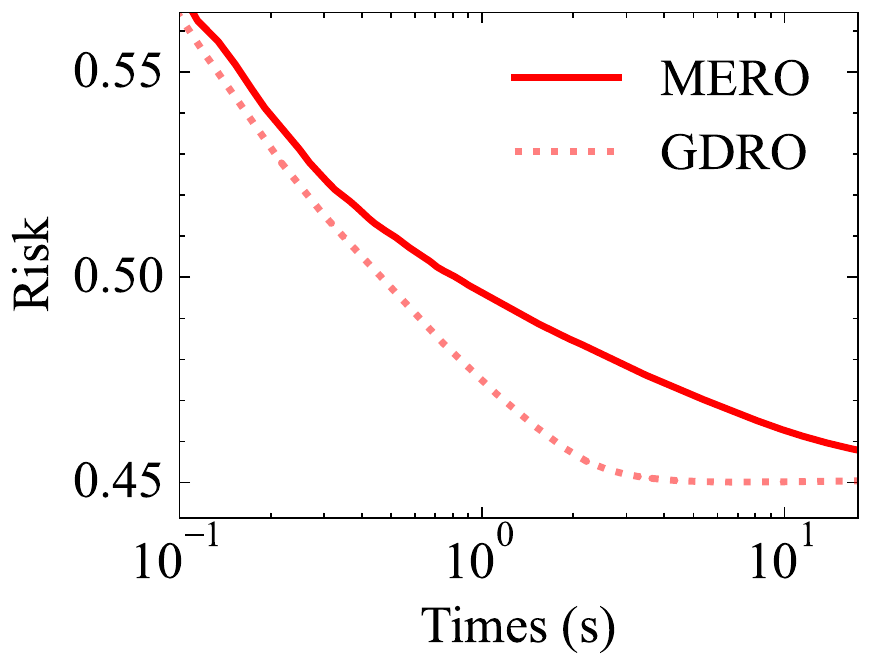}}%
\subfigure[Risk on $\P_{2}$]{
    \label{fig:6:b} 
    \includegraphics[width=0.33\textwidth]{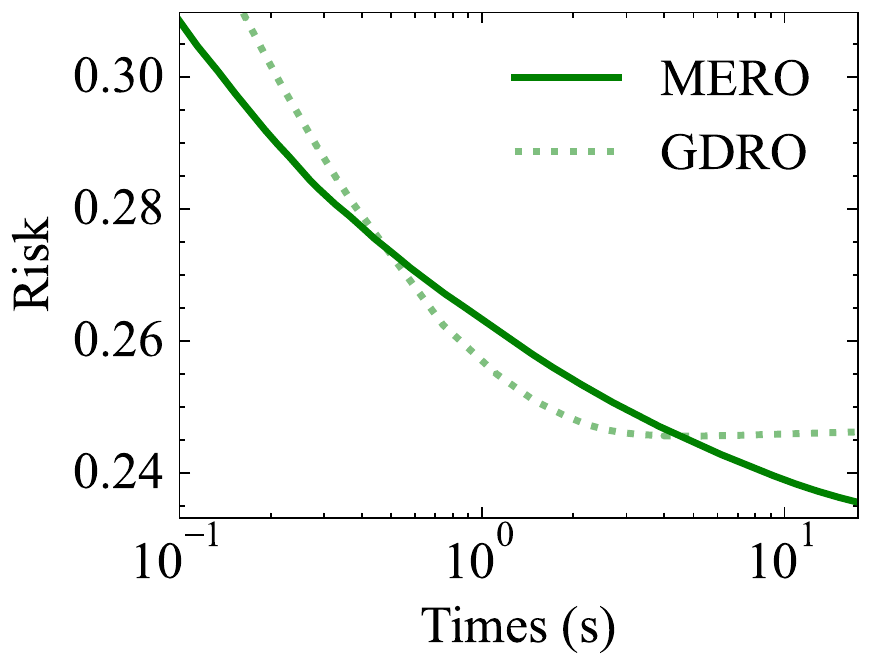}}%
    \subfigure[Risk on $\P_{3}$]{
    \label{fig:6:c} 
    \includegraphics[width=0.33\textwidth]{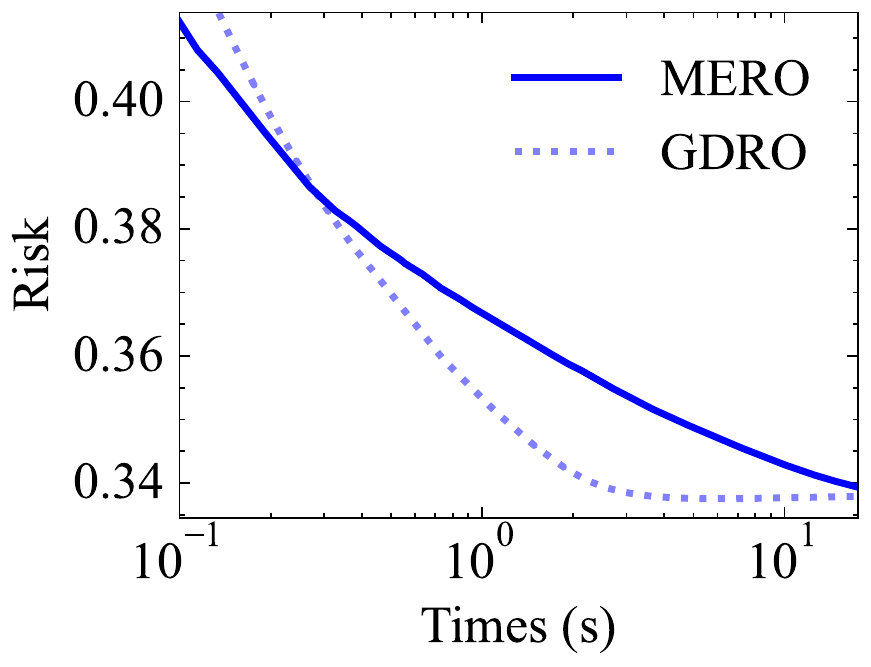}}\\%
        \subfigure[Risk on $\P_{4}$]{
    \label{fig:6:d} 
    \includegraphics[width=0.33\textwidth]{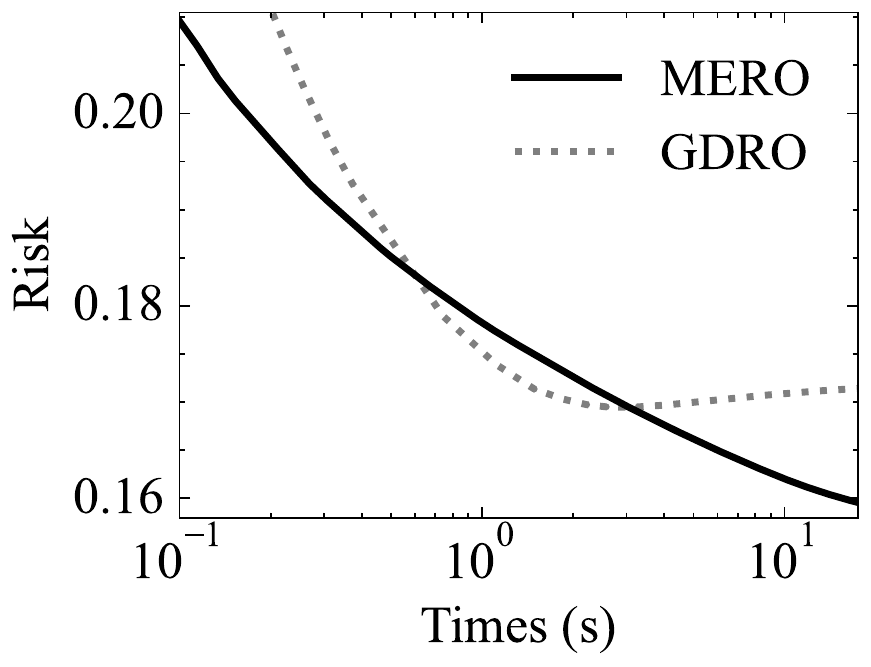}}%
    \subfigure[Risk on $\P_{5}$]{
    \label{fig:6:e} 
    \includegraphics[width=0.33\textwidth]{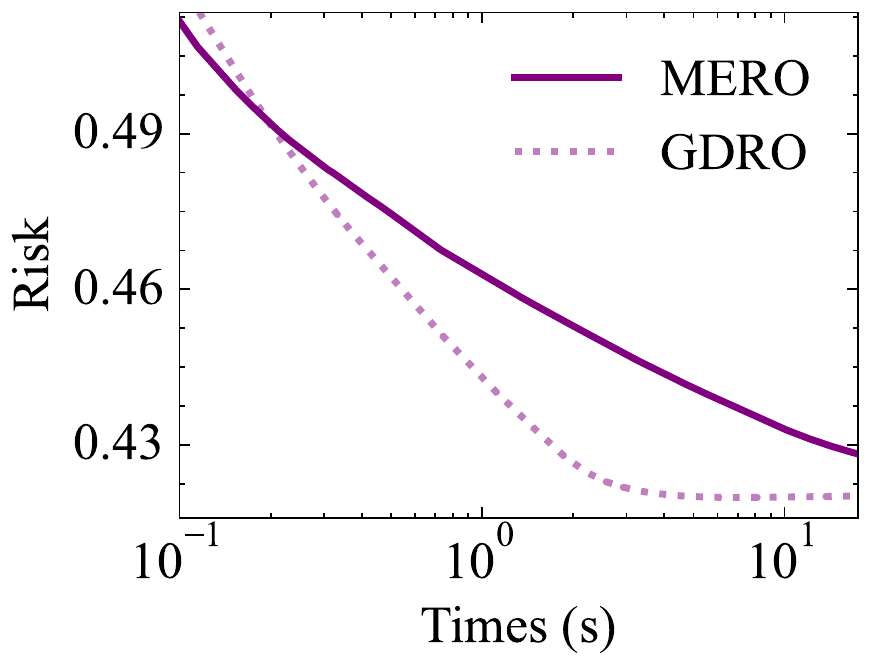}}%
    \subfigure[Risk on $\P_{6}$]{
    \label{fig:6:f} 
    \includegraphics[width=0.33\textwidth]{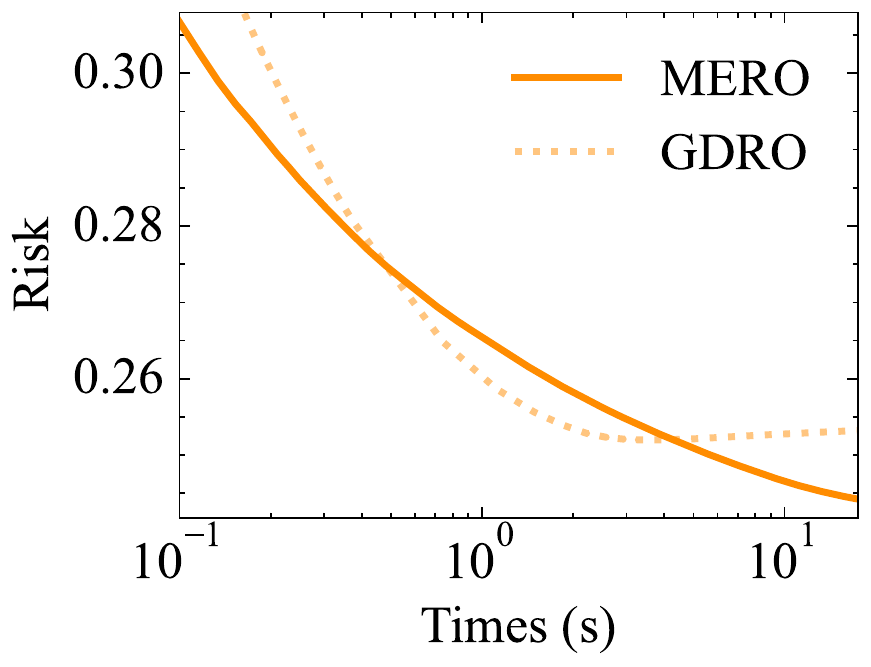}}\\%
  \caption{Individual risk versus the running time on the Adult dataset.}
  \label{fig:6}
\end{center}
\end{figure}

On the other hand, the stochastic algorithm for GDRO is more efficient than our Algorithm~\ref{alg:1}, as it does not require estimating the minimal risk for each distribution. To examine the difference, we plot the risk relative to the running time in Fig.~\ref{fig:5} and Fig.~\ref{fig:6}. We observe that the risk decreases more rapidly for GDRO compared to MERO, primarily because GDRO processes a greater number of samples in the same amount of time. Nevertheless, the ultimate performance of GDRO and MERO is consistent with previous findings presented in Fig.~\ref{fig:3} and Fig.~\ref{fig:4}, with GDRO generally performing better on harder distributions. Additionally, an upward trend in the final curve of GDRO is noted.  This is because Algorithm 1 of \citet{SA:GDRO}, like MS-MERO, is not anytime. As a result, its fixed step size becomes suboptimal as time progresses.

\afterpage{\clearpage}%

\subsection{Experiments on Imbalanced Data}
For experiments involving imbalanced data, we fix the number of samples per distribution, with the quantity varying among distributions. Mirroring the setup of  \citet{SA:GDRO}, we designate the sample size for the synthetic dataset as $n_i= 5000 \times (7-i)$, generating each sample as before. Pertaining to the Adult dataset, we randomly extract $364$ samples from each group, reserving them for subsequent risk assessment. The remaining number of samples across the $6$ groups is \{$26656$, $11519$, $1780$, $1720$, $999$, $364$\}. Furthermore, each sample within the groups is processed only \textit{once} to simulate the imbalanced scenario.  In this way,  $\P_i$ corresponds to the (unknown) underlying distribution from which the samples in the $i$-th group are drawn.

\begin{figure}[t]
\begin{minipage}[c]{0.5\textwidth} 
    \centering 
    \includegraphics[width=0.9\textwidth]{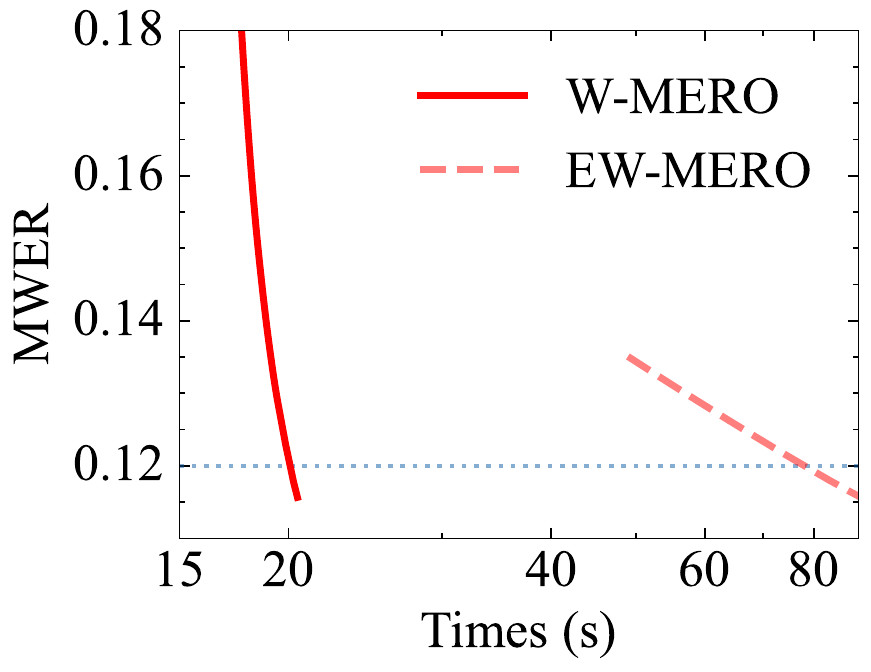} 
    \caption{The MWER versus the running time on the \emph{imbalanced} synthetic dataset.}   \label{fig:7}
  \end{minipage}%
  \begin{minipage}[c]{0.5\textwidth} 
      \tabcaption{Running times of W-MERO and EW-MERO on the \emph{imbalanced} synthetic dataset.}\label{tab:2}
    \centering
		\begin{tabular}{c|c|c}
			\hline
			  Algorithm &  MWER Value & Times (s) \\
			\hline
			  W-MERO & 0.12 & 20.3\\ \hline
		  EW-MERO & 0.12 & 78.5 \\			\hline
		\end{tabular}
\vspace{20ex}
  \end{minipage} 
\end{figure}

Note that the optimization procedure of \citet{Regret:RML:DS} is also applicable for minimizing the empirical counterpart of weighted MERO. Therefore, we initially conduct a comparison between our two-stage approach in Section~\ref{sec:two:stage}, and their method. For ease of reference, we label our algorithm and theirs as W-MERO and EW-MERO, respectively, to emphasize that the former is designed for weighted MERO, while the latter focuses on the empirical variant. In Fig.~\ref{fig:7}, we illustrate how the maximal weighted excess risk (MWER), denoted as $\max_{i\in[m]}\{p_i[R_i(\w)-R_i^*]\}$, changes with respect to the running time. Acknowledging that both algorithms have an initialization phase, so their curves do not start from zero. Consistent with previous experiments in Fig.~\ref{fig:1}, our W-MERO converges more rapidly than EW-MERO.\footnote{In this experiment, the imbalanced Adult dataset is omitted. This is because the limited number of samples prevents an accurate estimation of the value of $R_i^*$.}  For a detailed view, the precise times required for W-MERO and EW-MERO to reach a certain MWER value are listed in Table~\ref{tab:2}.

\begin{figure}[t] 
    \begin{center}
\subfigure[Risk on $\P_1$]{
    \label{fig:8:a} 
    \includegraphics[width=0.33\textwidth]{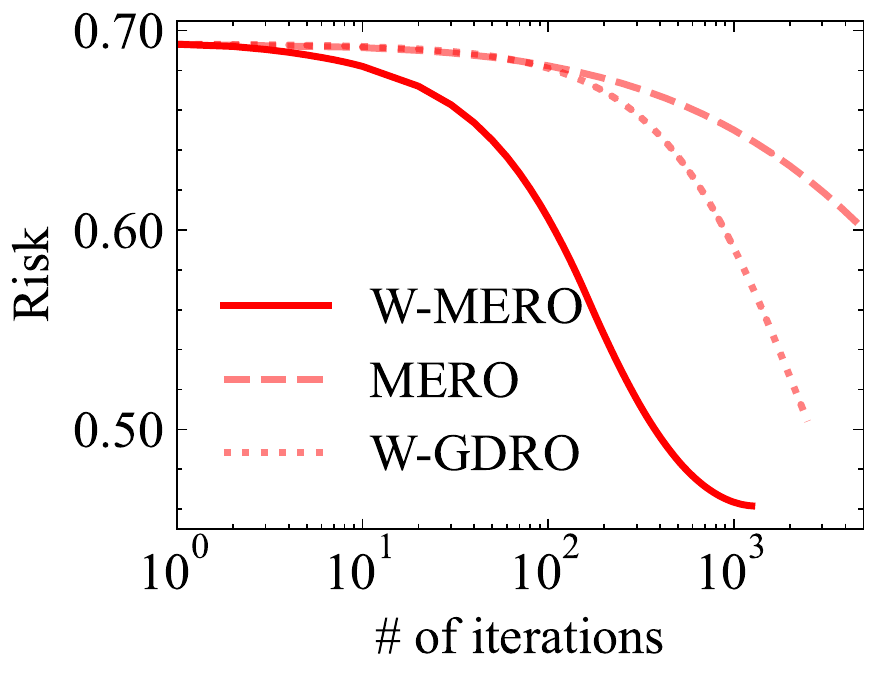}}%
\subfigure[Risk on $\P_{2}$]{
    \label{fig:8:b} 
    \includegraphics[width=0.33\textwidth]{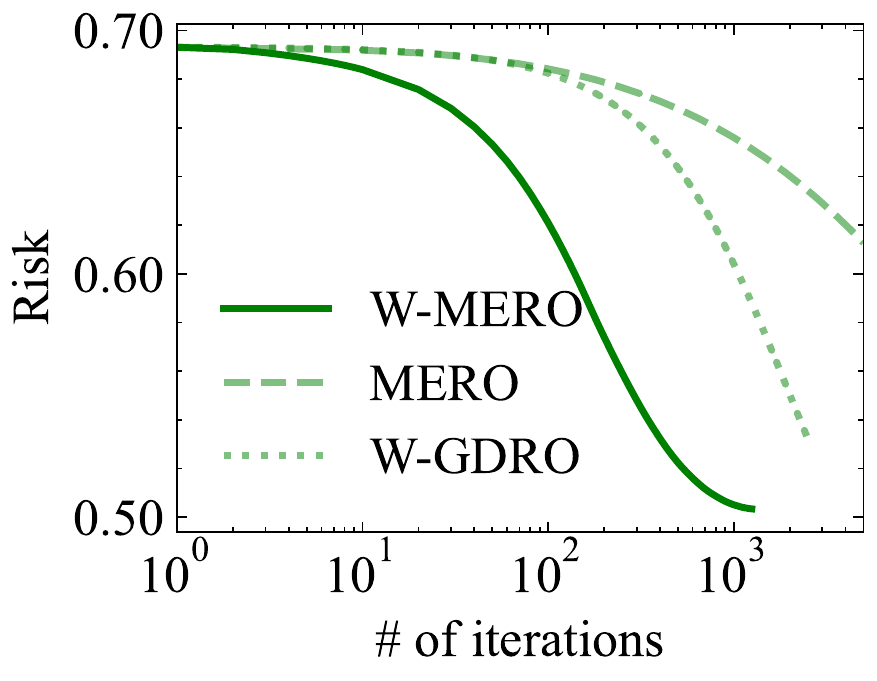}}%
    \subfigure[Risk on $\P_{3}$]{
    \label{fig:8:c} 
    \includegraphics[width=0.33\textwidth]{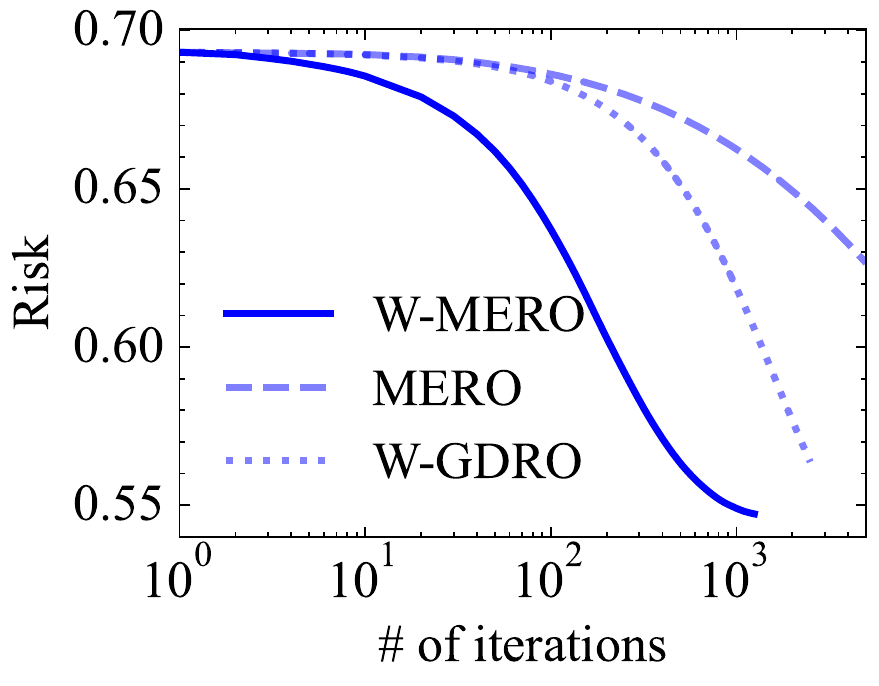}}\\%
        \subfigure[Risk on $\P_{4}$]{
    \label{fig:8:d} 
    \includegraphics[width=0.33\textwidth]{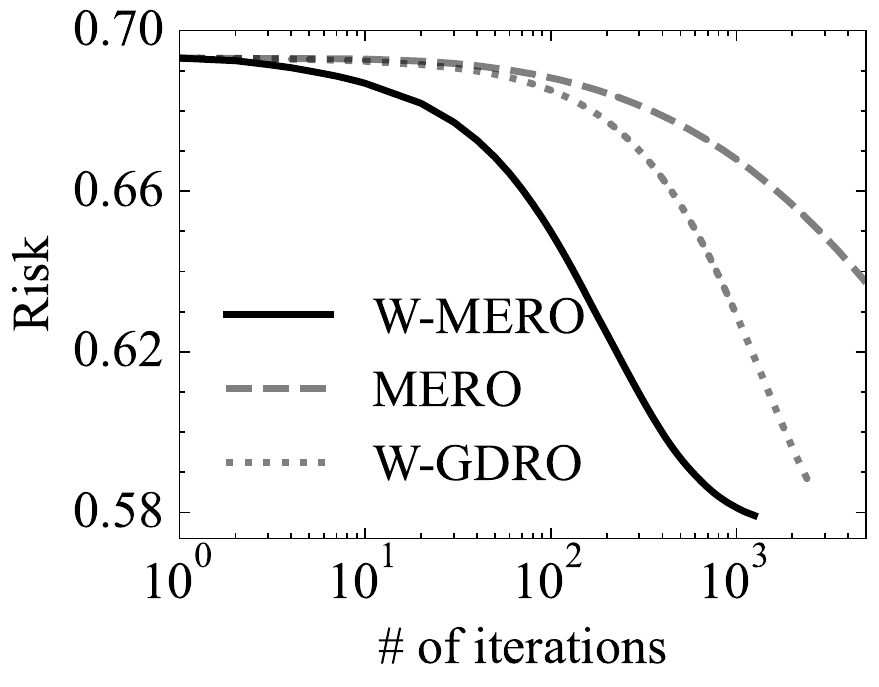}}%
    \subfigure[Risk on $\P_{5}$]{
    \label{fig:8:e} 
    \includegraphics[width=0.33\textwidth]{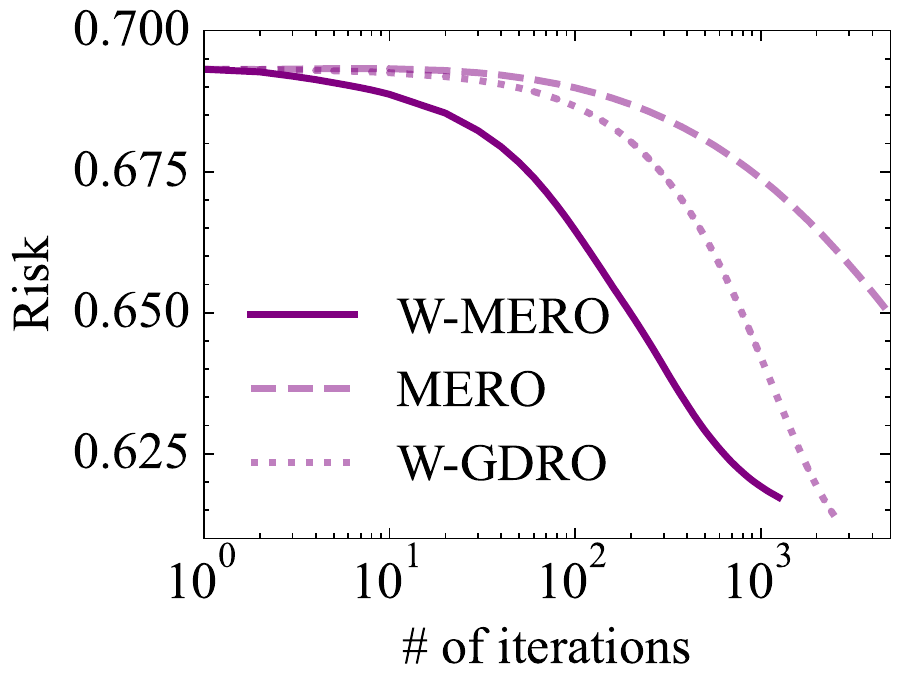}}%
    \subfigure[Risk on $\P_{6}$]{
    \label{fig:8:f} 
    \includegraphics[width=0.33\textwidth]{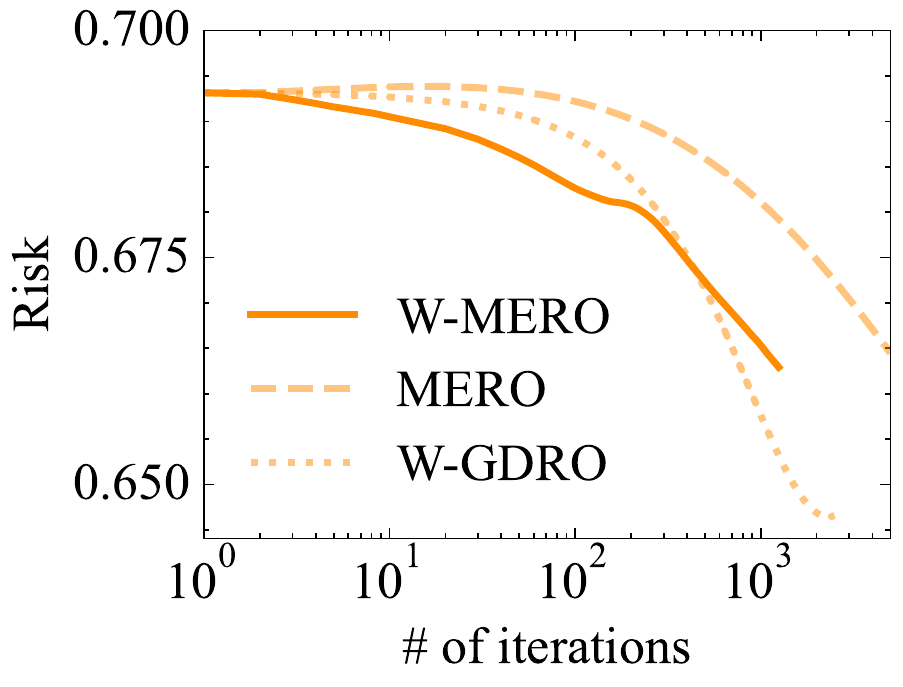}}\\%
  \caption{Individual risk versus the number of iterations on the \emph{imbalanced} synthetic dataset.}
  \label{fig:8} \vspace{1ex}%
\subfigure[Risk on $\P_1$]{
    \label{fig:9:a} 
    \includegraphics[width=0.33\textwidth]{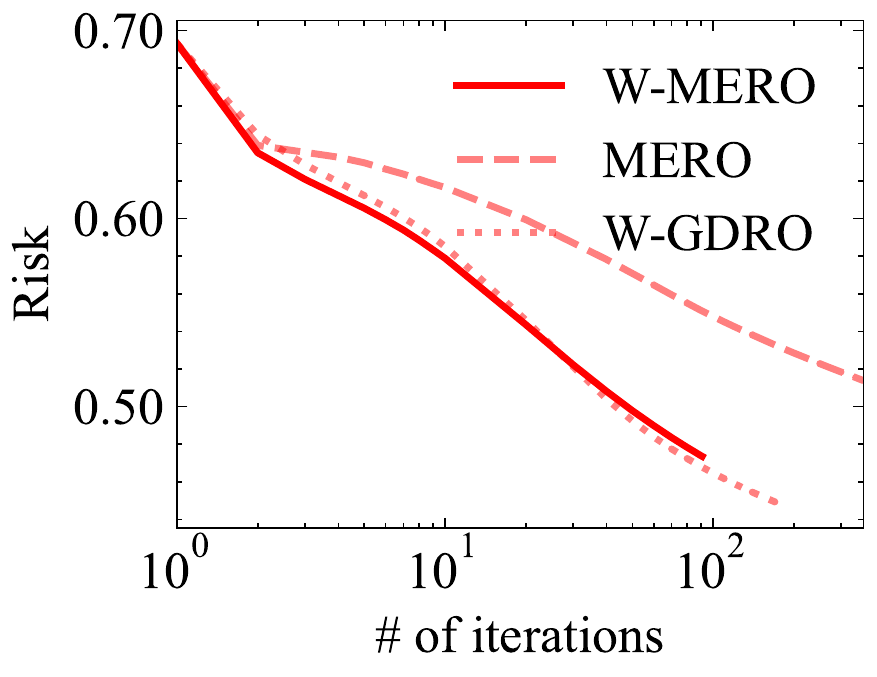}}%
\subfigure[Risk on $\P_{2}$]{
    \label{fig:9:b} 
    \includegraphics[width=0.33\textwidth]{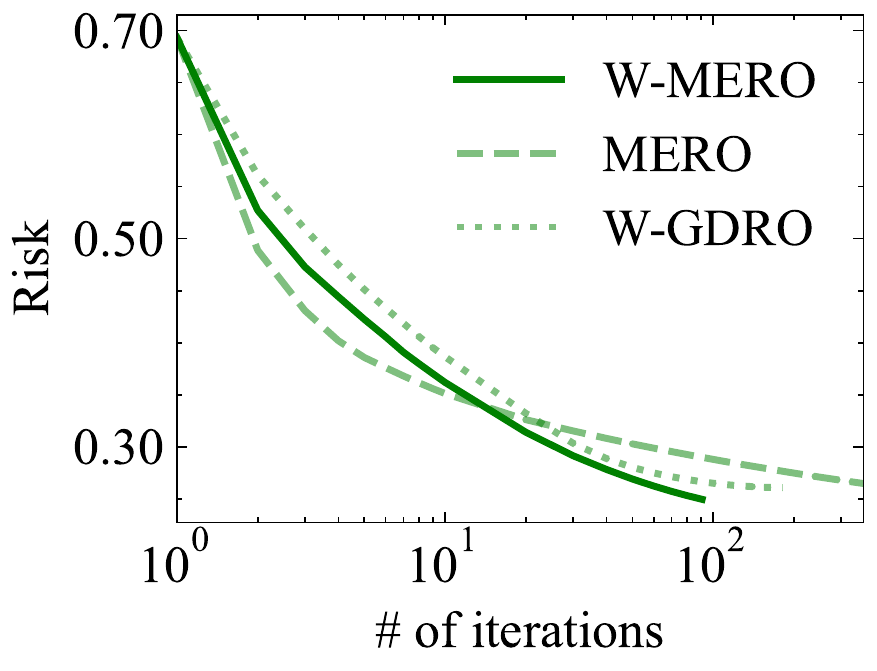}}%
    \subfigure[Risk on $\P_{3}$]{
    \label{fig:9:c} 
    \includegraphics[width=0.33\textwidth]{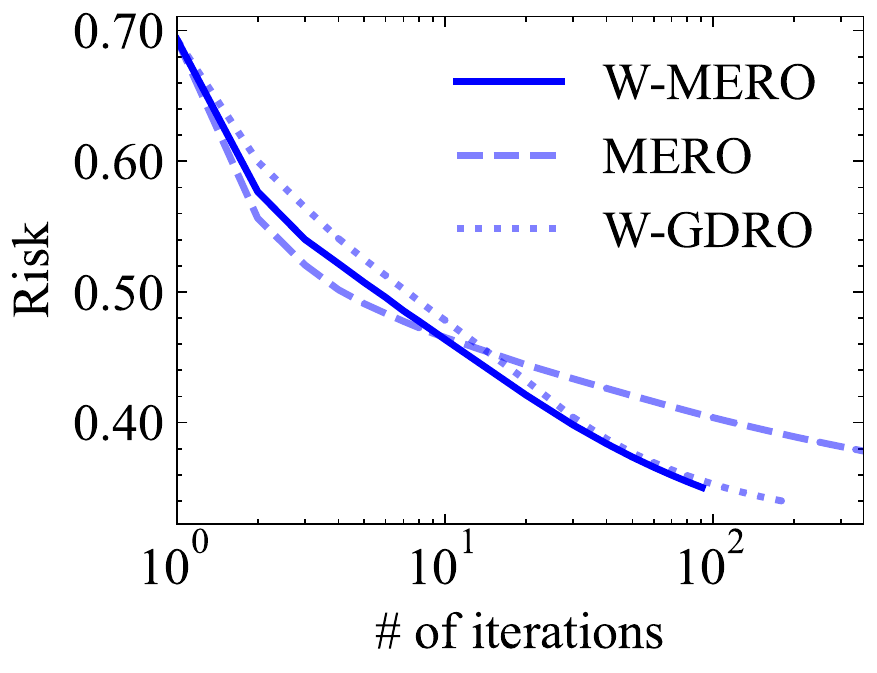}}\\%
        \subfigure[Risk on $\P_{4}$]{
    \label{fig:9:d} 
    \includegraphics[width=0.33\textwidth]{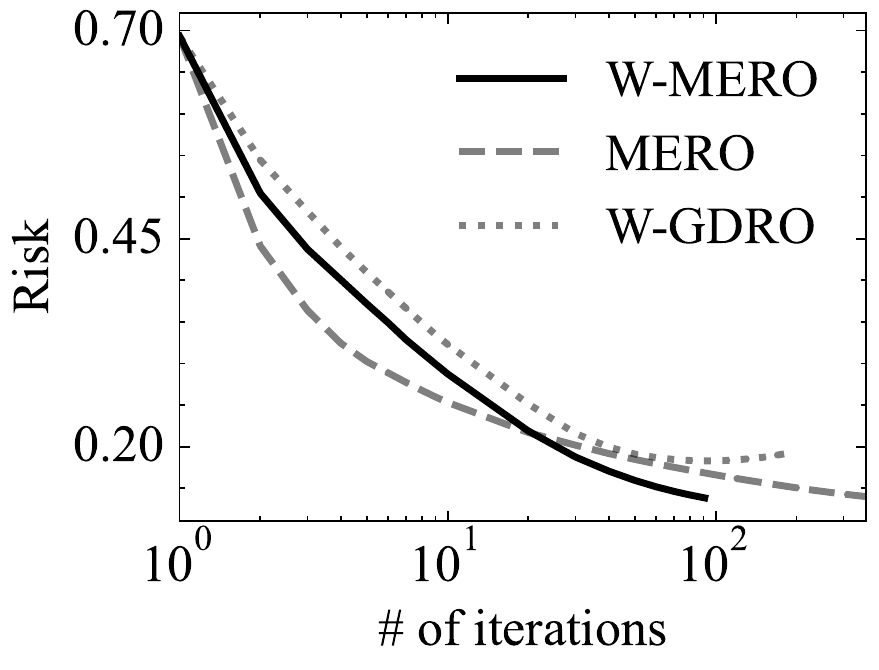}}%
    \subfigure[Risk on $\P_{5}$]{
    \label{fig:9:e} 
    \includegraphics[width=0.33\textwidth]{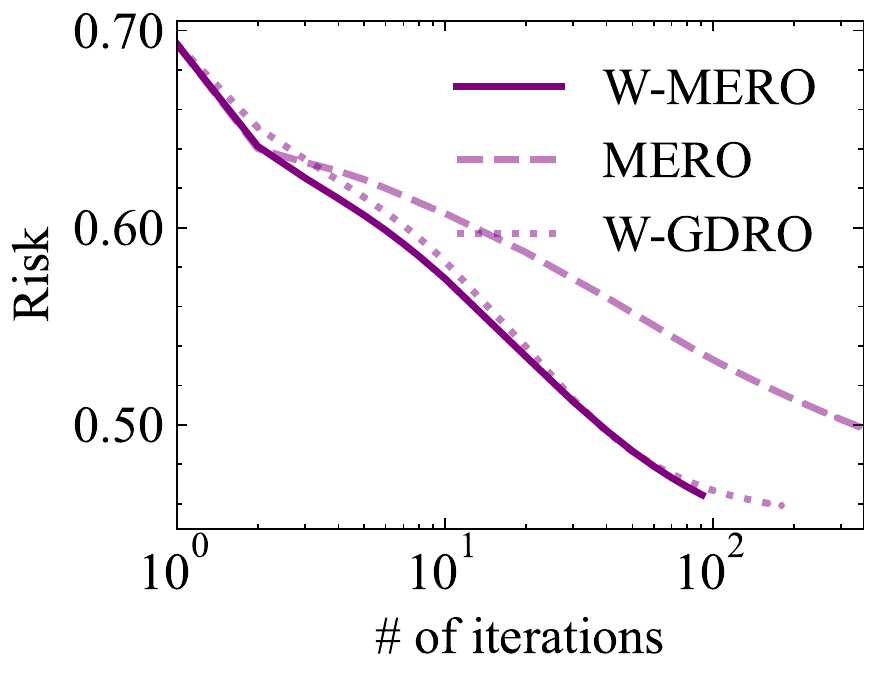}}%
    \subfigure[Risk on $\P_{6}$]{
    \label{fig:9:f} 
    \includegraphics[width=0.33\textwidth]{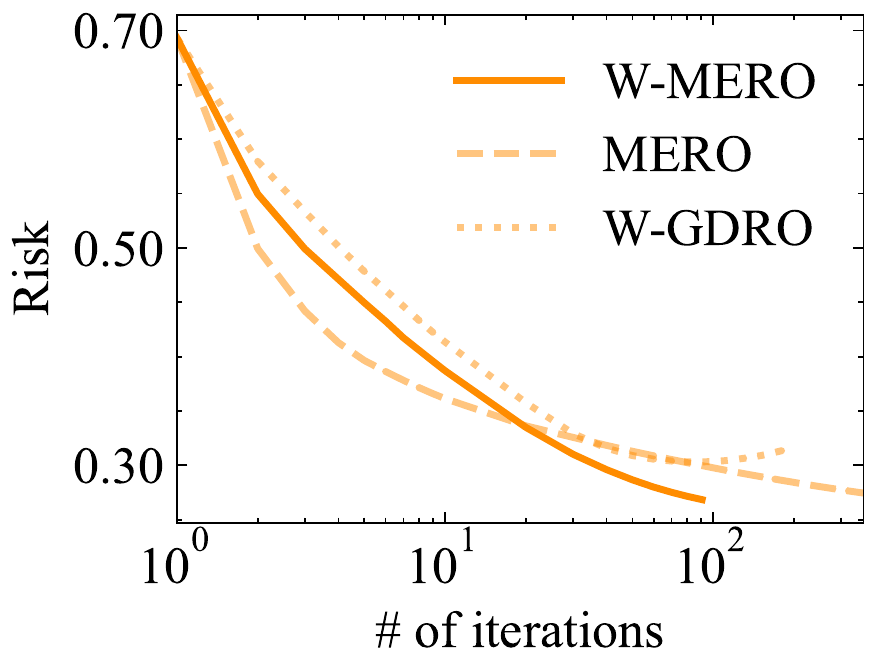}}\\%
  \caption{Individual risk versus the number of iterations on the \emph{imbalanced} Adult dataset.}
  \label{fig:9}
\end{center}
\end{figure}

Next, we examine the effectiveness of our W-MERO in handling imbalanced scenarios. To this end, we compare it with the original MERO---running Algorithm~\ref{alg:1} for $n_m$ iterations. Recall that W-MERO reduces the excess risk of the $i$-th distribution at an $O((\log m)/\sqrt{n_i})$ rate, and MERO attains an $\O(\sqrt{(\log m)/n_m})$ rate for all distributions. Additionally, we include the result of weighted GDRO (W-GDRO) \citep[Algorithm 4]{SA:GDRO} to reiterate the distinction between risk minimization and excess risk minimization. Experimental results on the synthetic dataset are provided in Fig.~\ref{fig:8}, and consistent with our theoretical expectations. Specifically, the final risk of W-MERO closely approaches that of MERO on distribution $\P_6$, which holds the smallest sample size, and W-MERO surpasses MERO on all other distributions. Moreover, the larger the number of samples is, the more pronounced the advantage of W-MERO becomes. In line with the experiments in Fig.~\ref{fig:3}, W-GDRO performs well on the last two distributions, characterized by their significantly high noise levels.  We present the outcomes on the Adult dataset in Fig.~\ref{fig:9}, and observe analogous patterns.

\afterpage{\clearpage}

\section{Conclusion and Future Work}
This paper aims to develop efficient stochastic approximation approaches for MERO. First, we design a  multi-stage stochastic algorithm, which attains a (nearly) optimal convergence rate of $O(\sqrt{(\log m)/T} )$ for a fixed number of iterations $T$. Then, we propose an anytime stochastic method, which reduces the error at an $\O(\sqrt{(\log m)/t})$ rate at every iteration $t$. Next, we delve into the setting where different distributions possess varying sample budgets, and develop a two-stage stochastic procedure that is endowed  with distribution-dependent convergence rates. Finally, we substantiate the efficiency and effectiveness of our methods through  experimental validation.

For special losses such as the squared loss, \citet[\S 4]{Regret:RML:DS} have established faster rates by minimizing the empirical MERO. In the future, we plan to investigate whether it is possible to attain similar improvements by optimizing MERO directly. Another area of future work involves applying stochastic approximation to empirical MERO to reduce the computational cost.  By leveraging the finite-sum structure \citep{NIPS:13:Mixed,NIPS13:ASGD,SVRG:Nonconvex}, we  have recently made progress in this direction \citep{VR:EMERO}.

\bibliography{E:/MyPaper/ref}

\appendix

\section{Supporting Lemmas}
\subsection{Proof of Lemma~\ref{lem:diff:opt}}
From the definition of $\varphi(\cdot,\cdot)$ and $\widehat{\varphi}(\cdot,\cdot)$, for any $\q$, we have
\begin{equation} \label{eqn:sam:w}
\argmin_{\w\in \W}  \varphi(\w,\q) = \argmin_{\w\in \W}\widehat{\varphi}(\w,\q)
\end{equation}
and for any $(\w, \q)$
\begin{equation}\label{eqn:diff:function}
\big|\varphi(\w,\q) - \widehat{\varphi}(\w,\q) \big|= \left|\sum_{i=1}^m q_i p_i \big[ R_i(\wb^{(i)})- R_i^*\big]\right| \leq \max_{i \in [m]} \left\{p_i \big[ R_i(\wb^{(i)})- R_i^*\big]\right\}.
\end{equation}

Let $\wh=\argmin_{\w\in \W}  \varphi(\w,\qb)$  and $\qh=\argmax_{\q\in \Delta_m} \varphi(\wb,\q)$. Then, we have
\[
\begin{split}
\epsilon_{\varphi}(\wb, \qb) = &\max_{\q\in \Delta_m}  \varphi(\wb,\q)- \min_{\w\in \W}  \varphi(\w,\qb)=\varphi(\wb,\qh) - \varphi(\wh,\qb)\\
\overset{\text{(\ref{eqn:diff:function})}}{\leq} & \widehat{\varphi}(\wb,\qh) - \widehat{\varphi}(\wh,\qb) + 2\max_{i \in [m]}\left\{ p_i \big[ R_i(\wb^{(i)})- R_i^*\big] \right\}\\
\leq &\max_{\q\in \Delta_m}  \widehat{\varphi}(\wb,\q) - \widehat{\varphi}(\wh,\qb) + 2\max_{i \in [m]} \left\{p_i \big[ R_i(\wb^{(i)})- R_i^*\big]\right\} \\
\overset{\text{(\ref{eqn:sam:w})}}{=} & \max_{\q\in \Delta_m}  \widehat{\varphi}(\wb,\q)- \min_{\w\in \W}  \widehat{\varphi}(\w,\qb)+ 2\max_{i \in [m]}\left\{ p_i \big[ R_i(\wb^{(i)})- R_i^*\big]\right\}\\
=& \epsilon_{\widehat{\varphi}}(\wb, \qb)+2\max_{i \in [m]} \left\{p_i \big[ R_i(\wb^{(i)})- R_i^*\big]\right\}.
\end{split}
\]
\subsection{Proof of Lemma~\ref{lem:1}}
We create a virtual sequence by performs SMD  with  $\E_{t-1} \left[ \g(\w_t,\q_t) \right]-\g(\w_t,\q_t)$ as the gradient:
\begin{equation} \label{eqn:virtual:smd}
\y_{t+1}= \argmin_{\x \in \W \times \Delta_m}  \Big\{ \eta_j \big \langle  \E_{t-1} [ \g(\w_t,\q_t) ]-\g(\w_t,\q_t) , \x -\y_t \big\rangle + B(\x,\y_t) \Big\}
\end{equation}
where $\y_1=\x_1$. Then, we further decompose the error term as
\begin{equation} \label{lem:1:1}
\begin{split}
&\max_{\x \in \W \times \Delta_m}  \sum_{j=1}^t \eta_j \big \langle  \E_{j-1} \left[ \g(\w_j,\q_j) \right]-\g(\w_j,\q_j) , \x_j - \x \big\rangle \\
\leq & \max_{\x \in \W \times \Delta_m}  \underbrace{\sum_{j=1}^t \eta_j \big \langle  \E_{j-1} \left[ \g(\w_j,\q_j) \right]-\g(\w_j,\q_j) , \y_j - \x \big\rangle}_{:=A}\\
& + \underbrace{\sum_{j=1}^t \eta_j \big \langle  \E_{j-1} \left[ \g(\w_j,\q_j) \right]-\g(\w_j,\q_j) , \x_j - \y_j \big\rangle}_{:=B}. \\
 \end{split}
\end{equation}

To bound term $A$, we repeat the analysis of (\ref{eqn:thm2:1}), and have
\begin{equation} \label{lem:1:2}
\begin{split}
& \eta_j \big \langle  \E_{j-1} \left[ \g(\w_j,\q_j) \right]-\g(\w_j,\q_j) , \y_j - \x \big\rangle \\
\leq & B(\x, \y_j)-B(\x, \y_{j+1})   + \frac{\eta_j^2 }{2} \big\|\E_{j-1} \left[ \g(\w_j,\q_j) \right]-\g(\w_j,\q_j) \big \|_*^2 \\
\leq &  B(\x, \y_j)-B(\x, \y_{j+1})   +   2M^2  \eta_j^2
\end{split}
\end{equation}
where the last step makes use of the following inequality
\begin{equation} \label{lem:1:3}
\begin{split}
&\big\|\E_{j-1} \left[ \g(\w_j,\q_j) \right] -\g(\w_j,\q_j) \big\|_* \leq \big\| \E_{j-1} \left[ \g(\w_j,\q_j) \right]\big\|_* + \|\g(\w_j,\q_j)\|_* \\
\leq &  \E_{j-1}\big[\|   \g(\w_j,\q_j)\|_*  \big] + \|\g(\w_j,\q_j)\|_* \overset{\text{(\ref{eqn:gradient:bound})}}{\leq} 2 M.
\end{split}
\end{equation}
Summing (\ref{lem:1:2}) over $j=1,\ldots,t$, we have
\begin{equation} \label{lem:1:4}
\begin{split}
 A =& \sum_{j=1}^t  \eta_j \big \langle  \E_{j-1} \left[ \g(\w_j,\q_j) \right]-\g(\w_j,\q_j) , \y_j - \x \big\rangle \\
\leq &B(\x, \y_1)   + 2M^2  \sum_{j=1}^t\eta_j^2 \overset{\text{(\ref{eqn:domain:merge})}}{\leq}  1+ 2M^2  \sum_{j=1}^t\eta_j^2 .
\end{split}
\end{equation}

To bound term $B$ in (\ref{lem:1:1}), we define
\[
\delta_j = \eta_j \big \langle  \E_{j-1} \left[ \g(\w_j,\q_j) \right]-\g(\w_j,\q_j) , \x_j - \y_j \big\rangle.
\]
As  $\x_j$ and $\y_j$ are independent of the random samples $\z_j^{(1)},\ldots,\z_j^{(m)}$  used to construct $\g(\w_j,\q_j)$ in (\ref{eqn:stochastic:gradient}), $\delta_1,\ldots,\delta_t$ forms a martingale difference sequence. Consequently, we have
\begin{equation} \label{lem:1:5}
\E[ B] = \E\left[ \sum_{j=1}^t \delta_j \right] = 0.
\end{equation}
Taking expectations over both sides of (\ref{lem:1:1}), we have
\[
\begin{split}
& \E \left[ \max_{\x \in \W \times \Delta_m}  \sum_{j=1}^t \eta_j \big \langle  \E_{j-1} \left[ \g(\w_j,\q_j) \right]-\g(\w_j,\q_j) , \x_j - \x \big\rangle \right] \\
\leq & \E \left[ \max_{\x \in \W \times \Delta_m} A \right] + \E[ B]\overset{\text{(\ref{lem:1:4}),(\ref{lem:1:5})}}{\leq}  1+ 2M^2  \sum_{j=1}^t\eta_j^2
\end{split}
\]
which proves (\ref{lem:1:exp}).

To establish a high probability bound, we follow the analysis in Section~\ref{sec:high:pro} and utilize Lemma~\ref{lem:azuma} to bound $B$. To this end, we first show that $|\delta_j|$ is bounded:
\[
\begin{split}
 |\delta_j| = &\Big | \eta_j \big \langle  \E_{j-1} \left[ \g(\w_j,\q_j) \right]-\g(\w_j,\q_j) , \x_j - \y_j \big\rangle \Big| \\
\leq & \eta_j \big\|\E_{j-1} \left[ \g(\w_j,\q_j) \right]-\g(\w_j,\q_j) \big\|_* \| \x_j - \y_j\| \\
\overset{\text{(\ref{lem:1:3})}}{\leq} & 2M \eta_j  \big(\| \x_j -\x_1\| + \|\x_1- \y_j\|\big) \\
\leq &  2 M \eta_j \left(\sqrt{2 B(\x_j,\x_1)} + \sqrt{2 B(\y_j,\x_1)}\right)  \overset{\text{(\ref{eqn:domain:merge})}}{\leq}  4 \sqrt{2} M\eta_j .
\end{split}
\]
From Lemma~\ref{lem:azuma} and the union bound, with probability at least $1-\delta$, we have
\begin{equation} \label{lem:1:6}
B=\sum_{j=1}^t \delta_j \leq 8 M \sqrt{\sum_{j=1}^t \eta_j^2 \left(\ln \frac{2 t^2}{\delta}\right)} , \quad \forall  t \in \zn.
\end{equation}
We obtain (\ref{lem:1:high}) by substituting (\ref{lem:1:4}) and (\ref{lem:1:6}) into (\ref{lem:1:1}).

\subsection{Proof of Lemma~\ref{lem:2}}
First, we have
\[
\begin{split}
 & \big \langle  F(\w_j,\q_j) -\E_{t-1} \left[ \g(\w_j,\q_j) \right], \x_j - \x \big \rangle  \\
 \overset{\text{(\ref{eqn:gradient:bias})}}{=}&  \left \langle  \Big[0 , -\big[R_1(\wb_j^{(1)})-R_1^*, \ldots, R_m(\wb_j^{(m)})-R_m^*\big]^\top \Big], \x_j - \x  \right\rangle  \\
=& -\Big \langle \big[R_1(\wb_j^{(1)})-R_1^*, \ldots, R_m(\wb_j^{(m)})-R_m^*\big] , \q_j - \q\Big \rangle \\
\leq & \Big\|\big[R_1(\wb_j^{(1)})-R_1^*, \ldots, R_m(\wb_j^{(m)})-R_m^*\big]\Big\|_\infty \|\q_j - \q\|_1 \\
\leq & 2 \Big\|\big[R_1(\wb_j^{(1)})-R_1^*, \ldots, R_m(\wb_j^{(m)})-R_m^*\big]\Big\|_\infty
\end{split}
\]
which implies
\begin{equation} \label{eqn:lem2:1}
\begin{split}
&\max_{\x \in \W \times \Delta_m}  \sum_{j=1}^t \eta_j \big \langle  F(\w_j,\q_j) -\E_{j-1} \left[ \g(\w_j,\q_j) \right], \x_j - \x \big \rangle \\
\leq &\sum_{j=1}^t 2 \eta_j  \Big\|\big[R_1(\wb_j^{(1)})-R_1^*, \ldots, R_m(\wb_j^{(m)})-R_m^*\big]\Big\|_\infty.
 \end{split}
\end{equation}

From (\ref{eqn:exp:high}) in Theorem~\ref{thm:1}, we know that with probability at least $1-\delta$,
\begin{equation} \label{eqn:lem2:2}
\Big\|\big[R_1(\wb_j^{(1)})-R_1^*, \ldots, R_m(\wb_j^{(m)})-R_m^*\big]\Big\|_\infty \leq \frac{DG\big[3+\ln j + 16\sqrt{(1+ \ln j) \ln (2 m j^2/\delta) } \big]}{4(\sqrt{j+1} - 1)}
\end{equation}
for all $j \in \zn$. We obtain (\ref{lem:2:high}) by substituting (\ref{eqn:lem2:2}) into (\ref{eqn:lem2:1}).

The proof of the expectation bound is more involved. Taking expectations over (\ref{eqn:lem2:1}), we have
\begin{equation} \label{eqn:lem2:3}
\begin{split}
&\E \left[\max_{\x \in \W \times \Delta_m}  \sum_{j=1}^t \eta_j \big \langle  F(\w_j,\q_j) -\E_{j-1} \left[ \g(\w_j,\q_j) \right], \x_j - \x \big \rangle \right] \\
\leq &\sum_{j=1}^t 2 \eta_j  \E\left[\Big\|\big[R_1(\wb_j^{(1)})-R_1^*, \ldots, R_m(\wb_j^{(m)})-R_m^*\big]\Big\|_\infty \right].
 \end{split}
\end{equation}
Then, one may attempt to make use of the expectation bound  (\ref{eqn:exp:risk}) in Theorem~\ref{thm:1}. However, due to the presence of the infinity norm in (\ref{eqn:lem2:3}), it is difficult to obtain a tight upper bound. As an alternative, we will exploit the fact that $R_i(\wb_t^{(i)})-R_i^*$ is sub-Gaussian \citep{HDP:Vershynin}, for all $i \in [m]$,  $t \in \zn$.

Recall the martingale difference sequence $\delta_1^{(i)},\ldots,\delta_t^{(i)}$ in (\ref{eqn:delta:mart}). From (\ref{eqn:bound:delta}) and Lemma~\ref{lem:azuma}, we have
\begin{equation} \label{eqn:lem2:4}
\Pr\left[ \sum_{j=1}^t \delta_j^{(i)} > x \right]  \leq \exp\left( - \frac{x^2}{64 D^2G^2 \sum_{j=1}^t(\eta_j^{(i)})^2 } \right).
\end{equation}
From  (\ref{eqn:mr:risk:3}), we have
\[
\begin{split}
&\Pr\left[ \left(R_i(\wb_t^{(i)}) - R_i(\w_*^{(i)}) \right) \left( \sum_{j=1}^t \eta_j^{(i)} \right)- \left( D^2 + \frac{ G^2 }{2} \sum_{j=1}^t  (\eta_j^{(i)})^2 \right) > x \right]  \\
\leq & \Pr\left[ \sum_{j=1}^t \delta_j^{(i)} > x \right] \overset{\text{(\ref{eqn:lem2:4})}}{\leq}  \exp\left( - \frac{x^2}{64 D^2G^2 \sum_{j=1}^t(\eta_j^{(i)})^2 } \right)
\end{split}
\]
which implies
\[
\begin{split}
&\Pr\left[ R_i(\wb_t^{(i)}) - R_i(\w_*^{(i)})  - \frac{1}{\sum_{j=1}^t \eta_j^{(i)} }\left( D^2 + \frac{ G^2 }{2} \sum_{j=1}^t  (\eta_j^{(i)})^2 \right) >  x \right] \\
 \leq & \exp\left( - \frac{x^2 \left(\sum_{j=1}^t \eta_j^{(i)}\right)^2}{64 D^2G^2 \sum_{j=1}^t(\eta_j^{(i)})^2 } \right).
\end{split}
\]

Since $\eta_j^{(1)}=\cdots=\eta_j^{(m)}$, we can invoke the following lemma to bound the expectation of  $\max_{i\in [m]} \{R_i(\wb_t^{(i)}) - R_i(\w_*^{(i)}) \}$.
\begin{lem}\label{lem:max:subG}
Suppose there are $m$ non-negative random variables $X_i$ such that
\[
\Pr\left[ X_i \geq \mu +t \right] \leq  \exp\left(-\frac{t^2}{\sigma^2}\right)
\]
for all $i \in [m]$.  Then, we have
\[
\E\left[\max_{i\in[m]} X_i\right] \leq \mu+\sigma\sqrt{2 \ln m} + \frac{\sigma \sqrt{2 \pi}}{2}.
\]
\end{lem}
From the above lemma, we have
\begin{equation} \label{eqn:lem2:5}
\begin{split}
& \E \left[\max_{i\in[m]}  \left\{R_i(\wb_t^{(i)}) - R_i(\w_*^{(i)}) \right\}\right]\\
 \leq & \frac{1}{\sum_{j=1}^t \eta_j^{(i)} }\left( D^2 + \frac{ G^2 }{2} \sum_{j=1}^t  (\eta_j^{(i)})^2 \right)+ \left(\sqrt{2 \ln m} + \frac{ \sqrt{2 \pi}}{2}\right) \sqrt{\frac{64 D^2G^2 \sum_{j=1}^t(\eta_j^{(i)})^2 }{\left(\sum_{j=1}^t \eta_j^{(i)}\right)^2}}\\
 =& \frac{2 D^2 + G^2 \sum_{j=1}^t  (\eta_j^{(i)})^2 + 16 DG(\sqrt{\ln m} + \sqrt{\pi}/2)  \sqrt{2\sum_{j=1}^t(\eta_j^{(i)})^2 }}{2\sum_{j=1}^t \eta_j^{(i)} } \\
 =& \frac{2DG + DG  \sum_{j=1}^t  \frac{1}{j} + 16 DG(\sqrt{\ln m} + \sqrt{\pi}/2)  \sqrt{2\sum_{j=1}^t \frac{1}{j} }}{2 \sum_{j=1}^t \frac{1}{\sqrt{j}} }\\
\overset{\text{(\ref{eqn:mr:risk:6})}}{\leq} &  \frac{DG\big[3+\ln t + 16 (\sqrt{\ln m} + \sqrt{\pi}/2) \sqrt{2(1+\ln t)}\big]}{4(\sqrt{t+1} - 1)}.
\end{split}
\end{equation}
We obtain (\ref{lem:2:exp}) by substituting (\ref{eqn:lem2:5}) into (\ref{eqn:lem2:3}) and noticing $\sqrt{\pi}/2 \leq 1$.

\subsection{Proof of Lemma~\ref{lem:max:subG}}
First, we have
\[
\Pr\left[\max_{i\in[m]} X_i  \geq \mu +t \right] \leq \sum_{i=1}^m  \Pr\left[ X_i  \geq \mu +t \right] \leq m \exp\left(-\frac{t^2}{\sigma^2}\right).
\]
To simplify the notation, we define $X=\max_{i\in[m]} X_i$. Since $X$ is non-negative, we have
\[
\begin{split}
\E\left[X\right] \leq &\int_{0}^\infty \Pr\left[ X  \geq x\right] d x = \int_{0}^{\mu+\sigma \sqrt{2 \ln m} } \Pr\left[ X  \geq x\right] d x  + \int_{\mu+\sigma\sqrt{2 \ln m}}^\infty \Pr\left[ X  \geq x\right] d x \\
\leq & \mu+\sigma\sqrt{2 \ln m} + \int_{\mu+\sigma\sqrt{2 \ln m}}^\infty m \exp\left(-\frac{(x-\mu)^2}{\sigma^2}\right) d x \\
=& \mu+\sigma\sqrt{2 \ln m} + \int_{\mu+\sigma\sqrt{2 \ln m}}^\infty m \exp\left(-\frac{(x-\mu)^2}{2\sigma^2}\right) \exp\left(-\frac{(x-\mu)^2}{2\sigma^2}\right) d x.
\end{split}
\]

When $x \geq \mu+\sigma \sqrt{2 \ln m}$, we have
\[
m \exp\left(-\frac{(x-\mu)^2}{2\sigma^2} \right) \leq 1.
\]
Thus
\[
\begin{split}
\E\left[X\right] \leq & \mu+\sigma\sqrt{2 \ln m} + \int_{\mu+\sigma\sqrt{2 \ln m}}^\infty \exp\left(-\frac{(x-\mu)^2}{2\sigma^2}\right) d x \\
\leq & \mu+\sigma\sqrt{2 \ln m} + \int_{\mu}^\infty \exp\left(-\frac{(x-\mu)^2}{2\sigma^2}\right) d x = \mu+\sigma\sqrt{2 \ln m} + \frac{\sigma \sqrt{2 \pi}}{2}.
\end{split}
\]

\end{document}